\documentclass{elsarticle}

\usepackage{bm}
\usepackage{amsmath}
\usepackage{amsfonts}
\usepackage{amsthm}
\usepackage{amssymb}
\usepackage{graphicx}
\usepackage{paralist}
\usepackage{algorithm}
\usepackage{pstool}
\usepackage{wrapfig}
\usepackage{subcaption}
\usepackage{algpseudocode}
\usepackage{mathtools}
\usepackage{fullpage}
\usepackage{etoolbox}

\usepackage{tikz}
\usepackage{pgfplots}
\usepackage{pgfplotstable, filecontents}
\pgfplotsset{compat=1.9}
\usetikzlibrary{pgfplots.groupplots}
\usepgfplotslibrary{fillbetween}
\usetikzlibrary{calc,fit,matrix,arrows,automata,positioning,shapes}
\usepackage{pgfplotstable, booktabs}
\usetikzlibrary{patterns,decorations.pathmorphing,decorations.markings,decorations.pathmorphing}
\pgfplotsset{select coords between index/.style 2 args={
    x filter/.code={
        \ifnum\coordindex<#1\fi
        \ifnum\coordindex>#2\fi
    }
}}

\pgfplotsset{compat=newest}

\usetikzlibrary{calc}

\newcommand{\logLogSlopeTriangle}[5]
{

    \pgfplotsextra
    {
        \pgfkeysgetvalue{/pgfplots/xmin}{\xmin}
        \pgfkeysgetvalue{/pgfplots/xmax}{\xmax}
        \pgfkeysgetvalue{/pgfplots/ymin}{\ymin}
        \pgfkeysgetvalue{/pgfplots/ymax}{\ymax}

        \pgfmathsetmacro{\xArel}{#1}
        \pgfmathsetmacro{\yArel}{#3}
        \pgfmathsetmacro{\xBrel}{#1-#2}
        \pgfmathsetmacro{\yBrel}{\yArel}
        \pgfmathsetmacro{\xCrel}{\xArel}

        \pgfmathsetmacro{\lnxB}{\xmin*(1-(#1-#2))+\xmax*(#1-#2)} 
        \pgfmathsetmacro{\lnxA}{\xmin*(1-#1)+\xmax*#1} 
        \pgfmathsetmacro{\lnyA}{\ymin*(1-#3)+\ymax*#3} 
        \pgfmathsetmacro{\lnyC}{\lnyA+#4*(\lnxA-\lnxB)}
        \pgfmathsetmacro{\yCrel}{\lnyC-\ymin)/(\ymax-\ymin)} 

        \coordinate (A) at (rel axis cs:\xArel,\yArel);
        \coordinate (B) at (rel axis cs:\xBrel,\yBrel);
        \coordinate (C) at (rel axis cs:\xCrel,\yCrel);

        \draw[#5]   (A)-- node[pos=0.5,anchor=north,scale=0.8] {1}
                    (B)-- 
                    (C)-- node[pos=0.5,anchor=west,scale=0.8] {#4}
                    cycle;
    }
}

\newcommand{\norm}[1]{\ensuremath{\left\| #1 \right\|}}

\newcommand{\pder}[2]{\ensuremath{\frac{\partial #1}{\partial #2}}} 

\newcommand{\oder}[2]{\ensuremath{\frac{\mathrm{d} #1}{\mathrm{d} #2}}} 

\newcommand{\Acal}{\ensuremath{\mathcal{A}}}

\newcommand{\Dcal}{\ensuremath{\mathcal{D}}}

\newcommand{\Fcal}{\ensuremath{\mathcal{F}}}
\newcommand{\Gcal}{\ensuremath{\mathcal{G}}}

\newcommand{\Ical}{\ensuremath{\mathcal{I}}}

\newcommand{\Lcal}{\ensuremath{\mathcal{L}}}

\newcommand{\Ocal}{\ensuremath{\mathcal{O}}}

\newcommand{\Ucal}{\ensuremath{\mathcal{U}}}

\newcommand{\Aboldcal}{\ensuremath{\boldsymbol{\mathcal{A}}}}

\newcommand{\Cboldcal}{\ensuremath{\boldsymbol{\mathcal{C}}}}

\newcommand{\Pboldcal}{\ensuremath{\boldsymbol{\mathcal{P}}}}

\newcommand{\Rbb}{\ensuremath{\mathbb{R} }}

\newcommand\Ibm{{\ensuremath{\bm{I}}}}

\newcommand\Mbm{{\ensuremath{\bm{M}}}}

\newcommand\cbm{{\ensuremath{\bm{c}}}}

\newcommand\fbm{{\ensuremath{\bm{f}}}}
\newcommand\gbm{{\ensuremath{\bm{g}}}}

\newcommand\kbm{{\ensuremath{\bm{k}}}}

\newcommand\pbm{{\ensuremath{\bm{p}}}}
\newcommand\qbm{{\ensuremath{\bm{q}}}}
\newcommand\rbm{{\ensuremath{\bm{r}}}}

\newcommand\ubm{{\ensuremath{\bm{u}}}}

\newcommand\xbm{{\ensuremath{\bm{x}}}}

\newcommand\mubold{{\ensuremath{\boldsymbol{\mu}}}}

\newcommand\Sigmabold{{\ensuremath{\boldsymbol{\Sigma}}}}

\newcommand\zerobold{\ensuremath{\mathbf{0}}}

\usepackage{xifthen}
\usepackage[thinlines]{easytable}

\usepackage{setspace}

\newbool{fastcompile}
\setbool{fastcompile}{false}

\newcommand{\dmu}[1][]{\ifthenelse{\isempty{#1}}{\partial_\mubold}{\pder{#1}{\mubold}}}
\newcommand{\mass}[1][]{\ifthenelse{\isempty{#1}}{\Mbm}{\Mbm^{#1}}}
\newcommand{\stvc}[1][]{\ifthenelse{\isempty{#1}}{\ubm}{\ubm^{#1}}}
\newcommand{\stvcdot}[1][]{\ifthenelse{\isempty{#1}}{\dot\ubm}{\dot\ubm^{#1}}}
\newcommand{\stvcbar}[1][]{\ifthenelse{\isempty{#1}}{\bar\ubm}{\bar\ubm^{#1}}}
\newcommand{\res}[1][]{\ifthenelse{\isempty{#1}}{\rbm}{\rbm^{#1}}}
\newcommand{\resimpl}[1][]{\ifthenelse{\isempty{#1}}{\gbm}{\gbm^{#1}}}
\newcommand{\resexpl}[1][]{\ifthenelse{\isempty{#1}}{\fbm}{\fbm^{#1}}}
\newcommand{\cpl}[1][]{\ifthenelse{\isempty{#1}}{\cbm}{\cbm^{#1}}}
\newcommand{\cplprd}[1][]{\ifthenelse{\isempty{#1}}{\tilde\cbm}{\tilde\cbm^{#1}}}

\newcommand{\allstvc}[0]{\stvc[1],\,\dots,\,\stvc[m]}

\newcommand{\pstp}[2][]{\ifthenelse{\isempty{#2}}{\ubm_{#1}}{\ubm_{#2}^{#1}}}
\newcommand{\pstgki}[3][]{\ifthenelse{\isempty{#3}}{\kbm_{#1,#2}}{\kbm_{#2,#3}^{#1}}}
\newcommand{\pstgke}[3][]{\ifthenelse{\isempty{#3}}{\hat\kbm_{#1,#2}}{\hat\kbm_{#2,#3}^{#1}}}
\newcommand{\pstgu}[3][]{\ifthenelse{\isempty{#3}}{\ubm_{#1,#2}}{\ubm_{#2,#3}^{#1}}}
\newcommand{\cplstp}[2][]{\ifthenelse{\isempty{#2}}{\cbm_{#1}}{\cbm_{#2}^{#1}}}
\newcommand{\cplprdstp}[2][]{\ifthenelse{\isempty{#2}}{\tilde\cbm_{#1}}{\tilde\cbm_{#2}^{#1}}}

\newcommand{\dt}[1]{\Delta t_{#1}}

\newcommand{\tstgi}[2]{t_{#1,#2}}

\newcommand{\qstg}[3][]{\ifthenelse{\isempty{#3}}{\qbm_{#1,#2}}{\qbm_{#2,#3}^{#1}}}
\newcommand{\prstp}[2][]{\ifthenelse{\isempty{#2}}{\pbm_{#1}}{\pbm_{#2}^{#1}}}

\newcommand{\eqnref}[1]{Eq. (\ref{#1})}
\newcommand{\figref}[1]{Figure \ref{#1}}
\newcommand{\tabref}[1]{Table \ref{#1}}

\title{High-order, linearly stable, partitioned solvers for general
       multiphysics problems based on implicit-explicit Runge-Kutta
       schemes}

\author[rvt1]{D.~Z.~Huang\fnref{fn1}\corref{cor1}}
\ead{zhengyuh@stanford.edu}

\author[rvt2,rvt3]{P.-O.~Persson\fnref{fn2}}
\ead{persson@berkeley.edu}

\author[rvt2,rvt4]{M.~J.~Zahr\fnref{fn3,fn4}\corref{cor1}}
\ead{mjzahr@lbl.gov}

\address[rvt1]{Institute for Computational and Mathematical Engineering,
               Stanford University, Stanford, CA, 94305, United States}
\address[rvt2]{Mathematics Group, Lawrence Berkeley National Laboratory,
               1 Cyclotron Road, Berkeley, CA 94720, United States}
\address[rvt3]{Department of Mathematics, University of California, Berkeley,
               Berkeley, CA 94720, United States}
\address[rvt4]{Department of Aerospace and Mechanical Engineering,
               University of Notre Dame, Notre Dame, IN, 46556, United States}
\cortext[cor1]{Corresponding author}

\fntext[fn1]{Graduate Student, Institute for Computational and Mathematical
             Engineering, Stanford University}
\fntext[fn2]{Associate Professor, Department of Mathematics, University of
             California, Berkeley}
\fntext[fn3]{Luis W. Alvarez Postdoctoral Fellow,
             Computational Research Division,
             Lawrence Berkeley National Laboratory}
\fntext[fn4]{Assistant Professor, Department of Aerospace and Mechanical
             Engineering, University of Notre Dame}

\begin{document}

\begin{abstract}
This work introduces a general framework for constructing high-order,
linearly stable, partitioned solvers for multiphysics problems from
a monolithic implicit-explicit Runge-Kutta (IMEX-RK) discretization
of the semi-discrete equations. The generic multiphysics problem is modeled as
a system of $n$ systems of partial differential equations where the $i$th
subsystem is coupled to the other subsystems through a coupling term that can
depend on the state of all the other subsystems. This coupled system of
partial differential equations reduces to a coupled system of ordinary
differential equations via the method of lines where an appropriate
spatial discretization is applied to each subsystem. The coupled system
of ordinary differential equations is taken as a monolithic system
and discretized using an IMEX-RK discretization with a specific
implicit-explicit decomposition that introduces the concept
of a \emph{predictor} for the coupling term. We propose four
coupling predictors that enable the monolithic system
to be solved in a partitioned manner, i.e., subsystem-by-subsystem,
and preserve the IMEX-RK structure and therefore the design order of
accuracy of the monolithic scheme. The four partitioned solvers that
result from these predictors are high-order accurate, allow for maximum re-use
of existing single-physics software, and two of the four solvers allow
the subsystems to be solved in parallel at a given stage and time step.
We also analyze the stability of a coupled, linear model problem 
with a specific coupling structure and show that one of the
partitioned solvers achieves unconditional linear stability
for this problem,
while the others are unconditionally stable only for certain values of the
coupling strength. We demonstrate the performance of the
proposed partitioned solvers on several classes of multiphysics problems
including a simple linear system of ODEs, advection-diffusion-reaction
systems, fluid-structure interaction problems, and particle-laden flows, where
we verify the design order of the IMEX schemes and study various stability
properties.
\end{abstract}

\maketitle

\section{Introduction} \label{SEC: INTRO}
The numerical simulation of multiphysics problems involving multiple
physical models or multiple simultaneous physical phenomena is significant
in many engineering and scientific applications, e.g., aircraft flutter
in transonic flows \cite{chen2007numerical}, biomedical flows in heart and
blood vessels \cite{griffith2007adaptive}, mixing and chemically reacting
flows \cite{day2000numerical}, reactor fuel performance
\cite{gaston2009moose}, turbomachinery \cite{carstens2003coupled},
magnetohydrodynamics \cite{toth2000b} and so on. These problems are generally
highly nonlinear, feature multiple scales and strong coupling effects, and
require heterogeneous discretizations for the various physics subsystems.
To balance the treatment of these features, solution strategies ranging from
a monolithic approach to partitioned procedures have been proposed.

In the monolithic approach
\cite{hubner2004monolithic, michler2004monolithic, hron2006monolithic},
all physical subsystems are solved simultaneously and is therefore preferred
in the case of strong interactions to ensure stability. However, when the
coupled subsystems are complex, the monolithic procedure can be
suboptimal and often requires significant implementation effort since only
small components of existing software can be re-used. An alternative is
the partitioned procedure
\cite{farhat2000two, piperno2001partitioned, badia2008fluid}, also known
as a staggered or a loosely coupled procedure, where different subsystems are
modeled and discretized separately, and the resulting equations are solved
independently. The coupling occurs through specific terms that are lagged
to previous time instances and communicated between solvers. This procedure
facilitates software modularity and mathematical modeling; however, these
schemes are often low-order accurate \cite{piperno2001partitioned} (second order accuracy)
and suffer from lack of stability \cite{causin2005added}. 

Recently, a partitioned solver based on implicit-explicit Runge-Kutta
schemes, first proposed to solve stiff additive ordinary differential
equations \cite{zhong1996additive, ascher1997implicit}, was proposed
\cite{van2007higher, froehle2014high} in the context of a specific
multiphysics system: fluid-structure interaction. This solver demonstrated
up to fifth-order accuracy without requiring the solution of the fully
coupled fluid-structure system. A key feature of this solver that
distinguishes it from other work on IMEX-RK methods for multiphysics
systems \cite{cyr2016implicit} is that both the fluid and the structure subsystems are
handled implicitly and only a correction to the predicted traction on the
structure is treated explicitly. Therefore the stability of this IMEX-RK
partitioned procedure is expected to be better than explicit schemes and
nearly as good as fully implicit schemes. Despite the advantages of the
partitioned fluid-structure interaction solver in
\cite{van2007higher, froehle2014high}, it is not directly applicable to
other multiphysics systems and their proposed traction predictor combines
stage information in a heuristic way, which may lead to accuracy reduction.

Inspired by these works, we built a general framework to construct
high-order partitioned solvers based on monolithic IMEX-RK discretizations
for general multiphysics systems. We consider a general model of multiphysics
problems as a system of $n$ systems of partial differential equations, coupled
through specific coupling terms that can depend on the state of all physical
subsystems, and which is reduced to a system of ODEs after semi-discretization.
An IMEX-RK discretization is applied to this monolithic system of ODEs,
with a specific implicit-explicit decomposition that introduces the
concept of a predictor. The implicit part of the decomposition is taken as
the entire multiphysics system where the coupling term is replaced by the
predictor and the explicit part is a correction to the system that accounts
for errors in the coupling predictor. Predictors that satisfy basic
properties outlined in this work will allow the monolithic discretization
to be solved in a partitioned manner, i.e., subsystem-by-subsystem.
Four consistent predictors are introduced, including weak and strong
Jacobi-type predictors and weak and strong Gauss-Seidel-type predictors, that
lead to different partitioned solvers that maintain the design order of
accuracy of the IMEX-RK scheme. However, the solvers resulting from these
four predictors have their own strengths and limitations by trading off
between implementation effort, stability, and efficiency.
The weak predictors require the least implementation effort since
they do not require any terms from the Jacobian of the coupling term,
while the strong predictors require the diagonal entries.
The Jacobi predictors allow for system-wise parallelization, while the
Gauss-Seidel predictors require the subsystems be solved sequentially. Despite
the simplicity and efficiency of the weak and the Jacobi predictors over the
strong and the Gauss-Seidel predictors, they have weaker linear stability
properties, which we demonstrate through linear stability analysis of
the four predictors applied to a chosen linear model problem
and provide numerical evidence. It is worth noting that, through our linear
stability analysis, we find the strong Gauss-Seidel predictor leads to an
unconditionally stable scheme when applied to the chosen model
problem, despite being a partitioned solver. The splitting choice implied
by the strong Gauss-Seidel predictor minimizes the explicit contribution to
the scheme and, in many cases, the implicit part appears to stabilize growing
modes produced by the explicit part.
Finally, we note that given the generality of this formulation,
             the proposed solver to can be applied to a vast number of
             multiphysics problems; however, it is well-known that partitioned
             solvers are unstable for certain physical regimes, including
             fluid-structure interaction at low mass ratios and
             magnetohydrodynamics with strong bidirectional coupling, which
             is not included by our linear stability analysis. We include a
             more general stability analysis to show the physical regimes,
             i.e., coupling strength, in which the strong Gauss-Seidel
             predictor is unconditionally stable.

The remainder of the paper is organized as follows. In
Section~\ref{SEC: GOVERN EQ}, the general form of the multiphysics problem
as a system of $n$ systems of partial differential equations and its
semi-discretization are introduced. In Section~\ref{SEC: PARTITIONED SOLVER},
an overview of IMEX-RK schemes is provided and a specific implicit-explicit
decomposition using the concept of a coupling predictor is introduced.
Additionally, four predictors are introduced that lead to different solvers
and their features such as accuracy, implementation effort, efficiency,
and stability are discussed. A slew of applications are provided in
Sections~\ref{SEC: APP ODE}-\ref{SEC: APP PARTICLE FLOWS}
that demonstrate the high-order accuracy, stability, and robustness of the
proposed solvers on an advection-diffusion-reaction system, fluid-structure
interaction problems, and particle-laden flows.

\section{Governing multiphysics equations and semi-discretization}
\label{SEC: GOVERN EQ}
Consider a general formulation of a mathematical model describing the
behavior of multiple interacting physical phenomena described by the
following coupled system of partial differential equations
\begin{equation} \label{EQ: GOVERN}
 \partial_t u^i = \Lcal^i(u^i,\,c^i,\,x,\,t),
 \quad x\in\Omega^i(c^i), \quad t \in (0,\,T)
\end{equation}
for $i = 1,\,\dots,\,m$, where $m$ represents the number of physical systems,
and boundary conditions are excluded for brevity.
The $i$th physical system is modeled as a partial differential equation
characterized by the generalized differential operator $\Lcal^i$ that
defines a conservation law or other type of balance law, the state variable
$u^i(x,\,t)$ that is the solution of the $i$th physical system on the
space-time domain $\Omega^i \times (0,\,T)$, and a \emph{coupling term} $c^i$
that, in general, couples the $i$th system to the other $m-1$ systems.
In the general case, the differential operator $\Lcal^i$, domain $\Omega^i$,
and boundary conditions depend on the coupling term.
The coupling term contains quantities usually considered \emph{data} required
to define the $i$th PDE, such as boundary conditions or material properties.
In a single-physics setting, these quantities would be prescribed, but in the
multiphysics setting they are determined from the state vectors of all $m$
systems, i.e.,
\begin{equation} \label{EQ: CPLTERM}
 c^i = c^i(u^1,\,\dots,\,u^m,\,x,\,t).
\end{equation}
The definition of the coupling term is problem-dependent and it will be shown
that special structure in the coupling term can be exploited to create a
better partitioned solver. While the form of (\ref{EQ: GOVERN}) is specific
to first-order temporal systems, it includes equations with higher-order
temporal derivatives, assuming they have been re-cast in first-order form.
The spatial domains $\Omega^i$ for the individual systems may or may not be
overlapping and in many cases are the same, i.e., $\Omega^i = \Omega$
for $i = 1,\,\dots,\,m$. 

As this work is focused on the development of high-order partitioned
schemes for evolving multiphysics problems, we introduce the semi-discrete
form of the coupled partial differential equations in (\ref{EQ: GOVERN}) that
arises from applying an appropriate spatial discretization to the $i$th PDE
system individually, which takes the form
\begin{equation} \label{EQ: GOVERN-SD-INDIV}
 \mass[i] \stvcdot[i] = \res[i](\stvc[i],\,\cpl[i],\,t),
 \quad t \in (0,\,T)
\end{equation}
where $\stvc[i](t)$ is the semi-discrete state vector corresponding to the
spatial discretization of $u^i(x,\,t)$, $\res^i$ is the spatial
discretization of the differential operator $\Lcal^i$ and called the
velocity of the ODE system in the remainder of the document, and $\cpl[i]$ is
the semi-discrete coupling term corresponding to the spatial discretization of
$c^i(u^1,\,\dots,\,u^m,\,x,\,t)$. In general, the coupling term depends on the
semi-discrete state vector of all $m$ systems
\begin{equation} \label{EQ: CPLTERM-SD-INDIV}
 \cpl[i] = \cpl[i](\stvc[1],\,\dots,\,\stvc[m],\,t).
\end{equation}
For convenience, we re-write the system of ordinary differential equations
in (\ref{EQ: GOVERN-SD-INDIV})-(\ref{EQ: CPLTERM-SD-INDIV}) as
\begin{equation} \label{EQ: GOVERN-SD}
 \mass\stvcdot = \res(\stvc,\,\cpl(\stvc,\,t),\,t),
 \quad t \in (0,\,T),
\end{equation}
where the combined mass matrix is a block diagonal matrix consisting of
the single-physics mass matrices
\begin{equation}
 \mass = \begin{bmatrix} \mass[1] && \\ & \ddots & \\ && \mass[m] \end{bmatrix}
\end{equation}
and the combined state vector, coupling term, and nonlinear residual are
vectors consisting of the corresponding single-physics term, concatenated
across all $m$ systems
\begin{equation}
 \stvc = \begin{bmatrix} \stvc^1 \\ \vdots \\ \stvc^m \end{bmatrix} \qquad
 \cpl(\stvc,\,t) = \begin{bmatrix}
                     \cpl^1(\stvc^1,\,\dots,\,\stvc^m,\,t) \\
                     \vdots \\
                     \cpl^m(\stvc^1,\,\dots,\,\stvc^m,\,t)
                  \end{bmatrix} \qquad
 \res(\stvc,\,\cpl,\,t) = \begin{bmatrix}
                            \res^1(\stvc^1,\,\cpl^1,\,t) \\
                            \vdots \\
                            \res^m(\stvc^m,\,\cpl^m,\,t)
                          \end{bmatrix}.
\end{equation}
The total derivative, or Jacobian, of the semi-discrete velocity
$D_{\stvc}\res$ is expanded as
\begin{equation}
 D_{\stvc}\res = \pder{\res}{\stvc} + \pder{\res}{\cpl}\pder{\cpl}{\stvc},
\end{equation}
where the individual terms take the form
\begin{equation} \label{EQ: drdu drdc dcdu}
 \pder{\res}{\stvc} =
 \begin{bmatrix} \displaystyle{\pder{\res[1]}{\stvc[1]}} & & \\
                 & \ddots & \\
                 & & \displaystyle{\pder{\res[m]}{\stvc[m]}}
 \end{bmatrix} \qquad
 \pder{\res}{\cpl} =
 \begin{bmatrix} \displaystyle{\pder{\res[1]}{\cpl[1]}} & & \\
                 & \ddots & \\
                 & & \displaystyle{\pder{\res[m]}{\cpl[m]}}
 \end{bmatrix} \qquad
 \pder{\cpl}{\stvc} =
 \begin{bmatrix} \displaystyle{\pder{\cpl[1]}{\stvc[1]}} & \cdots &
                 \displaystyle{\pder{\cpl[1]}{\stvc[m]}} \\
                 \vdots & \ddots & \vdots \\
                 \displaystyle{\pder{\cpl[m]}{\stvc[1]}} & \cdots &
                 \displaystyle{\pder{\cpl[m]}{\stvc[m]}}
 \end{bmatrix},
\end{equation}
and the dependencies have been dropped for brevity.
The first term in the Jacobian is block diagonal and accounts for the direct
contribution of a state to its own system while the second term accounts
for the coupling between systems. Several examples of this general multiphysics
formulation are provided in
Sections~\ref{SEC: APP ADR}-\ref{SEC: APP PARTICLE FLOWS} including
advection-diffusion-reaction systems, two-field and three-field fluid-structure
interaction problems, and particle-laden flows. The semi-discrete forms of the
multiphysics problem in (\ref{EQ: GOVERN-SD-INDIV}) and (\ref{EQ: GOVERN-SD})
will be the point of departure for the remainder of this document and the
starting point for the introduction of our proposed high-order partitioned
solvers.

\section{A high-order partitioned solver for multiphysics problems}
\label{SEC: PARTITIONED SOLVER}
In this section, we introduce our proposed high-order partitioned
time-integration scheme for multiphysics systems. As discussed in
Section~\ref{SEC: INTRO}, a partitioned solver combines individual,
single-physics solvers into an integration scheme for the multiphysics
problem, rather than considering the monolithic multiphysics system.
A partitioned solver can reduce computational complexity per time-step,
improve software maintainability, and exploit off-the-shelf software
components; however they tend to be limited to low-order accuracy and
have stringent stability requirements. Our partitioned time-integration
scheme mitigates most of these issues by combining high-order
implicit-explicit Runge-Kutta (IMEX) schemes for the monolithic multiphysics
system with a judicious implicit-explicit decomposition that partially
decouples the individual systems via a novel predictor for the coupling terms.

\subsection{Background: implicit-explicit Runge-Kutta schemes}
Implicit-explicit Runge-Kutta schemes, first proposed in \cite{zhong1996additive, ascher1997implicit}, define a
family of high-order discretizations for nonlinear differential equations
whose velocity term can be decomposed into a sum of a non-stiff $\resexpl$
and stiff $\resimpl$ velocity
\begin{equation} \label{EQ: ODEADDSPLIT}
 \mass\dot\stvc = \resexpl(\stvc,\,t) + \resimpl(\stvc,\,t).
\end{equation}
The non-stiff $\resexpl$ velocity is integrated with an $s$-stage explicit
Runge-Kutta scheme and the stiff term $\resimpl$ is integrated with an
$s$-stage diagonally implicit Runge-Kutta scheme.
IMEX Runge-Kutta schemes are compactly represented by a double tableau in
the usual Butcher notation (\tabref{TAB: IMEX}), where $\hat{A}$, $\hat{b}$,
$\hat{c}$ defines the Butcher tableau for the explicit Runge-Kutta
scheme used for $\resexpl$ and $A$, $b$, $c$ defines the diagonally implicit
Runge-Kutta scheme used for $\resimpl$.
In this work, we mainly consider 2nd-order 2-stage trapezoidal rule,
3rd-order 4-stage ARK3(2)4L[2]SA, and 4th-order 6 stage ARK4(3)6L[2]SA
proposed in \cite{christopher2001additive}. Theoretically, IMEX schemes can be
of arbitrarily high order accuracy, but as the number of stages increases,
inexactness of the coefficients may destroy the order of accuracy. The
implicit Runge-Kutta part of these IMEX schemes are L-stable,
stiffly-accurate, and have an explicit first stage ($a_{11} = 0$).

\begin{table}[!htb]
\begin{minipage}{.5\linewidth}
\caption*{Explicit Runge-Kutta coefficients}
\centering
\begin{tabular}{c|ccccc}
0\\
$\hat c_2$    & $\hat a_{21}$\\
$\hat c_3$    & $\hat a_{31}$ & $\hat a_{32}$\\
$\vdots$ & $\vdots$ &           & $\ddots$\\
$\hat c_s$    & $\hat  a_{s1}$ & $a_{s2}$  & $\cdots$ & $\hat a_{ss - 1}$\\
\hline
& $\hat b_1$    & $\hat b_2$     & $\cdots$ & $\hat b_{s-1}$     & $\hat b_s$\\
\end{tabular}
\end{minipage}
\begin{minipage}{.5\linewidth}
\centering
\caption*{Implicit Runge-Kutta coefficients}
 \begin{tabular}{c|ccccc}
$c_1$\\
$c_2$    & $a_{21}$ & $a_{22}$\\
$c_3$    & $a_{31}$ & $a_{32}$ & $a_{33}$\\
$\vdots$ & $\vdots$ &      &     & $\ddots$\\
$c_s$    & $a_{s1}$ & $a_{s2}$  & $\cdots$ & $a_{ss-1}$ & $a_{ss}$\\
\hline
& $b_1$    & $b_2$     & $\cdots$ & $b_{s-1}$     & $b_s$\\
\end{tabular}
 \end{minipage} 
\caption{Butcher Tableau for s-stage implicit-explicit Runge-Kutta scheme}   
\label{TAB: IMEX} 
\end{table}

Consider a discretization of the time domain $[0,\,T]$ into $N_t$ segments
with endpoints $\{t_0,\,\dots,\,t_{N_t}\}$, with the $n$th segment having
length $\Delta t_n = t_n - t_{n-1}$ for $n = 1,\,\dots,\,N_t$. Also,
let $\pstp{n}$ denote the approximation of the solution of the
differential equation in (\ref{EQ: ODEADDSPLIT}) at time step $n$, i.e.,
$\pstp{n} \approx \stvc(t_n)$. Then, given the explicit
$(\hat{A},\,\hat{b},\,\hat{c})$ and implicit $(A,\,b,\,c)$ Butcher tableaus,
the $s$-stage IMEX Runge-Kutta scheme that advances $\pstp{n-1}$ to
$\pstp{n}$ is given by
\begin{subequations}
\label{EQ: IMEX STEP IMP EXP STAGE}
 \begin{align}
  \pstp{n} &= \pstp{n-1} + \sum_{p=1}^s\hat{b}_p\pstgke{n}{p} +
                            \sum_{p=1}^s b_p\pstgki{n}{p},
                                                        \label{EQ: IMEX STEP} \\
  \mass\pstgki{n}{j} &= \dt{n}\resimpl(\pstgu{n}{j},\,t_{n-1}+c_j\dt{n}),
                                                         \label{EQ: IMEX IMP} \\
  \mass\pstgke{n}{j} &= \dt{n}\resexpl(\pstgu{n}{j},\,t_{n-1}+\hat{c}_j\dt{n}),
                                                         \label{EQ: IMEX EXP} \\
  \pstgu{n}{j} &= \pstp{n-1} + \sum_{p=1}^{j-1}\hat{a}_{jp}\pstgke{n}{p} +
                               \sum_{p=1}^j a_{jp}\pstgki{n}{p},
                                                         \label{EQ: IMEX STAGE}
 \end{align}
\end{subequations}
where $\pstgke{n}{p}$ and $\pstgki{n}{p}$ are the $p$th explicit and implicit
velocity stage, respectively, corresponding to time step $n$ and
$\pstgu{n}{p}$ is the approximation to $\pstp{n}$ at stage $p$ of time step $n$.
The complete algorithm to advance $\pstp{n-1}$ to $\pstp{n}$ using the IMEX
Runge-Kutta scheme is provided in Algorithm~\ref{ALG: IMEX}. For each
stage $j$, the nonlinear system of equations in (\ref{EQ: IMEX IMP}) must
be solved to compute the implicit stage $\pstgki{n}{j}$. Next, the explicit
stage can be computed directly from (\ref{EQ: IMEX EXP}) since the stage
approximation $\pstgu{n}{j}$ does not depend on the explicit stage
$\pstgke{n}{j}$. Finally, given the previous time step and all implicit
and explicit stages, the solution at time $n$ is determined from
(\ref{EQ: IMEX STEP}).

\begin{algorithm}
 \caption{Implicit-Explicit Runge-Kutta scheme}
 \label{ALG: IMEX}
 \begin{algorithmic}[1]
   \For{stages $j = 1,\,\dots,\,s$}
    \State Define stage solution according to (\ref{EQ: IMEX STEP}):
     $\displaystyle{\pstgu{n}{j} = \pstp{n-1} +
                                   \sum_{p=1}^{j-1}\hat{a}_{jp}\pstgke{n}{p} +
                                   \sum_{p=1}^j a_{jp}\pstgki{n}{p}}$
    \State Implicit solve (\ref{EQ: IMEX IMP}) for $\pstgki{n}{j}$:
     $\displaystyle{\mass\pstgki{n}{j} =
                    \dt{n}\resimpl(\pstgu{n}{j},\,t_{n-1}+c_j\dt{n})}$
    \State Explicit solve (\ref{EQ: IMEX EXP}) for $\pstgke{n}{j}$: 
     $\displaystyle{\mass\pstgke{n}{j} =
                    \dt{n}\resexpl(\pstgu{n}{j},\,t_{n-1}+\hat{c}_j\dt{n})}$
   \EndFor
   \State Set 
     $\displaystyle{\pstp{n} = \pstp{n-1} + \sum_{p=1}^s\hat{b}_p\pstgke{n}{p} +
                                            \sum_{p=1}^s b_p\pstgki{n}{p}}$
 \end{algorithmic}
\end{algorithm}

\subsection{A partitioned implicit-explicit Runge-Kutta scheme for
            multiphysics systems}
The proposed high-order partitioned scheme for integration of generic
time-dependent multiphysics problems of the form
(\ref{EQ: GOVERN-SD-INDIV})-(\ref{EQ: CPLTERM-SD-INDIV}) is built on an
IMEX Runge-Kutta discretization of the monolithic system.
A special choice of implicit-explicit decomposition, along with the
introduction of four predictors for the coupling term, creates a
\emph{diagonal} (uncoupled) or \emph{triangular} dependency between the
systems and allows the monolithic discretization to be solved in
a partitioned manner. The proposed decomposition handles a majority of
the relevant physics \emph{implicitly} to leverage the enhanced stability
properties of such schemes, while only the correction to the coupling
predictor is handled explicitly. This marks a key difference to previous
work on IMEX Runge-Kutta schemes for multiphysics systems \cite{cyr2016implicit} that choose
the implicit-explicit splitting based on stiffness of the physics.
It will be shown in Section~\ref{SEC: ACCURATE-STABILITY} that the proposed
predictors preserve the accuracy and stability properties of the IMEX scheme.

\subsubsection{Implicit-explicit decomposition and monolithic IMEX
               Runge-Kutta discretization}
To begin our construction, recall the semi-discrete form of the multiphysics
system (\ref{EQ: GOVERN-SD}) and consider the splitting of the
velocity term $\res(\stvc,\,\cpl(\stvc,\,t),\,t)$ as
\begin{equation} \label{EQ: RES SPLIT}
 \res(\stvc,\,\cpl(\stvc,\,t),\,t) = \resexpl(\stvc,\,\cplprd,\,t) +
                                     \resimpl(\stvc,\,\cplprd,\,t)
\end{equation}
where $\cplprd$ is an approximation, or \emph{predictor}, of the coupling
term $\cpl(\stvc,\,t)$ and the terms that will be handled explicitly
$\resexpl$ and implicitly $\resimpl$ in the IMEX discretization are defined as
\begin{subequations} \label{EQ: IMEX MP SPLIT}
 \begin{align}
  \resexpl(\stvc,\,\cplprd,\,t) &= \res(\stvc,\,\cpl(\stvc,\,t),\,t) - 
                                   \res(\stvc,\,\cplprd,\,t)
                                   \label{EQ: IMEX MP EXP} \\
  \resimpl(\stvc,\,\cplprd,\,t) &= \res(\stvc,\,\cplprd,\,t),
                                   \label{EQ: IMEX MP IMP}
 \end{align}
\end{subequations}
where the dependence on the predictor is explicitly included. In general,
the predictor depends on the instantaneous state vector $\stvc(t)$ and
data $\stvcbar$, likely from the history of the state vector
$\{\stvc(\tau)\,|\,\tau < t\}$
\begin{equation} \label{EQ: CPLPRD GEN}
 \cplprd = \cplprd(\stvc,\,\stvcbar,\,t).
\end{equation}

With this decomposition of the velocity of the semi-discrete multiphysics
system in (\ref{EQ: IMEX MP SPLIT}), the IMEX Runge-Kutta scheme in
(\ref{EQ: IMEX STEP IMP EXP STAGE}) applied to the monolithic multiphysics
system (\ref{EQ: GOVERN-SD}) becomes
\begin{equation}
 \begin{aligned}
  \pstp{n} &= \pstp{n-1} + \sum_{p=1}^s\hat{b}_p\pstgke{n}{p} +
                            \sum_{p=1}^s b_p\pstgki{n}{p}, \\
  \mass\pstgki{n}{j} &= \dt{n}\resimpl(\pstgu{n}{j},\,
                                   \cplprd(\pstgu{n}{j},\,\pstp{n-1},\,
                                           \tstgi{n}{j}),\,
                                   \tstgi{n}{j}),\\
  \mass\pstgke{n}{j} &= \dt{n}\resexpl(\pstgu{n}{j},\,
                                   \cplprd(\pstgu{n}{j},\,\pstp{n-1},\,
                                           \tstgi{n}{j}),\,
                                   \tstgi{n}{j}),\\
  \pstgu{n}{j} &= \pstp{n-1} + \sum_{p=1}^{j-1}\hat{a}_{jp}\pstgke{n}{p} +
                               \sum_{p=1}^j a_{jp}\pstgki{n}{p},
 \end{aligned}
\end{equation}
where the data used in the coupling predictor is taken from the previous
time step. This is the general form of the fully discrete, monolithic
multiphysics system where the coupling predictor is unspecified.
In the general setting where each coupling predictor depends on the state
of all systems, the Jacobian of the coupling predictor is block
dense with potentially sparse blocks
\begin{equation}
 \pder{\cplprd}{\stvc} =
 \begin{bmatrix} \displaystyle{\pder{\cplprd[1]}{\stvc[1]}} & \cdots &
                 \displaystyle{\pder{\cplprd[1]}{\stvc[m]}} \\
                 \vdots & \ddots & \vdots \\
                 \displaystyle{\pder{\cplprd[m]}{\stvc[1]}} & \cdots &
                 \displaystyle{\pder{\cplprd[m]}{\stvc[m]}}
 \end{bmatrix}.
\end{equation}
This implies the Jacobian of the implicit velocity
\begin{equation}
 D_{\stvc} \resimpl = \pder{\res}{\stvc} +
                      \pder{\res}{\cplprd} \pder{\cplprd}{\stvc}
\end{equation}
is also block dense, which highlights the fact that there is coupling across
all systems and a monolithic solver is required for the implicit step.
The next section will introduce four coupling predictors that reduce the
monolithic nature of the multiphysics IMEX-RK discretization to a
partitioned scheme.

\subsubsection{Four coupling predictors and reduction to partitioned schemes}
\label{SEC: 4 PREDICTORS}
To arrive at a scheme that can be solved in a \emph{partitioned} way,
we introduce four predictors that break the monolithic nature of the
multiphysics system. The proposed predictors will first be classified as
leading to a \emph{weak} or \emph{strong} coupling
depending on whether the diagonal of the coupling predictor Jacobian is
nonzero, i.e.,
\begin{equation}
 \pder{\cplprd[i]}{\stvc[i]} = 0 \quad \emph{weakly coupled}, \qquad
 \pder{\cplprd[i]}{\stvc[i]} \neq 0 \quad \emph{strongly coupled}
\end{equation}
for $i = 1,\,\dots,\,m$. In other words, for the weakly coupled predictor, the predicted interaction force $\cplprd[i]$ is constant with respect to the subsystem state $\stvc[i]$.
The predictors will further be classified based on whether they lead to a
Jacobi-type (diagonal) or Gauss-Seidel-type (triangular) coupling, i.e.,
\begin{equation}
 \pder{\cplprd[i]}{\stvc[j]} = 0 \quad i \neq j \quad \emph{Jacobi-type}, \qquad
 \pder{\cplprd[i]}{\stvc[j]} = 0 \quad i<j \quad \emph{Gauss-Seidel-type}.
\end{equation}
The remainder of this section is devoted to detailing the four predictors,
the partitioned IMEX schemes that arise, and the advantages and disadvantages
of each.

\subsubsection*{Weakly coupled Jacobi-type predictor}
The first and simplest predictor is the weakly coupled Jacobi-type predictor
that does not consider the instantaneous solution for \emph{any} of the
systems and only considers time history data, i.e.,
\begin{equation} \label{EQ: CPLPRD WEAK JACOBI}
 \cplprd(\stvc,\,\stvcbar,\,t) = \cpl(\stvcbar,\,t).
\end{equation}
At the fully discrete level, this predictor takes the form
\begin{equation} \label{EQ: CPLPRD WEAK JACOBI FD}
 \cplprd(\pstgu{n}{j},\,\pstp{n-1},\,t) = \cpl(\pstp{n-1},\,t),
\end{equation}
where $\pstp{n}$ is the multiphysics state vector at time step $n$ (the
previous time step) and $\pstgu{n}{j}$ is the approximation to $\pstp{n+1}$
at stage $j$ of time step $n$, as defined in
(\ref{EQ: IMEX STEP IMP EXP STAGE}). In the context of the IMEX-RK
discretization in (\ref{EQ: IMEX STEP IMP EXP STAGE}), this predictor
corresponds to lagging the coupling term to the previous time step throughout
all stages of the time step. With this predictor, the IMEX-RK discretization
of the multiphysics system in (\ref{EQ: IMEX STEP IMP EXP STAGE}) leads to
Algorithm~\ref{ALG: IMEX WEAK JACOBI}.
\begin{algorithm}[H]
 \caption{Implicit-Explicit Runge-Kutta partitioned multiphysics scheme:
          weak Jacobi predictor}
 \label{ALG: IMEX WEAK JACOBI}
 \begin{algorithmic}[1]
  \For{stages $j = 1,\,\dots,\,s$}
   \For {physical systems $i = 1,\,\dots,\,m$}
     \State Define stage solution according to (\ref{EQ: IMEX STEP}):
      $\displaystyle{\pstgu[i]{n}{j} = \pstp[i]{n-1} +
                                  \sum_{p=1}^{j-1}\hat{a}_{jp}\pstgke[i]{n}{p} +
                                  \sum_{p=1}^j a_{jp}\pstgki[i]{n}{p}}$
     \State Implicit solve (\ref{EQ: IMEX IMP}) for $\pstgki[i]{n}{j}$:
      $\displaystyle{\mass[i]\pstgki[i]{n}{j} =
                     \dt{n}\resimpl[i](\pstgu[i]{n}{j},\,
                                    \cpl[i](\pstp[1]{n-1},\,\dots,\,\pstp[m]{n-1},\,
                                            \tstgi{n}{j}),\,
                                    \tstgi{n}{j})}$
     \State Explicit solve (\ref{EQ: IMEX EXP}) for $\pstgke[i]{n}{j}$:
      $\displaystyle{\mass[i]\pstgke[i]{n}{j} =
                     \dt{n}\resexpl[i](\pstgu[i]{n}{j},\,
                                    \cpl[i](\pstp[1]{n-1},\,\dots,\,\pstp[m]{n-1},\,
                                            \tstgi{n}{j}),\,
                                    \tstgi{n}{j})}$
    \EndFor
   \EndFor
   \State Set
     $\displaystyle{\pstp{n} = \pstp{n-1} + \sum_{p=1}^s\hat{b}_p\pstgke{n}{p} +
                                            \sum_{p=1}^s b_p\pstgki{n}{p}}$
 \end{algorithmic}
\end{algorithm}

The IMEX-RK discretization with this coupling predictor is interpreted as, at
each stage, an implicit solve that simultaneously accounts for all physics
systems with the coupling term lagged one time step and corrected by an
explicit step that accounts for the error introduced due to this lagged
coupling term. Furthermore, this choice of predictor leads to a Jacobi-type
decoupling of the various systems during the implicit solve at a given stage,
thus allowing the monolithic system to be solved in a partitioned manner.
This can easily be seen from the fact that the implicit Jacobian
$D_{\stvc} \resimpl$ is block diagonal
\begin{equation}
 D_{\stvc} \resimpl = \pder{\res}{\stvc}
\end{equation}
since the Jacobian of the coupling predictor is zero and
$\displaystyle{\pder{\res}{\stvc}}$ is block diagonal from
(\ref{EQ: drdu drdc dcdu}).

There are a number of advantages surrounding the weak Jacobi-type coupling
predictor, mostly pertaining to simplicity and efficiency. First, the
implicit Jacobian (\ref{EQ: IMEX MP IMP}) does not require the Jacobian of
the coupling term, which can be cumbersome to implement, particularly when
used to patch together existing software to form a multiphysics tool.
Additionally, this simple predictor allows maximum re-use of single-physics
software since only the coupling term must be communicated between codes
to implement the multiphysics partitioned scheme. Once communication of
the coupling term is complete at the beginning of a time step, the
Jacobi-type coupling implies that, within a given time step, all
single-physics systems are independent and can be performed in parallel.
Finally, since the partitioned discretization is a special case of the
IMEX-RK discretization in (\ref{EQ: IMEX STEP IMP EXP STAGE}), it is
guaranteed to preserve the design order of the discretization; see
Section~\ref{SEC: ACCURATE-STABILITY} for a detailed discussion. The primary
disadvantage of this simple and efficient predictor is reduced stability
properties, which will be discussed further in
Section~\ref{SEC: ACCURATE-STABILITY}.

\subsubsection*{Strongly coupled Jacobi-type predictor}
A predictor that maintains a Jacobi-type coupling, i.e., block diagonal
implicit Jacobian, while incorporating additional instantaneous information
is defined as
\begin{equation} \label{EQ: CPLPRD STRONG JACOBI}
 \cplprd[i](\stvc,\,\stvcbar,\,t) = \cpl(\stvcbar[1],\,\dots,\,\stvcbar[i-1],\,
                                         \stvc[i],\,
                                         \stvcbar[i+1],\,\dots,\,\stvcbar[m],\,
                                         t).
\end{equation}
for $i = 1,\,\dots,\,m$. At the fully discrete level, this predictor takes
the form
\begin{equation} \label{EQ: CPLPRD STRONG JACOBI FD}
 \cplprd[i](\pstgu{n}{j},\,\pstp{n-1},\,t) =
                                     \cpl(\pstp[1]{n-1},\,\dots,\,\pstp[i-1]{n-1},\,
                                          \pstgu[i]{n}{j},\,
                                          \pstp[i+1]{n-1},\,\dots,\,\pstp[m]{n-1},\,
                                          t).
\end{equation}
This leads to a strong coupling where
$\displaystyle{\pder{\cplprd[i]}{\stvc[i]}\neq 0}$. In the context of
the IMEX-RK discretization in (\ref{EQ: IMEX STEP IMP EXP STAGE}), the strong
Jacobi predictor corresponds to using the instantaneous state for the $i$th
system in the $i$th coupling term and lagging the remaining states to the
previous time step. With this predictor, the IMEX-RK discretization of the
multiphysics system leads to
Algorithm~\ref{ALG: IMEX STRONG JACOBI}.
\begin{algorithm}[H]
 \caption{Implicit-Explicit Runge-Kutta partitioned multiphysics scheme:
          strong Jacobi predictor}
 \label{ALG: IMEX STRONG JACOBI}
 \begin{algorithmic}[1]
    \For{stages $j = 1,\,\dots,\,s$}
     \For {physical systems $i = 1,\,\dots,\,m$}
     \State Define stage solution according to (\ref{EQ: IMEX STEP}):
      $\displaystyle{\pstgu[i]{n}{j} = \pstp[i]{n-1} +
                                  \sum_{p=1}^{j-1}\hat{a}_{jp}\pstgke[i]{n}{p} +
                                  \sum_{p=1}^j a_{jp}\pstgki[i]{n}{p}}$
     \State Implicit solve (\ref{EQ: IMEX IMP}) for $\pstgki[i]{n}{j}$:
      $\displaystyle{\mass[i]\pstgki[i]{n}{j} =
                     \dt{n}\resimpl[i](\pstgu[i]{n}{j},\,
                                    \cpl[i](\pstp[1]{n-1},\,\dots,\,
                                            \pstp[i-1]{n-1},\,
                                            \pstgu[i]{n}{j},\,
                                            \pstp[i+1]{n-1},\,\dots,\,
                                            \pstp[m]{n-1},\,
                                            \tstgi{n}{j}),\,
                                    \tstgi{n}{j})}$
     \State Explicit solve (\ref{EQ: IMEX EXP}) for $\pstgke[i]{n}{j}$:
      $\displaystyle{\mass[i]\pstgke[i]{n}{j} =
                     \dt{n}\resexpl[i](\pstgu[i]{n}{j},\,
                                    \cpl[i](\pstp[1]{n-1},\,\dots,\,
                                            \pstp[i-1]{n-1},\,
                                            \pstgu[i]{n}{j},\,
                                            \pstp[i+1]{n-1},\,\dots,\,
                                            \pstp[m]{n-1},\,
                                            \tstgi{n}{j}),\,
                                    \tstgi{n}{j})}$
    \EndFor
   \EndFor
   \State Set
     $\displaystyle{\pstp{n} = \pstp{n-1} + \sum_{p=1}^s\hat{b}_p\pstgke{n}{p} +
                                            \sum_{p=1}^s b_p\pstgki{n}{p}}$
 \end{algorithmic}
\end{algorithm}

The interpretation of the IMEX-RK discretization with the strong Jacobi
predictor is similar to that of the weak Jacobi predictor with the exception
that the coupling term is not entirely lagged to the previous time step.
The Jacobian of the coupling predictor in (\ref{EQ: CPLPRD STRONG JACOBI})
is block diagonal
\begin{equation}
 \pder{\cplprd}{\stvc} = \begin{bmatrix}
                          \displaystyle{\pder{\cpl[1]}{\stvc[1]}} & & \\
                          & \ddots & \\
                          & & \displaystyle{\pder{\cpl[m]}{\stvc[m]}} \\
                         \end{bmatrix},
\end{equation}
which leads to a block diagonal implicit Jacobian
\begin{equation}
 D_{\stvc}\resimpl = \begin{bmatrix}
                      \displaystyle{\pder{\res[1]}{\stvc[1]} +
                                    \pder{\res[1]}{\cpl[1]}
                                    \pder{\cpl[1]}{\stvc[1]}} & & \\
                      & \ddots & \\
                      & & \displaystyle{\pder{\res[m]}{\stvc[m]} +
                                        \pder{\res[m]}{\cpl[m]}
                                        \pder{\cpl[m]}{\stvc[m]}}
                     \end{bmatrix}.
\end{equation}

The strong Jacobi predictor shares some of the advantages as the weak Jacobi
predictor such as a block diagonal implicit Jacobian that allows all systems
to be solved simultaneously and the ability to re-use single physics software
since only the coupling term must be communicated between codes. However, the
$i$th system now requires the Jacobian of its own coupling term with respect to
its own state, a term that may not be readily available or have an obvious
data structure. The strong Jacobi predictor is guaranteed to preserve the
design order of the IMEX-RK discretization and has better stability properties
than the weak Jacobi predictor; see Section~\ref{SEC: ACCURATE-STABILITY} for
a detailed discussion.

\subsubsection*{Weakly coupled Gauss-Seidel-type predictor}
The Gauss-Seidel-type (triangular) predictors for the multiphysics system
assume the individual systems are \emph{ordered} in a physically relevant
manner. The preferred ordering is problem-dependent and
a number of examples are provided in Sections~\ref{SEC: APP ODE}-\ref{SEC: APP PARTICLE FLOWS}.
The weakly coupled Gauss-Seidel-type predictor for the $i$th system is
defined as
\begin{equation} \label{EQ: CPLPRD WEAK GS}
 \cplprd[i](\stvc,\,\stvcbar) = \cpl(\stvc[1],\,\dots,\,\stvc[i-1],\,
                                     \stvcbar[i],\,\dots,\,\stvcbar[m])
\end{equation}
for $i = 1,\,\dots,\,m$. At the fully discrete level, this predictor takes
the form
\begin{equation} \label{EQ: CPLPRD WEAK GS FD}
 \cplprd[i](\pstgu{n}{j},\,\pstp{n-1},\,t) =
                             \cpl(\pstgu[1]{n}{j},\,\dots,\,\pstgu[i-1]{n}{j},\,
                                  \pstp[i]{n-1},\,\dots,\,\pstp[m]{n-1}).
\end{equation}
In the context of the IMEX-RK discretization in
(\ref{EQ: IMEX STEP IMP EXP STAGE}), the $i$th predictor
lags the state of systems $i,\,\dots,\,m$ to the previous time step in the
evaluation of the coupling term throughout all stages of the time step.
The IMEX-RK discretization of the multiphysics system in
(\ref{EQ: IMEX STEP IMP EXP STAGE}) with this form of the predictor leads to
Algorithm~\ref{ALG: IMEX WEAK GS}.
\begin{algorithm}[H]
 \caption{Implicit-Explicit Runge-Kutta partitioned multiphysics scheme:
          weak Gauss-Seidel predictor}
 \label{ALG: IMEX WEAK GS}
 \begin{algorithmic}[1]
   \For{stages $j = 1,\,\dots,\,s$}
    \For{physical systems $i = 1,\,\dots,\,m$}
     \State Define stage solution according to (\ref{EQ: IMEX STEP}):
      $\displaystyle{\pstgu[i]{n}{j} =
                                  \pstp[i]{n-1} +
                                  \sum_{p=1}^{j-1}\hat{a}_{jp}\pstgke[i]{n}{p} +
                                  \sum_{p=1}^j a_{jp}\pstgki[i]{n}{p}}$
     \State Implicit solve (\ref{EQ: IMEX IMP}) for $\pstgki[i]{n}{j}$:
      $\displaystyle{\mass[i]\pstgki[i]{n}{j} =
                     \dt{n}\resimpl[i](\pstgu[i]{n}{j},\,
                                    \cpl[i](\pstgu[1]{n}{j},\,\dots,\,
                                            \pstgu[i-1]{n}{j},\,
                                            \pstp[i]{n-1},\,\dots,\,
                                            \pstp[m]{n-1},\,
                                            \tstgi{n}{j}),\,
                                    \tstgi{n}{j})}$
     \State Explicit solve (\ref{EQ: IMEX EXP}) for $\pstgke[i]{n}{j}$:
      $\displaystyle{\mass[i]\pstgke[i]{n}{j} =
                     \dt{n}\resexpl[i](\pstgu[i]{n}{j},\,
                                    \cpl[i](\pstgu[1]{n}{j},\,\dots,\,
                                            \pstgu[i-1]{n}{j},\,
                                            \pstp[i]{n-1},\,\dots,\,
                                            \pstp[m]{n-1},\,
                                            \tstgi{n}{j}),\,
                                    \tstgi{n}{j})}$
    \EndFor
   \EndFor
   \State Set
     $\displaystyle{\pstp{n} = \pstp{n-1} + \sum_{p=1}^s\hat{b}_p\pstgke{n}{p} +
                                            \sum_{p=1}^s b_p\pstgki{n}{p}}$
 \end{algorithmic}
\end{algorithm}

In this case, the Jacobian of the coupling predictor is block strictly lower
triangular
\begin{equation}
 \pder{\cplprd}{\stvc} =
 \begin{bmatrix} 0 & & & \\
                 \displaystyle{\pder{\cpl[2]}{\stvc[1]}} & 0 & & \\
                 \vdots & \ddots & \ddots & \\
                 \displaystyle{\pder{\cpl[m]}{\stvc[1]}} & \cdots &
                 \displaystyle{\pder{\cpl[m]}{\stvc[m-1]}} & 0
 \end{bmatrix},
\end{equation}
which implies the Jacobian of the monolithic implicit system is
block lower triangular
\begin{equation}
 D_{\stvc[j]} \resimpl[i] =
 \begin{cases}
  \displaystyle{\pder{\res[i]}{\stvc[i]}} & i = j \\
  \displaystyle{\pder{\res[i]}{\cpl[i]}
                               \pder{\cpl[i]}{\stvc[j]}} & i > j \\
  \zerobold & i < j.
 \end{cases}
 \label{EQ: JAC WEAK GS}
\end{equation}
This block lower triangular nature of the monolithic implicit system implies
that the individual systems can be solved sequentially beginning with system $1$
and yields a partitioned scheme. 

The implicit Jacobian of the monolithic implicit system of the weak
Gauss-Seidel predictor (\ref{EQ: JAC WEAK GS}) involves the entire lower
triangular portion of the coupling predictor; however, it is
{\it not required} for the implementation. From inspection of
Algorithm~\ref{ALG: IMEX WEAK GS}, the implicit phase at stage $j$ for the $i$th
physical system requires the solution of a nonlinear system of equations in
the variable $\pstgu[i]{n}{j}$, with
$\pstgu[1]{n}{j},\,\dots,\,\pstgu[i-1]{n}{j}$ available from the implicit
solve corresponding to previous physical systems at the current stage.
Therefore, only the \emph{diagonal terms}
$\displaystyle{\frac{D \resimpl[i]}{D \stvc[i]} = \pder{\res[i]}{\stvc[i]}}$
of the monolithic implicit Jacobian are required, which shows that the Jacobians of
the coupling terms are not required for the weak Gauss-Seidel predictor.
Compared with Jacobi-type predictors, the Gauss-Seidel-type predictor also
requires the systems be solved serially within each Runge-Kutta stage and
therefore forfeits the opportunity to parallelize across systems. This
predictor is guaranteed to preserve the design order of the IMEX-RK
discretization and possesses similar stability properties to the weak Jacobi
predictor; see Section~\ref{SEC: ACCURATE-STABILITY}. In
Section~\ref{SEC: APP FSI} we show some desirable properties of the weak
Gauss-Seidel predictor that arise in practice.

\subsubsection*{Strongly coupled Gauss-Seidel-type predictor}
A strong Gauss-Seidel-type coupling is obtained if the $i$th coupling
predictor considers the instantaneous solution for systems $1,\,\dots,\,i$
and the time history for the remaining systems, i.e.,
\begin{equation} \label{EQ: CPLPRD STRONG GS}
 \cplprd[i](\stvc,\,\stvcbar,\,t) = \cpl(\stvc[1],\,\dots,\,\stvc[i],\,
                                        \stvcbar[i+1],\,\dots,\,\stvcbar[m],\,t)
\end{equation}
for $i = 1,\,\dots,\,m$. At the fully discrete level, this predictor takes
the form
\begin{equation} \label{EQ: CPLPRD STRONG GS FD}
 \cplprd[i](\pstgu{n}{j},\,\pstp{n-1},\,t) =
                             \cpl(\pstgu[1]{n}{j},\,\dots,\,\pstgu[i]{n}{j},\,
                                  \pstp[i+1]{n-1},\,\dots,\,\pstp[m]{n-1}).
\end{equation}
In the context of the IMEX-RK discretization in
(\ref{EQ: IMEX STEP IMP EXP STAGE}), the $i$th predictor lags the state of
systems $i+1,\,\dots,\,m$ to the previous time step in the evaluation of the
coupling term throughout all stages of the time step. The IMEX-RK
discretization of the multiphysics system in
(\ref{EQ: IMEX STEP IMP EXP STAGE}) with the strong Gauss-Seidel predictor
becomes Algorithm~\ref{ALG: IMEX STRONG GS}.
\begin{algorithm}[H]
 \caption{Implicit-Explicit Runge-Kutta partitioned multiphysics scheme:
          strong Gauss-Seidel predictor}
 \label{ALG: IMEX STRONG GS}
 \begin{algorithmic}[1]
   \For{stages $j = 1,\,\dots,\,s$}
    \For{physical systems $i = 1,\,\dots,\,m$}
     \State Define stage solution according to (\ref{EQ: IMEX STEP}):
      $\displaystyle{\pstgu[i]{n}{j} =
                                  \pstp[i]{n-1} +
                                  \sum_{p=1}^{j-1}\hat{a}_{jp}\pstgke[i]{n}{p} +
                                  \sum_{p=1}^j a_{jp}\pstgki[i]{n}{p}}$
     \State Implicit solve (\ref{EQ: IMEX IMP}) for $\pstgki[i]{n}{j}$:
      $\displaystyle{\mass[i]\pstgki[i]{n}{j} =
                     \dt{n}\resimpl[i](\pstgu[i]{n}{j},\,
                                    \cpl[i](\pstgu[1]{n}{j},\,\dots,\,
                                            \pstgu[i]{n}{j},\,
                                            \pstp[i+1]{n-1},\,\dots,\,
                                            \pstp[m]{n-1},\,
                                            \tstgi{n}{j}),\,
                                    \tstgi{n}{j})}$
     \State Explicit solve (\ref{EQ: IMEX EXP}) for $\pstgke[i]{n}{j}$:
      $\displaystyle{\mass[i]\pstgke[i]{n}{j} =
                     \dt{n}\resexpl[i](\pstgu[i]{n}{j},\,
                                    \cpl[i](\pstgu[1]{n}{j},\,\dots,\,
                                            \pstgu[i]{n}{j},\,
                                            \pstp[i+1]{n-1},\,\dots,\,
                                            \pstp[m]{n-1},\,
                                            \tstgi{n}{j}),\,
                                    \tstgi{n}{j})}$
    \EndFor
   \EndFor
   \State Set
     $\displaystyle{\pstp{n} = \pstp{n-1} + \sum_{p=1}^s\hat{b}_p\pstgke{n}{p} +
                                            \sum_{p=1}^s b_p\pstgki{n}{p}}$
 \end{algorithmic}
\end{algorithm}

The Jacobian of the coupling predictor is block lower triangular
\begin{equation}
 \pder{\cplprd}{\stvc} =
 \begin{bmatrix} \displaystyle{\pder{\cpl[1]}{\stvc[1]}} & & \\
                 \vdots & \ddots & \\
                 \displaystyle{\pder{\cpl[m]}{\stvc[1]}} & \cdots &
                 \displaystyle{\pder{\cpl[m]}{\stvc[m]}}
 \end{bmatrix},
\end{equation}
which implies that the Jacobian of the monolithic implicit system is also
block lower triangular
\begin{equation}
 D_{\stvc[j]} \resimpl[i] =
 \begin{cases}
  \displaystyle{\pder{\res[i]}{\stvc[i]} +
                               \pder{\res[i]}{\cpl[i]}
                               \pder{\cpl[i]}{\stvc[i]}} & i = j \\
  \displaystyle{\pder{\res[i]}{\cpl[i]}
                               \pder{\cpl[i]}{\stvc[j]}} & i > j \\
  \zerobold & i < j.
 \end{cases}
\end{equation}
Similar to the weak Gauss-Seidel-type predictor, this block lower triangular
nature of the monolithic implicit system implies that the individual systems can
be solved sequentially beginning with system $1$ and yields a partitioned
scheme. 

The strong Gauss-Seidel predictor uses as much current information as
possible while guaranteeing a partitioned scheme and the design accuracy
of the IMEX-RK discretization is not reduced. Similar to the weak Gauss-Seidel
predictor, only the diagonal terms
$\displaystyle{\frac{D \resimpl[i]}{D \stvc[i]} =
               \pder{\res[i]}{\stvc[i]} +
               \pder{\res[i]}{\cpl[i]}\pder{\cpl[i]}{\stvc[i]}}$
of the monolithic implicit Jacobian are required. The implementation effort
is only slightly higher than the weak Gauss-Seidel predictor given that
the diagonal of the coupling Jacobian is required. It will be shown in
Section~\ref{SEC: ACCURATE-STABILITY} that the inclusion of these diagonal
terms leads to enhanced stability properties.

In general, strong coupling predictors include contributions to the block
diagonal from the coupling term, which improves the stability of the
resulting partitioned scheme, but requires more implementation effort than
the weak coupling counterparts. Gauss-Seidel-type predictors lead to
partitioned algorithms where the individual physical subsystems must be
solved sequentially, which reduces their efficiency compared to Jacobi-type
predictors.

\subsubsection{A special case of the coupling structure}
The aforementioned implicit-explicit decomposition with the coupling predictor
(Section~\ref{SEC: 4 PREDICTORS}) is the most general form of the splitting; however, it does
not take advantage of any special structure in the multiphysics problems since we
must predict all $m$ interactions to decouple the multiphysics system. For
many multiphysics problems, such as two-field coupling problems, three-field
fluid-structure-interaction, and magnetohydrodynamics, the following coupling
structure exists
\begin{equation}
\label{EQ: SPE STRU}
 \begin{aligned}
  \cpl[1] &= \cpl[1](\allstvc,\,t) \\
  \cpl[i] &= \cpl[i](\stvc[1],\,\dots,\,\stvc[i],\,t) \qquad i = 2,\,\dots,\,m.
 \end{aligned}
\end{equation}
The strong Gauss-Seidel coupling predictor applied to a coupling term with
the above structure yields
\begin{equation}
 \begin{aligned}
  \cplprd[1](\stvc,\,\stvcbar,\,t) &=
  \cpl(\stvc[1],\,\stvcbar[2],\,\dots,\,\stvcbar[m],\,t) \\
  \cplprd[i](\stvc,\,\stvcbar,\,t) &=
  \cpl(\stvc[1],\,\dots,\,\stvc[i],\,t) \qquad i = 2,\,\dots,\,m
 \end{aligned}
\end{equation}
and therefore the coupling predictors for systems $2,\,\dots,\,m$ are exact,
i.e., identical to the true coupling term. This implies only a single
predictor $\cplprd[1]$ is needed to decouple the multiphysics system and
arrive at a partitioned scheme. In this case, the explicit and implicit terms
of the IMEX-RK scheme reduce to
 \begin{equation}
 \resexpl(\stvc,\,\cplprd,\,t) =
 \begin{bmatrix} \res[1](\stvc[1],\,\cpl[1],\,t) -
                 \res[1](\stvc[1],\,\cplprd[1],\,t) \\ \zerobold \\ \vdots  \\ \zerobold
                  \end{bmatrix}, \qquad
 \resimpl(\stvc,\,\cplprd,\,t) =
 \begin{bmatrix}
   \res[1](\stvc[1],\,\cplprd[1],\,t) \\  \res[2](\stvc[2],\,\cpl[2],\,t) \\
   \vdots \\
   \res[m](\stvc[m],\,\cpl[m],\,t)
 \end{bmatrix}
\end{equation}
In Section ~\ref{SEC: APP ODE}-\ref{SEC: APP PARTICLE FLOWS}, a series of applications that possess this special
coupling structure are presented.

\subsubsection{Accuracy and stability analysis}
\label{SEC: ACCURATE-STABILITY}
The accuracy of implicit-explicit Runge-Kutta schemes is analyzed in detail in
\cite{ascher1997implicit, pareschi2000implicit,
      christopher2001additive, koto2008imex}, where
order conditions are derived from the Taylor expansion of the exact and
numerical solution. Generally, $p$th order IMEX schemes have local truncation
error of $\Ocal(\dt{}^{p+1})$ during one time step
$[t_{n-1},\,t_{n-1}+\dt{}]$ and therefore global temporal error
$\Ocal(\dt{}^p)$. Great care was taken in Section~\ref{SEC: 4 PREDICTORS} to
introduce the proposed predictor-based, partitioned multiphysics scheme as
an implicit-explicit Runge-Kutta discretization to emphasize that the design
order of the IMEX-RK scheme applied to the monolithic multiphysics system is
inherited. This is only possible because the chosen predictors have an
interpretation at the \emph{semi-discrete} level; predictors that combine
stage information in a heuristic way \cite{van2007higher, froehle2014high} may
in general suffer from order reduction, which will be demonstrated in
Section~\ref{SEC: APP ODE}. Therefore, incorporating any of the four proposed
predictors into a $p$th order IMEX-RK schemes leads to the same
$\Ocal(\dt{}^{p+1})$ local truncation error and the same $\Ocal(\dt{}^{p})$
global temporal error.

To study the linear stability of the partitioned IMEX-RK schemes, we consider
the coupled, stable, linear model problem
\begin{equation}
\label{EQ: MODEL-2WAY-SIMPLE}
\begin{aligned}
\partial_t u^1 &= \lambda_1(u^1 + u^2) \\
\partial_t u^2 &= \lambda_2(u^1 + u^2),
\end{aligned}
\end{equation}
where $\Re(\lambda_1) < 0$ and $\Re(\lambda_2) < 0$, that will exhibit the
crux of the linear stability issues for first-order systems
such as advection-diffusion-reaction and
particle-laden flow. However, it does not model complex
bi-directional coupling, e.g.,  characteristic of many problems in
magnetohydrodynamics; see \ref{SEC: APPENDIX-A} for
analysis of a general linear system of ODEs. This system
can be written compactly as
\begin{equation} \label{EQ: ANALYSIS ODE}
 \Mbm \dot\ubm = \rbm(\ubm,\,\cbm(\ubm)),
\end{equation}
where
\begin{equation} \label{EQ: ANALYSIS ODE COMP}
 \Mbm = \begin{bmatrix} 1 & \\ & 1 \end{bmatrix}, \quad
 \ubm = \begin{bmatrix} u^1 \\ u^2 \end{bmatrix}, \quad
 \cbm(\ubm) = \begin{bmatrix}
                c^1(u^1,\,u^2) \\
                c^2(u^1,\,u^2)
              \end{bmatrix}, \quad
 \rbm(\ubm,\,\cbm) =
   \begin{bmatrix}
     (1-\alpha)\lambda_1u^1 + \lambda_1c^1 \\
     (1-\alpha)\lambda_2u^2 + \lambda_2c^2
   \end{bmatrix}
\end{equation}
The coupling terms are chosen as
$c^1(u^1,\,u^2) = \alpha u^1 + u^2$,
$c^2(u^1,\,u^2) = u^1 + \alpha u^2$,
and $\alpha \in \Rbb$ is a coupling parameter that
varies the extent to which an evolution equation
depends on its own state directly through the velocity term
or the coupling term, an important distinction when comparing
the weak and strong predictors. For values of $\alpha$ near unity,
the $i$th evolution equation depends on $u^i$ mostly through the
coupling term, whereas $\alpha$ near zero implies the dependence
is directly through the uncoupled velocity term. The four predictor-based
IMEX schemes introduced in \ref{SEC: 4 PREDICTORS} are applied to this model
problem. The predictor and associated implicit-explicit partition for each
are provided in the \tabref{TAB: MODEL-2WAY-SIMPLE PARTITION}.

\begin{table}[H]
\centering
\begin{tabular}{|>{\centering\arraybackslash}m{2.5cm}|>{\centering\arraybackslash}m{2cm}|>{\centering\arraybackslash}m{5cm}|>{\centering\arraybackslash}m{5cm}|}
\hline
&$\tilde{\cbm}$ &$\gbm$ & $\fbm$\\
Weak Jacobi & \[\left[ \begin{array}{cc} \alpha \bar{u}^1 + \bar{u}^2 \\ \alpha\bar{u}^2 + \bar{u}^1 \end{array}\right]\] & \[ \left[ \begin{array}{cc} \lambda_1((1 -\alpha)u^1 + \alpha \bar{u}^1 + \bar{u}^2) \\ \lambda_2((1 - \alpha)u^2   + \alpha \bar{u}^2 + \bar{u}^1) \end{array}\right]\]  & 
\[\left[ \begin{array}{cc} \lambda_1\alpha(u^1 - \bar{u}^1)  + \lambda_1 (u^2 - \bar{u}^2)  \\ \lambda_2\alpha(u^2 - \bar{u}^2) + \lambda_2(u^1 - \bar{u}^1)   \end{array}\right]\]\\

Strong Jacobi& \[\left[ \begin{array}{cc} \alpha u^1 + \bar{u}^2 \\ \alpha u^2 + \bar{u}^1 \end{array}\right]\]&\[\left[ \begin{array}{cc} \lambda_1(u^1 + \bar{u}^2) \\ \lambda_2(\bar{u}^1 + u^2) \end{array}\right]\] & \[\left[ \begin{array}{cc} \lambda_1(u^2 - \bar{u}^2) \\ \lambda_2(u^1 - \bar{u}^1) \end{array}\right]\]\\

Weak Gauss-Seidel& \[\left[ \begin{array}{cc} \alpha \bar{u}^1 + \bar{u}^2 \\ \alpha\bar{u}^2 + u^1 \end{array}\right]\] &\[\left[ \begin{array}{cc} \lambda_1((1 - \alpha)u^1   + \alpha \bar{u}^1 + \bar{u}^2)\\ \lambda_2((1 - \alpha)u^2 + \alpha \bar{u}^2 + u^1)  \end{array}\right]\] & \[\left[ \begin{array}{cc}  \lambda_1\alpha(u^1 - \bar{u}^1)  + \lambda_1(u^2 - \bar{u}^2)\\ \lambda_2\alpha(u^2 - \bar{u}^2)  \end{array}\right]\]\\

Strong Gauss-Seidel &\[\left[ \begin{array}{cc} \alpha u^1 + \bar{u}^2 \\ \alpha u^2 + u^1 \end{array}\right]\]&\[\left[ \begin{array}{cc} \lambda_1(u^1 + \bar{u}^2) \\ \lambda_2(u^1 + u^2) \end{array}\right]\] & \[\left[ \begin{array}{cc} \lambda_1(u^2 - \bar{u}^2) \\ 0 \end{array}\right]\]\\
\hline
\end{tabular}
\caption{The partition of \eqnref{EQ: MODEL-2WAY-SIMPLE}  based on weak/strong predictors and Jacobi/Gauss-Seidel strategies}
\label{TAB: MODEL-2WAY-SIMPLE PARTITION}
\end{table}

In this section, we consider the 1st-order forward-backward Euler
IMEX scheme \cite{ascher1997implicit}; the linear stability analysis of the
other IMEX schemes considered in this work is provided in \ref{SEC: APPENDIX-A}.
\begin{table}[!htb]
\centering
\begin{minipage}{.35\linewidth}
\centering
\caption*{Explicit Runge-Kutta coefficients}
\begin{tabular}{c|cc}
$0$&   $0$  &  $0$\\
$1$&   $1$  & $0$ \\
\hline
$0$  & $1$ &$0$     \\
\end{tabular}
\end{minipage}
\begin{minipage}{.35\linewidth}
\centering
\caption*{Implicit Runge-Kutta coefficients}
\begin{tabular}{c|cc}
$0$&   $0$  &  $0$\\
$1$&   $0$  &  $1$\\
\hline
$0$ &  $0$ & $1$    \\
\end{tabular}
 \end{minipage} 
\end{table} \\
The forward-backward Euler IMEX scheme applied to the system in
(\ref{EQ: ANALYSIS ODE})-(\ref{EQ: ANALYSIS ODE COMP}) yields the one-step
update equation
\begin{equation}
\label{EQ: IMEX-O1}
\ubm_n = \ubm_{n-1} + \dt{}(\fbm(\ubm_{n-1}) + \gbm(\ubm_n))
\end{equation}
that can be re-written as
\begin{equation}
\label{EQ: ITERATIVE IMEX}
\ubm_n = \Cboldcal(\dt{}, \lambda_1, \lambda_2, \alpha) \ubm_{n-1}. 
\end{equation}
once the partitions in \tabref{TAB: MODEL-2WAY-SIMPLE PARTITION} are introduced.
An update equation of this form is stable if the spectral radius
of the matrix satisfies
\begin{equation}
\rho(\Cboldcal) \leq 1,
\end{equation}
and the multiplicity of any eigenvalues of magnitude 1 is equal to the dimension of its eigenspace.

The spectral radius and the region of unconditional stability for the 1st-order
IMEX scheme based on the partitions in
\tabref{TAB: MODEL-2WAY-SIMPLE PARTITION} are provided in
\tabref{TAB: MODEL-2WAY-SIMPLE O1 IMEX STABILITY}.
Both strong predictors lead to unconditional stability, regardless of the
value of $\alpha$, while the stable regions for the weak predictors depend
on the coupling strength. The weak Gauss-Seidel predictor has a larger
$\alpha$-range of unconditional stability than the weak Jacobi predictor
which is only stable for $\alpha \leq 0$. The high-order IMEX-RK schemes
considered in this
work are analyzed in \ref{SEC: APPENDIX-A} and the strong Gauss-Seidel
predictors are shown to be unconditionally stable, regardless of $\alpha$,
while the strong Jacobi predictor is not.
\begin{table}[H]
\centering
\begin{tabular}{|>{\centering\arraybackslash}m{3cm}|>{\centering\arraybackslash}m{8cm}|>{\centering\arraybackslash}m{3cm}|}
\hline
&Spectral radius&  Unconditional stability range \\
Weak Jacobi&\[\max\Big\{ 1, \Big|\frac{(1 + \alpha \dt{}  \lambda_1)(1 + \alpha \dt{}  \lambda_2) - \dt{}^2 \lambda_1\lambda_2}{(1 - (1 - \alpha)\dt{} \lambda_2)(1 - (1 - \alpha)\dt{} \lambda_1)} \Big| \Big\}\] &$\alpha \leq 0$\\
Strong Jacobi &\[\max\Big\{ 1, \Big|\frac{1 - \dt{}^2 \lambda_1\lambda_2}{(1 - \dt{} \lambda_1)(1 - \dt{} \lambda_2)} \Big| \Big\}\] &$\forall \alpha$ \\
Weak Gauss-Seidel &\[\max\Big\{ 1, \Big|\frac{(1 + \alpha \dt{}  \lambda_1)(1 + \alpha \dt{}  \lambda_2)}{(1 - (1 - \alpha)\dt{} \lambda_1)(1 - (1-\alpha) \dt{} \lambda_2)} \Big| \Big\}\] & $\alpha \leq 0.5$\\
Strong Gauss-Seidel &\[\max\Big\{ 1, \Big|\frac{1 }{(1 - \dt{} \lambda_1)(1 - \dt{} \lambda_2)} \Big| \Big\}\]  &$\forall \alpha$\\
\hline
\end{tabular}
\caption{The iterative matrix spectrum radius and unconditional stability range of the 1st-order  IMEX based on  weak/strong predictors and Jacobi/Gauss-Seidel strategies}
\label{TAB: MODEL-2WAY-SIMPLE O1 IMEX STABILITY}
\end{table}

We close this section by demonstrating that the unconditional linear stability
result requires analyzing the chosen numerical scheme applied to the entire
coupled system, rather than using the exact solution of one subsystem to reduce to
the problem to a single system with added-mass \cite{causin2005added}.
The exact solution of the model problem \eqnref{EQ: MODEL-2WAY-SIMPLE} with
the initial condition
\begin{equation}
 u^2(0) = \frac{\lambda_2}{\lambda_1}u^1(0)
\end{equation}
takes the form
\begin{equation} \label{EQ: EXACT SOL}
 u^2(t) = \frac{\lambda_2}{\lambda_1}u^1(t).
\end{equation}
From \tabref{TAB: MODEL-2WAY-SIMPLE PARTITION}, the partitioned scheme
resulting from the 1st-order, 2-stage IMEX-RK scheme with strong Gauss-Seidel
predictor leads to the one-step update equation
\begin{subequations}
 \begin{align}
  u^1_n &= u^1_{n-1} + \dt{}\lambda_1(u^1_{n} +  u^2_{n-1})
                                                        \label{EQ: DECOUPLE1} \\
  u^2_n &= u^2_{n-1} + \dt{}\lambda_2(u^1_{n} +  u^2_{n}),
 \end{align}
\end{subequations}
which was determined to be unconditionally stable; see
\tabref{TAB: MODEL-2WAY-SIMPLE O1 IMEX STABILITY}.
However, this scheme can be analyzed solely in terms of the equation for
$u^1_n$ by substituting the exact solution for $u^2$ into
\eqnref{EQ: DECOUPLE1}
\begin{equation}
u^1_n = u^1_{n-1} + \dt{}\lambda_1u^1_{n} + \dt{}\lambda_2u^1_{n-1}.
\end{equation}
In this case, the scheme is only conditionally stable
\begin{equation}
\left|\frac{1 + \dt{}\lambda_2}{1 - \dt{} \lambda_1}\right| \leq 1,
\end{equation}
which illustrates how the artificial decoupling brought into the analysis
through the exact solution can degenerate stability.

\section{Application to a coupled system of ordinary differential equations}
\label{SEC: APP ODE}
In this section, we study the proposed high-order partitioned solvers
and predictors on a $3 \times 3$ system of linear Ordinary Differential
Equations (ODEs)
\begin{equation} \label{EQ: ODE SYS}
\dot\stvc = \Aboldcal \stvc\,, \qquad
\Aboldcal = 
\begin{bmatrix}
 1 & 1 & 1 \\
 1 & 1 & 0 \\
 1 & 1 & 1
\end{bmatrix}\,, \qquad
\ubm = \begin{bmatrix} \stvc[1] \\ \stvc[2] \\ \stvc[3] \end{bmatrix}
\end{equation}
with initial condition $\stvc(0) = (1,\,0,\,2)^{T}$ and consider the time
domain $t \in (0,\,2]$. The exact solution at
any time $t$ is given in terms of the initial condition and the eigenvalue
decomposition of the coefficient matrix,
$\Aboldcal\Pboldcal = \Pboldcal\Sigmabold$, as
\begin{equation}
\stvc(t) = \Pboldcal  e^{t \Sigmabold} \Pboldcal^{-1} \stvc(0).
\end{equation}
To conform to the multiphysics formulation in (\ref{EQ: GOVERN-SD-INDIV})
the ODE system is treated as a coupled system with three subsystems. The
mass matrix is identity, the velocity term is taken as
\begin{equation}
\res = (\stvc[1] + \cpl[1], \stvc[2] + \cpl[2], \stvc[3] + \cpl[3])^T,
\end{equation}
and the coupling terms are defined as
\begin{equation}
 \cpl[1] = \stvc[2] + \stvc[3]\,,\quad
 \cpl[2] = \stvc[1]\,,\quad
 \cpl[3] = \stvc[1] + \stvc[2].
\end{equation}
This decomposition of the velocity is non-unique. In fact,
many other choices exist that will lead to different schemes; however, the
above choice is the most sensible since it mimics the multiphysics applications
we are targeting and possesses special structure. In particular,
$\cpl[i]$ does not depend on $\stvc[i]$ and therefore the
strong predictors and weak predictors are equivalent. Additionally, the
coupling term possesses the same structure as \eqnref{EQ: SPE STRU}, which
implies only $\cplprd[1]$ is required for the strong Gauss-Seidel predictor.

To validate the temporal convergence of the high-order partitioned scheme, we
apply the 2nd-order 2-stage trapezoidal rule, 3rd-order 4-stage ARK3(2)4L[2]SA,
and 4th-order 6 stage ARK4(3)6L[2]SA \cite{christopher2001additive} to the
ODE system in (\ref{EQ: ODE SYS}). These schemes will be abbreviated by
IMEX2, IMEX3, and IMEX4, respectively. For a given IMEX-RK scheme, three
different predictors are tested. The first two are the Jacobi and Gauss-Seidel
predictors proposed in Section~\ref{SEC: 4 PREDICTORS}; recall there is
no distinction between weak and strong predictors given the structure of
the coupling term. We also consider the predictor proposed in
\cite{froehle2014high, van2007higher} for two-field fluid-structure interaction problems,
which is a Gauss-Seidel-type predictor
\begin{equation}
 \cplprd_{n,j}^{1} =
 \sum_{k=1}^{j-1}\frac{\hat{a}_{jk} - a_{jk}}{\hat{a}_{jj}}\cpl_{n,k}^{1}
\end{equation}
for stages $j = 2,\,\dots,\,s$, where
$\cpl_{n,k}^1 = \pstgu[2]{n}{k} + \pstgu[3]{n}{k}$, and $\hat{a}_{ki}$ and
$a_{ki}$ are the Butcher coefficients of the ERK and ESDIRK schemes,
respectively, in \tabref{TAB: IMEX}. Unlike the predictors proposed in this
work, this predictor is stage-dependent and does not have an interpretation
at the ODE level. Therefore it is not guaranteed to preserve the design
order of the IMEX scheme, even though it does so empirically in
\cite{van2007higher, froehle2014high}.

The accuracy is quantified via the $L_\infty$-norm of the error in the numerical
solution at time $t = 2.0$
\begin{equation}
 e_\text{ODE} = \max_{1\leq i \leq3} |\pstp[i]{N} - \stvc[i](2)|,
\end{equation}
where $\stvc[i](2)$ is the exact solution at $t = 2.0$ and $\pstp[i]{N}$ the
numerical solution at the final time step for the $i$th subsystem. The
error $e_\text{ODE}$ as a function of the time step size for the second, third,
and fourth order IMEX-RK methods are shown in \figref{FIG: ODE 3FIELD}. Note
that in \figref{FIG: ODE 3FIELD JAC} and \figref{FIG: ODE 3FIELD GS}, the
schemes exhibit convergence at the design rate of the IMEX scheme and the error
with the Gauss-Seidel predictor is several times smaller than
that of the Jacobi predictor due to different error constants.
However, the stage-variant predictor in \figref{FIG: ODE 3FIELD OLD} results
in a scheme with an order of accuracy one less than the design order.
\ifbool{fastcompile}{}{
\begin{figure}
 \begin{subfigure}[b]{0.33\textwidth}
  \centering
  \begin{tikzpicture}
\begin{loglogaxis}[
    width=0.98\textwidth,
    height=0.98\textwidth,
    xlabel={Time step ($\Delta t$)},
    ylabel={$e_\text{ODE}$}]

\addplot [red, solid, thick, mark=square*, mark size=1, mark options={solid}]  table[x index=0, y index=4] {data/ode3/ode_3fields_pred_Jac_err_RK2.dat};
\addplot [green!75!black, solid, thick, mark=triangle*, mark size=1, mark options={solid}]  table[x index=0, y index=4] {data/ode3/ode_3fields_pred_Jac_err_RK3.dat};
\addplot [blue, solid, thick, mark=diamond*, mark size=1, mark options={solid}]  table[x index=0, y index=4] {data/ode3/ode_3fields_pred_Jac_err_RK4.dat};

\logLogSlopeTriangle{0.4}{0.1}{0.72}{2}{red};
\logLogSlopeTriangle{0.4}{0.1}{0.42}{3}{green!75!black};
\logLogSlopeTriangle{0.4}{0.1}{0.2}{4}{blue};

\end{loglogaxis}

\end{tikzpicture}
  \caption{Jacobi predictor}
  \label{FIG: ODE 3FIELD JAC}
 \end{subfigure}
 \begin{subfigure}[b]{0.33\textwidth}
  \centering
  \begin{tikzpicture}
\begin{loglogaxis}[
    width=0.98\textwidth,
    height=0.98\textwidth,
    xlabel={Time step ($\Delta t$)},
    ylabel={~~}]

\addplot [red, solid, thick, mark=square*, mark size=1, mark options={solid}]  table[x index=0, y index=4] {data/ode3/ode_3fields_pred_GS_err_RK2.dat};\label{line:ode_3field:rk2}
\addplot [green!75!black, solid, thick, mark=triangle*, mark size=1, mark options={solid}]  table[x index=0, y index=4] {data/ode3/ode_3fields_pred_GS_err_RK3.dat};\label{line:ode_3field:rk3}
\addplot [blue, solid, thick, mark=diamond*, mark size=1, mark options={solid}]  table[x index=0, y index=4] {data/ode3/ode_3fields_pred_GS_err_RK4.dat};\label{line:ode_3field:rk4}

\logLogSlopeTriangle{0.4}{0.1}{0.7}{2}{red};
\logLogSlopeTriangle{0.4}{0.1}{0.44}{3}{green!75!black};
\logLogSlopeTriangle{0.4}{0.1}{0.2}{4}{blue};

\end{loglogaxis}

\end{tikzpicture}
  \caption{Gauss-Seidel predictor}
  \label{FIG: ODE 3FIELD GS}
 \end{subfigure}
 \begin{subfigure}[b]{0.33\textwidth}
  \centering
  \begin{tikzpicture}
\begin{loglogaxis}[
    width=0.98\textwidth,
    height=0.98\textwidth,
    xlabel={Time step ($\Delta t$)},
    ylabel={~~}]

\addplot [red, solid, thick, mark=square*, mark size=1, mark options={solid}]  table[x index=0, y index=4] {data/ode3/ode_3fields_pred_old_err_RK2.dat};
\addplot [green!75!black, solid, thick, mark=triangle*, mark size=1, mark options={solid}]  table[x index=0, y index=4] {data/ode3/ode_3fields_pred_old_err_RK3.dat};
\addplot [blue, solid, thick, mark=diamond*, mark size=1, mark options={solid}]  table[x index=0, y index=4] {data/ode3/ode_3fields_pred_old_err_RK4.dat};

\logLogSlopeTriangle{0.4}{0.1}{0.78}{1}{red};
\logLogSlopeTriangle{0.4}{0.1}{0.51}{2}{green!75!black};
\logLogSlopeTriangle{0.4}{0.1}{0.2}{3}{blue};

\end{loglogaxis}

\end{tikzpicture}
  \caption{Stage-variant predictor \cite{van2007higher}}
  \label{FIG: ODE 3FIELD OLD}
 \end{subfigure}
 \caption{Convergence of the IMEX2 (\ref{line:ode_3field:rk2}),
          IMEX3 (\ref{line:ode_3field:rk3}), and
          IMEX4 (\ref{line:ode_3field:rk4}) schemes with various
          predictors as applied to the ODE system.}
 \label{FIG: ODE 3FIELD}
\end{figure}
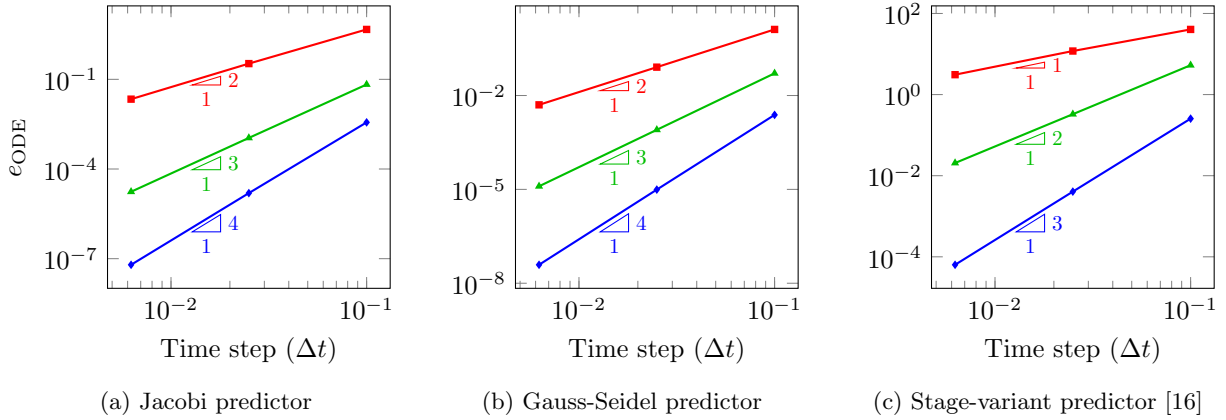
}

\section{Application to time-dependent advection-diffusion-reaction equations}
\label{SEC: APP ADR}
In this section, we consider time-dependent coupled advection-diffusion-reaction
(ADR) systems that have applications in the modeling of chemical reactions
\cite{tezduyar1986discontinuity}, the description for superconductivity of
liquids \cite{estep2000estimating}, and biological predator-prey models
\cite{estep2000using}. The governing equation for the $i$th species in a
general ADR system with $n$ components in $d$-dimensions is
\begin{equation} \label{EQ: CLAW ADR}
  \pder{u^i}{t} + (v^i \cdot \nabla) u^i - \nabla( D^i \cdot \nabla u^i) =
  f^i(u,\,x,\,t), \quad (x,\,t) \in v\times (0,\,1],
  \quad 1 \leq i \leq n.
\end{equation}
Here, $u = \begin{bmatrix} u^1 & \cdots & u^n \end{bmatrix}^T$ contains the
$n$ conserved quantities modeled by the ADR equations, $\Omega \subset \Rbb^d$
is the computational domain, $D^i \in \Rbb^{d \times d}$ is the diffusivity
matrix and $v^i(x) \in \Rbb^d$ is the velocity field for the $i$th species.
In this work, we consider the predator-prey model from \cite{estep2000using},
which involves $n = 2$ coupled systems with
\begin{equation}
 \begin{aligned}
  f^1(u,\,x,\,t) = u^1(-(u^1 - a^1)(u^1 - 1) - a^2u^2)\,\qquad
  f^2(u,\,x,\,t) = u^2(-a^3 - a^4u^2 + a^2 u^1)
  \label{EQ: PREDATOR PREY SOURCE}
 \end{aligned}
\end{equation}
where $a^1 = 0.25$, $a^2 = 2$, $a^3 = 1$, $a^4 = 3.4$, and the diffusivity
matrices are constant, isotropic $D^1 = D^2 = 0.01\Ibm_2$ and
$\Ibm_2$ is the $2 \times 2$ identity matrix. The computational domain is the
two-dimensional unit square $\Omega = [-0.5,\,0.5]\times[-0.5,\,0.5]$ with the
prey initially uniformly distributed, and predators initially gathered
near $(x_0,y_0) = (-0.25,-0.25)$ 
\begin{equation}
 \begin{aligned}
  u^1(x,\,y,\, 0) = 1.0   \quad \textrm{and} \quad  u^2(x,\,y,\,0) =
                \begin{cases}
                  0 & r > d\\
                  e^{-\frac{d^2}{d^2 - r^2}} & r \leq d\\
                \end{cases},
  \end{aligned}
\end{equation}
where $d=0.2$, $r = \sqrt{(x - x_0)^2 + (y-y_0)^2}$. The boundary conditions
are all Neumann conditions $\displaystyle{\pder{u}{n} = 0}$ and the velocity
fields are constant $v^1(x) = (0,0)$ and $v_2(x) = (0.5,0.5)$.
The equations are discretized with a standard high-order discontinuous
Galerkin method using Roe's upwind flux \cite{roe1981approximate} for the
inviscid numerical flux and the Compact DG flux \cite{peraire2008compact}
for the viscous numerical flux on a $40 \times 40$ structured mesh of
quadratic simplex elements.

The governing equations in (\ref{EQ: CLAW ADR}) reduce to the following
system of ODEs after the DG discretization is applied
\begin{equation}
\label{EQ: ADR0}
 \mass[i]\stvcdot[i] = \res[i](\stvc[i]) +
                       \cpl[i](\stvc[1],\,\stvc[2]),
\end{equation}
where $\mass[i]$ is the fixed mass matrix, $\stvc[i](t)$ is the semi-discrete
state vector, i.e., the discretization of $u$ on $\Omega$,
$\res[i](\stvc[i])$ is the spatial discretization of the advection and
diffusion terms on $\Omega$, and $\cpl[i]$ is the coupling term
that contains the DG discretization of the $i$th reaction source term in
(\ref{EQ: PREDATOR PREY SOURCE}). This non-unique decomposition of the
governing equation (\ref{EQ: CLAW ADR}) implies that,
once the high-order partitioned solver is applied, various terms of the
reaction source term will be predicted. The solution of (\ref{EQ: ADR0})
using the IMEX4 scheme with strong Gauss-Seidel predictor is provided in
\figref{FIG: PREDATOR_PREY} using the time step size $\dt{}=0.1$.
The predators are diffused quickly and migrate diagonally
upward, while the prey are mostly affected by the coupled reaction near the
extent of the predator population.

\ifbool{fastcompile}{}{
\begin{figure}[ht]
\centering
 \includegraphics[width=0.3\textwidth]{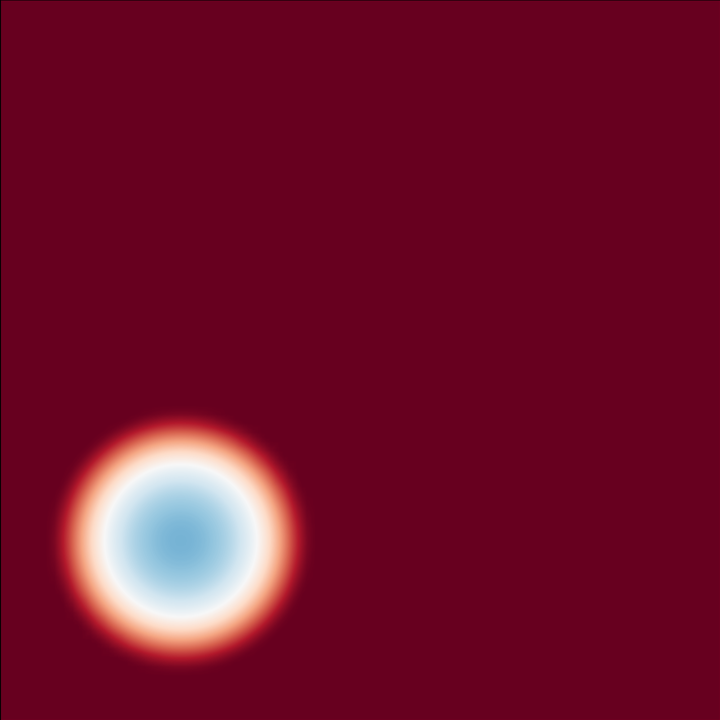} \quad
 \includegraphics[width=0.3\textwidth]{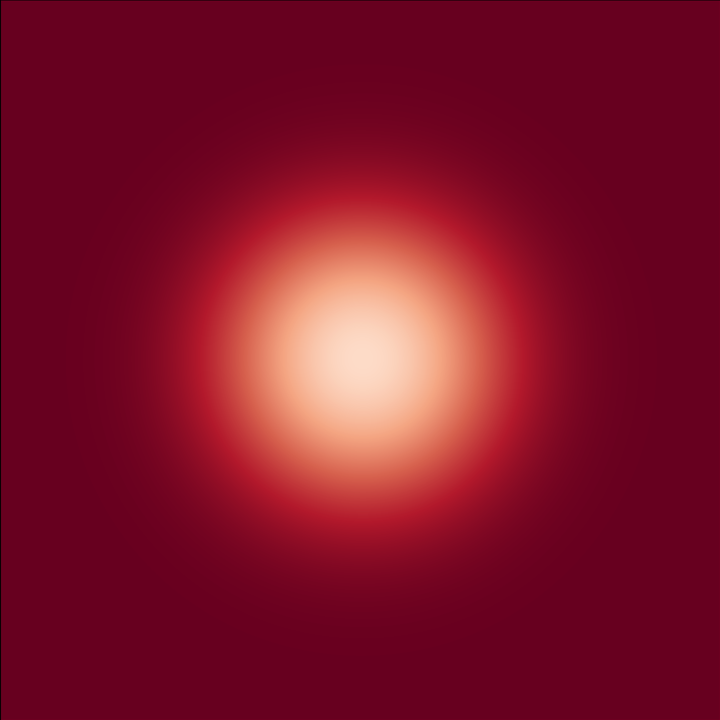} \quad
 \includegraphics[width=0.3\textwidth]{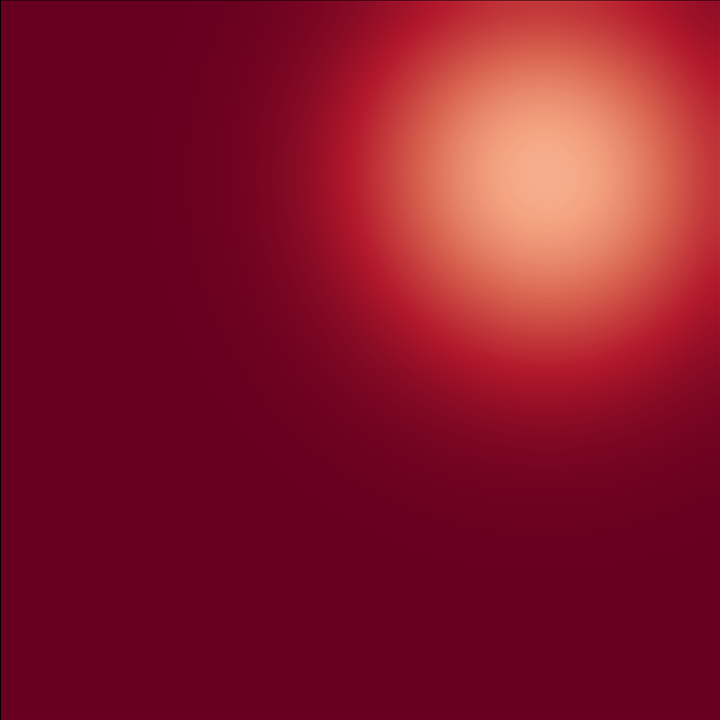} \\\vspace{2mm}
 \includegraphics[width=0.3\textwidth]{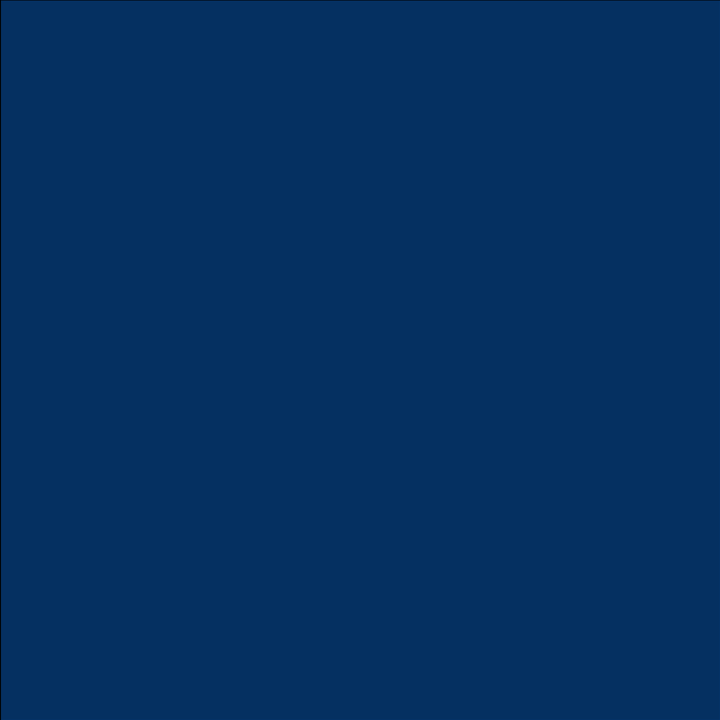} \quad
 \includegraphics[width=0.3\textwidth]{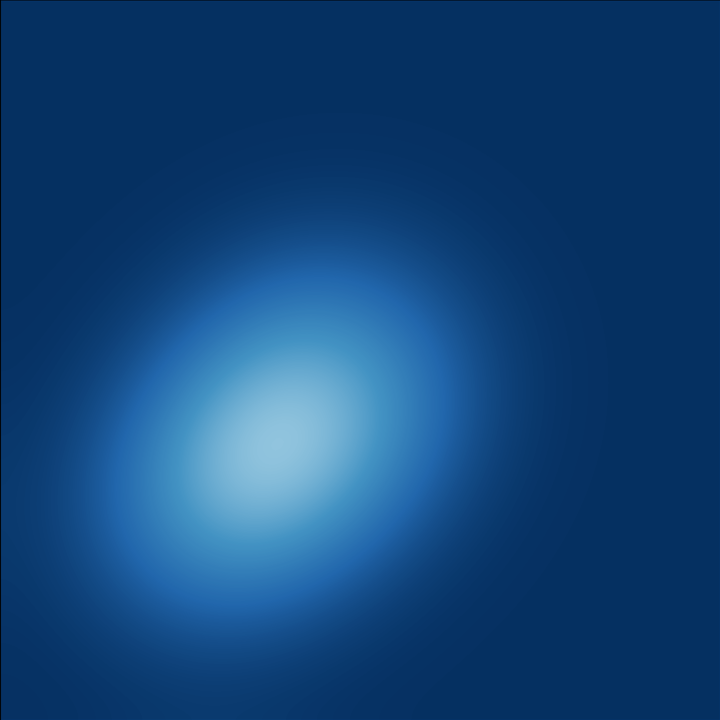} \quad
 \includegraphics[width=0.3\textwidth]{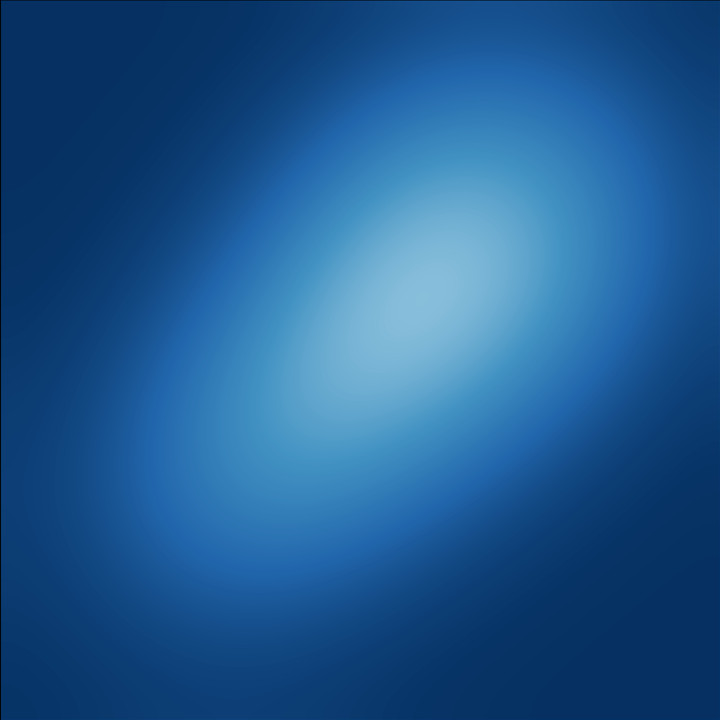}
 \caption{Predator (\emph{top}) and prey (\emph{bottom}) populations at
          various snapshots in time: $t = 0.0$ (\emph{left}), $t = 0.5$
          (\emph{center}), and $t = 1.0$ (\emph{right}).}
 \label{FIG: PREDATOR_PREY}
\end{figure}
}

To validate the temporal convergence of the high-order partitioned scheme, we
apply the 2nd-order 2-stage trapezoidal rule, 3rd-order 4-stage ARK3(2)4L[2]SA,
and 4th-order 6 stage ARK4(3)6L[2]SA \cite{christopher2001additive} to the
ADR system in (\ref{EQ: ADR0}). Similar to the previous section, these schemes
will be abbreviated by IMEX2, IMEX3, and IMEX4, respectively. For a
given IMEX-RK scheme, the four predictors proposed in
Section~\ref{SEC: 4 PREDICTORS} are tested. Similar to the previous section,
we use the $L_{\infty}$-error between a reference solution and the numerical
solution provided by a particular solver at an instant in time $t = 1.0$
to quantify the error
\begin{equation}
 e_\text{ADR} = \norm{\stvc[1](1.0) - \pstp[1]{N}}_{\infty},
\end{equation}
where $\stvc[1](1.0)$ is a reference solution of the first subsystem
at $t = 1.0$ obtained by using the IMEX4 scheme with
$\dt{} = 6.25 \times 10^{-3}$ and the strong Gauss-Seidel predictor and
$\pstp[i]{N}$ is the numerical solution at the final time step for the first
subsystem. The error $e_\text{ADR}$ as a function of time step size for
the second, third, and fourth order IMEX-RK methods are provided in
\figref{FIG: PREDATOR_PREY_ERR}.
From this figure we see the design order of accuracy of the scheme is
obtained for all four proposed predictors. Unlike the ODE system in the
previous system, there is not a significant difference between the accuracy
at a given time step between the Jacobi and Gauss-Seidel predictors. This
figure also shows that no stability issues were observed for any of the
predictors, even for the coarsest time step $\Delta t = 0.1$.
\ifbool{fastcompile}{}{
\begin{figure}[ht]
 \begin{subfigure}[b]{0.245\textwidth}
  \centering
  \begin{tikzpicture}
\begin{loglogaxis}[
    width=1.05\textwidth,
    height=1.05\textwidth,
    xlabel={Time step ($\Delta t$)},
    ylabel={$e_\text{ADR}$}]

\addplot [red, solid, thick, mark=square*, mark size=1, mark options={solid}]  table[x index=0, y index=1] {data/predator_prey/predator_prey_err_RK2_weak_Jac.dat};\label{line:predator_prey_weak_Jac:rk2}
\addplot [green!75!black, solid, thick, mark=triangle*, mark size=1, mark options={solid}]  table[x index=0, y index=1] {data/predator_prey/predator_prey_err_RK3_weak_Jac.dat};\label{line:predator_prey_weak_Jac:rk3}
\addplot [blue, solid, thick, mark=diamond*, mark size=1, mark options={solid}]  table[x index=0, y index=1] {data/predator_prey/predator_prey_err_RK4_weak_Jac.dat};\label{line:predator_prey_weak_Jac:rk4}

\logLogSlopeTriangle{0.25}{0.1}{0.48}{2}{red};
\logLogSlopeTriangle{0.37}{0.1}{0.28}{3}{green!75!black};
\logLogSlopeTriangle{0.6}{0.1}{0.15}{4}{blue};

\end{loglogaxis}

\end{tikzpicture}
  \caption{weak Jacobi}
 \end{subfigure}
 \begin{subfigure}[b]{0.245\textwidth}
  \centering
  \begin{tikzpicture}
\begin{loglogaxis}[
    width=1.05\textwidth,
    height=1.05\textwidth,
    xlabel={Time step ($\Delta t$)},
    ylabel={~~}]

\addplot [red, solid, thick, mark=square*, mark size=1, mark options={solid}]  table[x index=0, y index=1] {data/predator_prey/predator_prey_err_RK2_strong_Jac.dat};\label{line:predator_prey_strong_Jac:rk2}
\addplot [green!75!black, solid, thick, mark=triangle*, mark size=1, mark options={solid}]  table[x index=0, y index=1] {data/predator_prey/predator_prey_err_RK3_strong_Jac.dat};\label{line:predator_prey_strong_Jac:rk3}
\addplot [blue, solid, thick, mark=diamond*, mark size=1, mark options={solid}]  table[x index=0, y index=1] {data/predator_prey/predator_prey_err_RK4_strong_Jac.dat};\label{line:predator_prey_strong_Jac:rk4}

\logLogSlopeTriangle{0.25}{0.1}{0.48}{2}{red};
\logLogSlopeTriangle{0.37}{0.1}{0.30}{3}{green!75!black};
\logLogSlopeTriangle{0.6}{0.1}{0.15}{4}{blue};

\end{loglogaxis}

\end{tikzpicture}
  \caption{strong Jacobi}
 \end{subfigure}
 \begin{subfigure}[b]{0.245\textwidth}
  \centering
  \begin{tikzpicture}
\begin{loglogaxis}[
    width=1.05\textwidth,
    height=1.05\textwidth,
    xlabel={Time step ($\Delta t$)},
    ylabel={~~}]

\addplot [red, solid, thick, mark=square*, mark size=1, mark options={solid}]  table[x index=0, y index=1] {data/predator_prey/predator_prey_err_RK2_weak_GS.dat};\label{line:predator_prey_weak_GS:rk2}
\addplot [green!75!black, solid, thick, mark=triangle*, mark size=1, mark options={solid}]  table[x index=0, y index=1] {data/predator_prey/predator_prey_err_RK3_weak_GS.dat};\label{line:predator_prey_weak_GS:rk3}
\addplot [blue, solid, thick, mark=diamond*, mark size=1, mark options={solid}]  table[x index=0, y index=1] {data/predator_prey/predator_prey_err_RK4_weak_GS.dat};\label{line:predator_prey_weak_GS:rk4}

\logLogSlopeTriangle{0.25}{0.1}{0.48}{2}{red};
\logLogSlopeTriangle{0.37}{0.1}{0.28}{3}{green!75!black};
\logLogSlopeTriangle{0.6}{0.1}{0.15}{4}{blue};

\end{loglogaxis}

\end{tikzpicture}
  \caption{weak Gauss-Seidel}
 \end{subfigure}
 \begin{subfigure}[b]{0.245\textwidth}
  \centering
  \begin{tikzpicture}
\begin{loglogaxis}[
    width=1.05\textwidth,
    height=1.05\textwidth,
    xlabel={Time step ($\Delta t$)},
    ylabel={~~}]

\addplot [red, solid, thick, mark=square*, mark size=1, mark options={solid}]  table[x index=0, y index=1] {data/predator_prey/predator_prey_err_RK2_strong_GS.dat};\label{line:predator_prey_strong_GS:rk2}
\addplot [green!75!black, solid, thick, mark=triangle*, mark size=1, mark options={solid}]  table[x index=0, y index=1] {data/predator_prey/predator_prey_err_RK3_strong_GS.dat};\label{line:predator_prey_strong_GS:rk3}
\addplot [blue, solid, thick, mark=diamond*, mark size=1, mark options={solid}]  table[x index=0, y index=1] {data/predator_prey/predator_prey_err_RK4_strong_GS.dat};\label{line:predator_prey_strong_GS:rk4}

\logLogSlopeTriangle{0.25}{0.1}{0.48}{2}{red};
\logLogSlopeTriangle{0.37}{0.1}{0.30}{3}{green!75!black};
\logLogSlopeTriangle{0.6}{0.1}{0.15}{4}{blue};

\end{loglogaxis}

\end{tikzpicture}
  \caption{strong Gauss-Seidel}
 \end{subfigure}
 \caption{Convergence of the IMEX2 (\ref{line:predator_prey_weak_GS:rk2}),
          IMEX3 (\ref{line:predator_prey_weak_GS:rk3}), and
          IMEX4(\ref{line:predator_prey_weak_GS:rk4}) schemes with various
          predictors applied to the predator-prey ADR system. For this problem,
          all predictors achieve the design order of the IMEX schemes, with
          small differences in the accuracy.}
 \label{FIG: PREDATOR_PREY_ERR}
\end{figure}
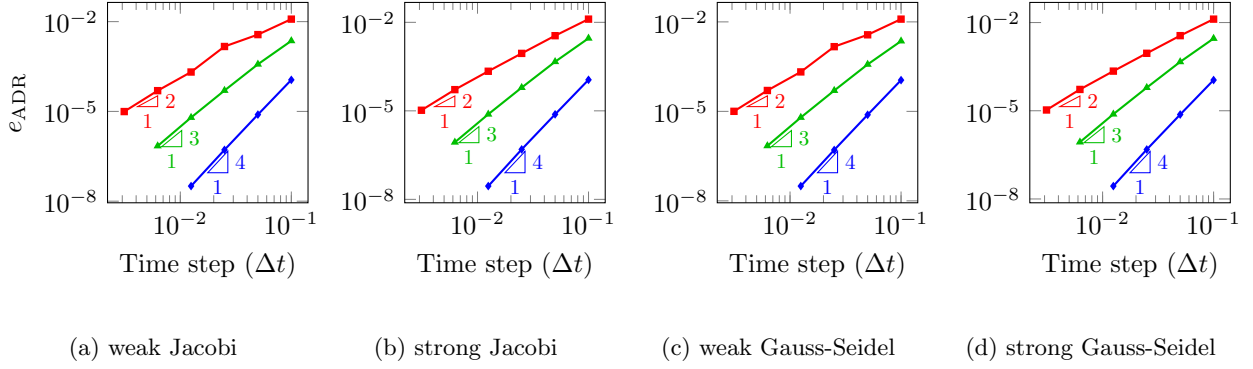
}

\section{Application to fluid-structure interaction}
\label{SEC: APP FSI}
In this section, we demonstrate the proposed high-order IMEX-based partitioned
solver on fluid-structure interaction (FSI) problems. Partitioned solution
procedures are widely used to solve such problems given that they allow
maximal re-use of individual fluid and structure software. The FSI problem
is usually formulated in an Arbitrary Lagrangian-Eulerian framework using
three fields: the deformation of the solid, the fluid flow, and the motion
of the fluid mesh. The deformation of the fluid mesh is commonly assumed
quasi-static \cite{farhat1995mixed, farhat1998torsional} or interpolation, e.g., via radial basis functions,
is used to transfer the boundary displacement of the fluid mesh into the
interior \cite{van2007higher, froehle2014high}. In both cases, the
formulation effectively reduces to a two-field system involving the
structural displacements and the fluid flow. Both two- and three-field FSI
formulations are considered in this section.

\subsection{Governing equations and semi-discretization}
This section introduces the governing partial differential equations for
the 2- and 3-field FSI formulation as a coupled multiphysics system
(\ref{EQ: GOVERN}) and their semi-discretization to yield a system
of ODEs connected via coupling terms (\ref{EQ: GOVERN-SD-INDIV}).

\subsubsection{Compressible fluid flow}
\label{SUBSUBSEC: ALE Fluid}
The governing equations for compressible fluid flow, defined on a
deformable fluid domain $\Omega(t)$, can be written as a viscous conservation
law
\begin{equation}
\label{EQ: FSI GOVERN}
\frac{\partial U}{\partial t} + \nabla \cdot \Fcal^{inv}(U) + \nabla \cdot \Fcal^{vis}(U, \nabla U)= 0 \quad \text{in} \quad \Omega(t),
\end{equation}
where $U$ is the conservative state variable vector and the physical
flux consists of inviscid part $\Fcal^{inv}(U)$ and a viscous part
$\Fcal^{vis}(U,\,\nabla U)$. The conservation law in (\ref{EQ: FSI GOVERN})
is transformed to a fixed reference domain $\Omega_0$ by defining a
time-dependent diffeomorphism $\Gcal$ between the reference domain and
physical domain; see Figure~\ref{FIG: DOM MAP}. At each time $t$, a point
$X$ in the reference domain $\Omega_0$ is mapped to $x(X,t) = \Gcal(X,t)$
in the physical domain $\Omega(t)$.
\ifbool{fastcompile}{}{
\begin{figure}
  \centering
  \includegraphics[width=2.5in]{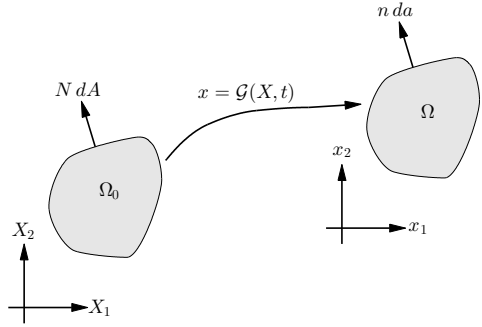}
  \caption{Mapping between reference and physical domains.}
  \label{FIG: DOM MAP}
\end{figure}
}
The deformation gradient $G$, velocity $v_G$, and Jacobian $g$ of
the mapping are defined as
\begin{equation}
 G = \nabla_X\Gcal\,,\quad
 v_G = \pder{\Gcal}{t}\,,\quad
 g = \det G.
\end{equation}
Following the procedure in \cite{persson2009discontinuous, zahr2016adjoint},
the governing equation
(\ref{EQ: FSI GOVERN}) can be written in the reference domain as
\begin{equation}\label{EQ: FSI GOVERN TRANSF}
\frac{\partial U_X}{\partial t} + \nabla_X \cdot \Fcal_X^{inv}(U_X) +
\nabla_X\cdot \Fcal_X^{vis}(U_X, \nabla_X U_X)= 0
\quad \text{in} \quad \Omega_0, 
\end{equation}
where $\nabla_X$ defines the spatial derivative with respect to the reference
domain, conserved quantities and its derivatives in the reference domain are
written as
\begin{equation}\label {EQ: FSI GOVERN TRANSF2}
 U_X = gU\,,\qquad
 \nabla_X U_X = g\nabla U_X \cdot G + g^{-1}U_X \frac{\partial g}{\partial X}.
\end{equation}
The inviscid and viscous fluxes are transformed to the reference domain as
\begin{equation}
\begin{aligned}
 \Fcal_X^{inv}(U_X) &= g\Fcal^{inv}(g^{-1} U_X)G^{-T} -
                       U_X \mathbin{\mathop{\otimes}} G^{-1}v_G, \\
 \Fcal_X^{vis}(U_X) &= g\Fcal^{vis}\left(g^{-1} U_X,
                       g^{-1}\left[\nabla_X U_X -
                       g^{-1} U_X \pder{g}{X}\right]G^{-1}\right)G^{-T}.
 \end{aligned}
\end{equation}

The governing equations in (\ref{EQ: FSI GOVERN TRANSF}) reduce to the following
system of ODEs after an appropriate spatial discretization, such as a
discontinuous Galerkin or finite volume method, is applied
\begin{equation}
\label{EQ: FSI3 FLUID COUPLING}
 \mass[f]\stvcdot[f] = \res[f](\stvc[f],\,\cpl[f]),
\end{equation}
where $\mass[f]$ is the fixed mass matrix, $\stvc[f](t)$ is the semi-discrete
fluid state vector, i.e., the discretization of $U_X$ on $\Omega_0$,
$\res[f](\stvc[f],\,\cpl[f])$ is the spatial discretization of the transformed
inviscid and viscous fluxes on $\Omega_0$, and $\cpl[f]$ is the coupling term
that contains information about the domain mapping $\Gcal(X,\,t)$. In
particular, the coupling term contains the position and velocities of the
nodal coordinates of the computational mesh.
The domain mapping is defined using an element-wise nodal (Lagrangian)
polynomial basis on the mesh with coefficients from the nodal positions and
velocities.

\subsubsection{Simple structure model}
\label{SEC: APP FSI STRUCT}
In general, the governing equations for the structure will be given by
a system of partial differential equation such as the continuum equations in
total Lagrangian form with an arbitrary constitutive law. However, in this
work, we only consider simple structures such as mass-spring-damper systems
that can directly be written as a second-order system of ODEs
\begin{equation} \label{EQ: SIMP STRUCT}
 m_s\ddot{u}_s + c_s\dot{u}_s + k_s u_s = f_{ext}(t),
\end{equation}
where $m_s$ is the mass of the (rigid) object, $c_s$ is the damper resistance
constant, $k_s$ is the spring stiffness, and $f_{ext}(t)$ is a time-dependent
external load, which will be given by integrating the pointwise force the
fluid exerts on the object. These simple structures allow us to study the
stability and accuracy properties of the proposed high-order partitioned
solver for this class of multiphysics problems without the distraction of
transferring solution fields across the fluid-structure interface.

To conform to the notation in this document and encapsulate the
semi-discretization of PDE-based structure models, the equations in
(\ref{EQ: SIMP STRUCT}) are re-written in a first-order form as
\begin{equation}
 \mass[s]\stvcdot[s] = \res[s](\stvc[s],\,\cpl[s]).
\end{equation}
In the case of the simple structure in (\ref{EQ: SIMP STRUCT}), the mass
matrix, state vector, residual, and coupling term are
\begin{equation}
 \mass[s] = \begin{bmatrix} m_s & \\ & 1 \end{bmatrix}, \qquad
 \stvc[s] = \begin{bmatrix} \dot{u}_s \\ u_s \end{bmatrix}, \qquad
 \cpl[s] = f_{ext}, \qquad
 \res[s](\stvc[s],\,\cpl[s]) =
              \begin{bmatrix} f_{ext}-c_s\dot{u}_s-k_s u_s \\ u_s \end{bmatrix}.
\label{EQ: FSI3 STRUCT COUPLING}
\end{equation}

\subsubsection{Deformation of the fluid domain}
In the three-field fluid-structure interaction formulation pioneered in
\cite{farhat1995mixed, farhat1998torsional}, the fluid mesh is considered a pseudo-structure driven solely
by Dirichlet boundary conditions provided by the displacement of the structure
at the fluid-structure interface. The governing equations are given by the
continuum mechanics equations in total Lagrangian form with an arbitrary
constitutive law
\begin{equation} 
 \begin{aligned}
  \pder{\bar{p}}{t} - \nabla \cdot P(G) &= 0
    \qquad\qquad &&\text{in}~\Omega_0 \\
  x &= x_b
    \qquad\qquad &&\text{on}~\partial \Omega_0^D \\
  \dot{x} &= \dot{x}_b
    \qquad\qquad &&\text{on}~\partial \Omega_0^D,
 \end{aligned}
\label{EQ: FSI3 MESH}
\end{equation}
where $\bar{p}(X,\,t) = \rho_m\dot{x}$ is the linear momentum, $\rho_m$
is the density, and $P$ is the first Piola-Kirchhoff stress of the
pseudo-structure. The deformation gradient $G$ is the mapping that
defines the deformation of the reference fluid domain $\Omega_0$ to physical
fluid domain $\Omega(t)$. The position and velocity of the fluid domain are
prescribed along $\partial \Omega_0^D$, the union of the fluid-structure
interface and the fluid domain boundary.

The governing equations in (\ref{EQ: FSI3 MESH}) reduce to the following
system of ODEs after an appropriate spatial discretization, such as the
finite element method, is applied and recast in first-order form
\begin{equation}
 \Mbm^x\dot\ubm^x = \rbm^x(\ubm^x,\,\cbm^x)
\end{equation}
where $\mass[x]$ is the fixed mass matrix, $\stvc[x](t)$ is the semi-discrete
state vector consisting of the displacements and velocities of the
mesh nodes, $\res[x](\stvc[x],\,\cpl[x])$ is the spatial discretization
of the continuum equations and boundary conditions on the reference domain
$\Omega_0$, and $\cpl[x]$ is the coupling term that contains information about
the motion of the fluid structure interface. This model of the mesh motion
leads to a three-field FSI formulation when coupled to the fluid and
structure equations.

Alternatively, the motion of the fluid mesh can be described through a
parametrized mapping such as radial basis functions \cite{van2007higher, rendall2008unified, froehle2014high} or
blending maps \cite{persson2009discontinuous}. That is, the domain mapping
$x = \Gcal(X,\,t)$ is given by an analytical function, parametrized by the
deformation and velocity of the fluid-structure interface, that can be
analytically differentiated to obtain the deformation gradient $G(X,\,t)$
and velocity $v_G(X,\,t)$. Since the fluid mesh motion is no longer included
in the system of time-dependent partial differential equations, this leads to
a two-field FSI formulation in terms of the fluid and structure states only.

\subsubsection{Two-field and three-field fluid-structure coupling}
In the three-field fluid-structure interaction setting
\begin{equation}
\label{EQ: FSI3 SEMI DISC}
\mass[s]\stvcdot[s] = \res[s](\stvc[s],\,\cpl[s]), \quad 
\mass[x]\stvcdot[x] = \res[x](\stvc[x],\,\cpl[x]) , \quad
\mass[f]\stvcdot[f] = \res[f](\stvc[f],\,\cpl[f])
\end{equation}
introduced in \cite{farhat1995mixed}, the coupling terms have the following
dependencies
\begin{equation} 
\cpl[s] = \cpl[s](\stvc[s],\,\stvc[x],\,\stvc[f]), \quad
\cpl[x] = \cpl[x](\stvc[s]), \quad
\cpl[f] = \cpl[f](\stvc[s],\, \stvc[x]).
\label{EQ: FSI3 COUPLE}
\end{equation}
From \eqnref{EQ: FSI3 STRUCT COUPLING}, the structure coupling term is the external force applied
to the structure that comes from integrating the fluid stresses over the
fluid-structure interface. The mesh coupling term is the position and
velocity of the fluid-structure interface and therefore depends solely on
the state of the structure. From Eq.~(\ref{EQ: FSI GOVERN TRANSF})-(\ref{EQ: FSI GOVERN TRANSF2}), 
the fluid coupling term is
the position and velocity of the entire fluid mesh and therefore depends
on the state of the structure and the mesh.

In the two-field FSI setting
\begin{equation}
\label{EQ: FSI2 SEMI DISC}
\mass[s]\stvcdot[s] = \res[s](\stvc[s],\,\cpl[s]), \quad 
\mass[f]\stvcdot[f] = \res[f](\stvc[f],\,\cpl[f])
\end{equation}
the mesh motion is given by an analytical function and the coupling terms have
the following dependencies
\begin{equation} \label{EQ: FSI2 COUPLE}
\cpl[s] = \cpl[s](\stvc[s],\,\stvc[f]), \quad
\cpl[f] = \cpl[f](\stvc[s]).
\end{equation}
In this case, the structure coupling term is determined from the fluid and
structure state since the external force depends on the traction integrated
over the fluid-structure interface. The fluid coupling term, i.e., the
position and velocity of the fluid mesh, is determined from the structure
state. Finally, the ordering of the subsystems implied in
(\ref{EQ: FSI3 SEMI DISC}) and (\ref{EQ: FSI2 SEMI DISC}) is used throughout
the remainder of this section, which plays an important role when
defining the Gauss-Seidel predictors.

\subsection{1D Fluid-structure-mesh three-field coupling piston problem}
We begin our investigation into the performance of the proposed high-order,
partitioned multiphysics solver in the FSI context with the canonical FSI
model problem: a one-dimensional piston (\figref{FIG: 1D PISTON}).
\ifbool{fastcompile}{}{
\begin{figure}
\centering
    \begin{tikzpicture}[scale=1.3]
      \fill [pattern = north east lines] (-0.3,0) rectangle (0.3,1.0);
      \fill[blue!15!white] (0,0) rectangle (4,1.0);
      \draw [fill=red!15!white,line width=0.5mm] (4,0) rectangle (4.5,1);
      \draw node at (1.5,0.5) {inviscid flow};
      
      \draw[decoration={aspect=0.3, segment length=1.5mm, amplitude=3mm,coil},decorate] (4.5,0.5) -- (5.5,0.5);
      \fill [pattern = north east lines] (5.5,0) rectangle (6.0,1.0);
      \draw[line width=0.5mm](0, 0)-- (0, 1);
      \draw[line width=0.5mm](5.5, 0)-- (5.5, 1);
      \draw[thick,->,line width=0.5mm,](0,-0.3) -- (4.5,-0.3) node[anchor=north west][scale=1.2] {x};
      \draw node at (4.25,0.5) {$m$};
      \draw node at (5.0,1) {$k$};
    \end{tikzpicture}
    \caption{One-dimensional piston system}
    \label{FIG: 1D PISTON}
\end{figure}
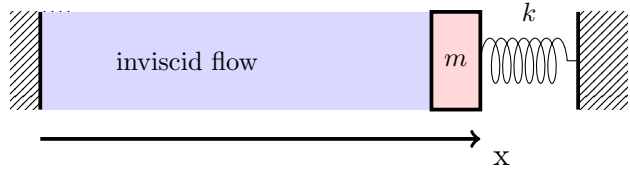
}
The inviscid fluid is governed by the one-dimensional Euler equations
\begin{equation}
 \begin{aligned}
  \pder{\rho}{t} + \pder{\rho u}{x} &= 0 \\
  \pder{\rho u}{t} + \pder{}{x}(\rho u^2 + p)  &= 0 \\
  \pder{\rho E}{t} + \pder{}{x}(u(\rho E + p)) &= 0
 \end{aligned}
\end{equation}
for $x \in \Omega(t) = [0,\,1.0-u_s]$, where $u_s$ is the displacement of the
piston, $\rho$ is the fluid density, $u$ is the fluid velocity, $E$ is the total
energy, the pressure $p$ is given by the ideal gas law
\begin{equation}
 p = (\gamma-1) \rho (E - \frac{1}{2}u^2),
\end{equation}
and the adiabatic gas constant is $\gamma = 1.4$. The fluid is
initially at rest $u = 0$ with a density $\rho = 1.0$ and pressure $p = 0.4$.
After transformation to the reference domain $\Omega_0 = [0,\,1]$ following the
procedure in Section~\ref{SUBSUBSEC: ALE Fluid}, the equations are semi-discretized by a standard
first-order finite volume method using Roe's flux \cite{roe1981approximate}
with $128$ elements.

The deformation of the fluid mesh is handled by considering the fluid domain
to be a pseudo-structure governed by the continuum equations in
\eqnref{EQ: FSI3 MESH}, restricted to the one-dimensional case with a
linear, isotropic constitutive law and infinitesimal strains assumed
\begin{equation}
\rho_m\ddot{u}_x = E_m\frac{\partial^2 u_x}{\partial X^2} - c_m \dot{u}_x,
\end{equation}
where $u_x(X,\,t)$ is the mesh displacement vector defined over the reference
domain $X \in \Omega_0$ and the density, Young's modulus, and damping
coefficient are $\rho_m = 1.0$, $E_m = 1.0$, $c_m = 0.0$, respectively. The
governing equation for the mesh deformation is discretized in space using the
finite difference method.

Finally, the structure is modeled by a linear mass-spring system as
\eqnref{EQ: SIMP STRUCT} with piston mass $m_s = 1.0$, spring stiffness
$k_s = 1.0$, and no damper $c_s = 0$. The piston is initially displaced a
distance of $u_s = -0.3$. Once the piston is released, it immediately begins
to recede due to the combination of the spring being perturbed from its
equilibrium configuration and the flow pressure, which causes a $C^0$
rarefaction wave near the interface.

To validate the temporal convergence of the scheme, the proposed high-order
partitioned framework is applied to solve the three-field coupled FSI problem.
In this case, we only consider the weak and strong Gauss-Seidel predictors.
The accuracy of a given simulation is
quantified by considering the error in fluid, mesh, and structure states
between a reference solution and the numerical solution at time $t = 5.0$
\begin{equation}
 \begin{aligned}
  e_\text{FSI3}^f &= \norm{\pstp[f]{N} - \stvc[f](5.0)}_\infty \\
  e_\text{FSI3}^x &= \norm{\pstp[x]{N} - \stvc[x](5.0)}_\infty \\
  e_\text{FSI3}^s &= \norm{\pstp[s]{N} - \stvc[s](5.0)}_\infty,
 \end{aligned}
\end{equation}
where $\stvc[f](5.0)$, $\stvc[x](5.0)$, $\stvc[s](5.0)$ are the fluid,
mesh, and structure states, respectively, from the reference solution,
computed by using the IMEX4 scheme with $\dt{} = 9.765625\times 10^{-5}$ and
strong Gauss-Seidel predictor at $t = 5.0$ and
$\pstp[f]{N}$, $\pstp[x]{N}$, $\pstp[s]{N}$
are the corresponding states from the numerical solution at the final time
step. The convergence plots are provided in \figref{FIG: 1D PISTON ERR WEAK}
and \figref{FIG: 1D PISTON ERR STRONG} and show the partitioned solver with
both predictors attain the design order of accuracy of the IMEX-RK scheme,
despite the fact that the solution is not $C^1$ continuous due to the
rarefaction wave.
\ifbool{fastcompile}{}{
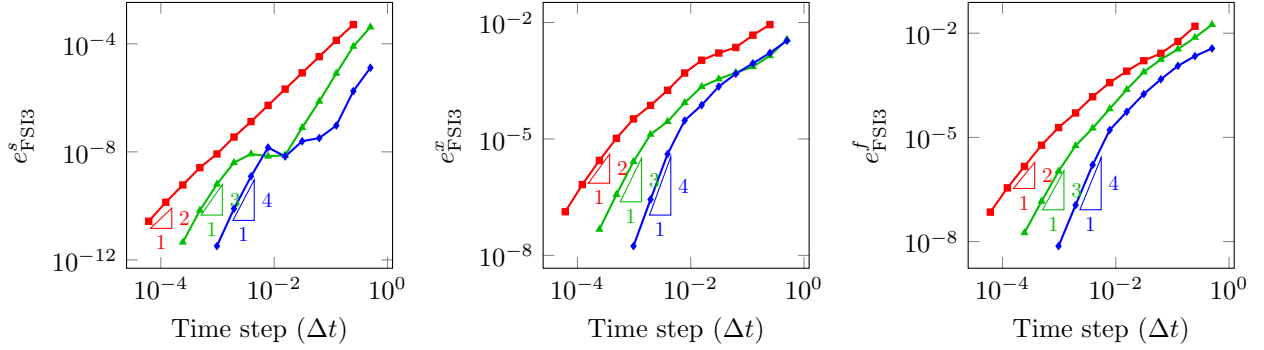
\begin{figure}
\begin{tikzpicture}
\begin{loglogaxis}[
    width=0.31\textwidth,
    height=0.31\textwidth,
    xlabel={Time step ($\Delta t$)},
    ylabel={$e_\text{FSI3}^s$}]

\addplot [red, solid, thick, mark=square*, mark size=1, mark options={solid}]  table[x index=0, y index=7] {data/piston_euler/euler_1d_3fields_err_RK2_weak.dat};\label{line:1d_fsi:rk2:structure:weak}
\addplot [green!75!black, solid, thick, mark=triangle*, mark size=1, mark options={solid}]  table[x index=0, y index=7] {data/piston_euler/euler_1d_3fields_err_RK3_weak.dat};\label{line:1d_fsi:rk3:structure:weak}
\addplot [blue, solid, thick, mark=diamond*, mark size=1, mark options={solid}]  table[x index=0, y index=7] {data/piston_euler/euler_1d_3fields_err_RK4_weak.dat};\label{line:1d_fsi:rk4:structure:weak}

\logLogSlopeTriangle{0.17}{0.08}{0.15}{2}{red};
\logLogSlopeTriangle{0.36}{0.08}{0.20}{3}{green!75!black};
\logLogSlopeTriangle{0.48}{0.08}{0.18}{4}{blue};

\end{loglogaxis}

\end{tikzpicture}
\begin{tikzpicture}
\begin{loglogaxis}[
    width=0.31\textwidth,
    height=0.31\textwidth,
    xlabel={Time step ($\Delta t$)},
    ylabel={$e_\text{FSI3}^x$}]

\addplot [red, solid, thick, mark=square*, mark size=1, mark options={solid}]  table[x index=0, y index=4] {data/piston_euler/euler_1d_3fields_err_RK2_weak.dat};\label{line:1d_fsi:rk2:mesh:weak}
\addplot [green!75!black, solid, thick, mark=triangle*, mark size=1, mark options={solid}]  table[x index=0, y index=4] {data/piston_euler/euler_1d_3fields_err_RK3_weak.dat};\label{line:1d_fsi:rk3:mesh:weak}
\addplot [blue, solid, thick, mark=diamond*, mark size=1, mark options={solid}]  table[x index=0, y index=4] {data/piston_euler/euler_1d_3fields_err_RK4_weak.dat};\label{line:1d_fsi:rk4:mesh:weak}

\logLogSlopeTriangle{0.25}{0.08}{0.32}{2}{red};
\logLogSlopeTriangle{0.37}{0.08}{0.25}{3}{green!75!black};
\logLogSlopeTriangle{0.48}{0.08}{0.20}{4}{blue};

\end{loglogaxis}

\end{tikzpicture}
\begin{tikzpicture}
\begin{loglogaxis}[
    width=0.31\textwidth,
    height=0.31\textwidth,
    xlabel={Time step ($\Delta t$)},
    ylabel={$e_\text{FSI3}^f$}]

\addplot [red, solid, thick, mark=square*, mark size=1, mark options={solid}]  table[x index=0, y index=1] {data/piston_euler/euler_1d_3fields_err_RK2_weak.dat};\label{line:1d_fsi:rk2:fluid:weak}
\addplot [green!75!black, solid, thick, mark=triangle*, mark size=1, mark options={solid}]  table[x index=0, y index=1] {data/piston_euler/euler_1d_3fields_err_RK3_weak.dat};\label{line:1d_fsi:rk3:fluid:weak}
\addplot [blue, solid, thick, mark=diamond*, mark size=1, mark options={solid}]  table[x index=0, y index=1] {data/piston_euler/euler_1d_3fields_err_RK4_weak.dat};\label{line:1d_fsi:rk4:fluid:weak}

\logLogSlopeTriangle{0.25}{0.08}{0.3}{2}{red};
\logLogSlopeTriangle{0.36}{0.08}{0.22}{3}{green!75!black};
\logLogSlopeTriangle{0.50}{0.08}{0.22}{4}{blue};

\end{loglogaxis}

\end{tikzpicture}
\caption{Convergence of the IMEX2 (\ref{line:1d_fsi:rk2:fluid:weak}), IMEX3 (\ref{line:1d_fsi:rk3:fluid:weak}), and IMEX4 (\ref{line:1d_fsi:rk4:fluid:weak}) with the \textbf{weak} Gauss-Seidel predictor as applied to the three-field coupling piston problem.}
\label{FIG: 1D PISTON ERR WEAK}
\end{figure}
}
\ifbool{fastcompile}{}{
 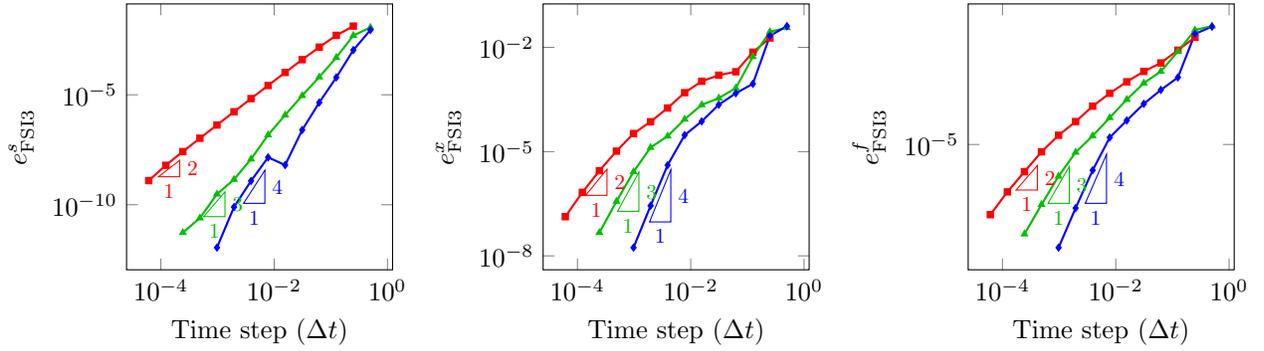
\begin{figure}
 \begin{tikzpicture}
\begin{loglogaxis}[
    width=0.31\textwidth,
    height=0.31\textwidth,
    xlabel={Time step ($\Delta t$)},
    ylabel={$e_\text{FSI3}^s$}]

\addplot [red, solid, thick, mark=square*, mark size=1, mark options={solid}]  table[x index=0, y index=7] {data/piston_euler/euler_1d_3fields_err_RK2_strong.dat};\label{line:1d_fsi:rk2:structure:strong}
\addplot [green!75!black, solid, thick, mark=triangle*, mark size=1, mark options={solid}]  table[x index=0, y index=7] {data/piston_euler/euler_1d_3fields_err_RK3_strong.dat};\label{line:1d_fsi:rk3:structure:strong}
\addplot [blue, solid, thick, mark=diamond*, mark size=1, mark options={solid}]  table[x index=0, y index=7] {data/piston_euler/euler_1d_3fields_err_RK4_strong.dat};\label{line:1d_fsi:rk4:structure:strong}

\logLogSlopeTriangle{0.2}{0.08}{0.35}{2}{red};
\logLogSlopeTriangle{0.37}{0.08}{0.20}{3}{green!75!black};
\logLogSlopeTriangle{0.52}{0.08}{0.25}{4}{blue};

\end{loglogaxis}

\end{tikzpicture}
 \begin{tikzpicture}
\begin{loglogaxis}[
    width=0.31\textwidth,
    height=0.31\textwidth,
    xlabel={Time step ($\Delta t$)},
    ylabel={$e_\text{FSI3}^x$}]

\addplot [red, solid, thick, mark=square*, mark size=1, mark options={solid}]  table[x index=0, y index=4] {data/piston_euler/euler_1d_3fields_err_RK2_strong.dat};\label{line:1d_fsi:rk2:mesh:strong}
\addplot [green!75!black, solid, thick, mark=triangle*, mark size=1, mark options={solid}]  table[x index=0, y index=4] {data/piston_euler/euler_1d_3fields_err_RK3_strong.dat};\label{line:1d_fsi:rk3:mesh:strong}
\addplot [blue, solid, thick, mark=diamond*, mark size=1, mark options={solid}]  table[x index=0, y index=4] {data/piston_euler/euler_1d_3fields_err_RK4_strong.dat};\label{line:1d_fsi:rk4:mesh:strong}

\logLogSlopeTriangle{0.24}{0.08}{0.28}{2}{red};
\logLogSlopeTriangle{0.36}{0.08}{0.22}{3}{green!75!black};
\logLogSlopeTriangle{0.48}{0.08}{0.18}{4}{blue};

\end{loglogaxis}

\end{tikzpicture}
 \begin{tikzpicture}
\begin{loglogaxis}[
    width=0.31\textwidth,
    height=0.31\textwidth,
    xlabel={Time step ($\Delta t$)},
    ylabel={$e_\text{FSI3}^f$}]

\addplot [red, solid, thick, mark=square*, mark size=1, mark options={solid}]  table[x index=0, y index=1] {data/piston_euler/euler_1d_3fields_err_RK2_strong.dat};\label{line:1d_fsi:rk2:fluid:strong}
\addplot [green!75!black, solid, thick, mark=triangle*, mark size=1, mark options={solid}]  table[x index=0, y index=1] {data/piston_euler/euler_1d_3fields_err_RK3_strong.dat};\label{line:1d_fsi:rk3:fluid:strong}
\addplot [blue, solid, thick, mark=diamond*, mark size=1, mark options={solid}]  table[x index=0, y index=1] {data/piston_euler/euler_1d_3fields_err_RK4_strong.dat};\label{line:1d_fsi:rk4:fluid:strong}

\logLogSlopeTriangle{0.26}{0.08}{0.3}{2}{red};
\logLogSlopeTriangle{0.38}{0.08}{0.25}{3}{green!75!black};
\logLogSlopeTriangle{0.52}{0.08}{0.25}{4}{blue};

\end{loglogaxis}

\end{tikzpicture}
 \caption{Convergence of the IMEX2 (\ref{line:1d_fsi:rk2:fluid:strong}), IMEX3 (\ref{line:1d_fsi:rk3:fluid:strong}), and IMEX4 (\ref{line:1d_fsi:rk4:fluid:strong}) with the \textbf{strong} Gauss-Seidel predictor as applied to the three-field coupling piston problem.}
\label{FIG: 1D PISTON ERR STRONG}
\end{figure}
}

\subsection{2D Fluid-structure two-field coupling foil damper problem}
We continue our investigation into the performance of the proposed
high-order, partitioned multiphysics solvers on FSI problems with
a two-dimensional energy-harvesting model problem \cite{peng2009energy, zahr2016adjoint}
that uses a two-field FSI formulation. Consider the mass-damper system
in Figure~\ref{FIG: FOIL DAMPER} suspended in an isentropic, viscous
flow where the rotational motion is a prescribed periodic motion
$\theta(t) = \frac{\pi}{4}\cos(2\pi f t)$ with frequency $f = 0.2$ and
the vertical displacement is determined by balancing the forces exerted
on the airfoil by fluid and damper.
\ifbool{fastcompile}{}{
\begin{figure}
 \centering
 \begin{tikzpicture}
\begin{axis}[
axis equal,
axis lines=none,
width=8cm,
ymax=0.6,
xmax=0.85,
xmin=-0.85,
ymin=-0.6]
\addplot [white, forget plot]
coordinates {
( -0.85000000,  -0.60000000)
(  0.85000000,  -0.60000000)
(  0.85000000,   0.60000000)
( -0.85000000,   0.60000000)
( -0.85000000,  -0.60000000)};

\draw [smooth, ultra thin] plot coordinates {(axis cs:-1.000000, -0.580000) (axis cs:1.000000, -0.580000)
};

\draw [smooth, ultra thin, opacity=0.5, fill=black!30!white] plot coordinates {(axis cs:-0.385791, 0.103372) (axis cs:-0.371603, 0.117293) (axis cs:-0.360141, 0.121044) (axis cs:-0.349136, 0.123087) (axis cs:-0.338373, 0.124230) (axis cs:-0.327767, 0.124784) (axis cs:-0.317276, 0.124909) (axis cs:-0.306874, 0.124701) (axis cs:-0.296545, 0.124224) (axis cs:-0.286276, 0.123522) (axis cs:-0.276059, 0.122625) (axis cs:-0.265887, 0.121559) (axis cs:-0.255756, 0.120343) (axis cs:-0.245660, 0.118992) (axis cs:-0.235597, 0.117520) (axis cs:-0.225564, 0.115935) (axis cs:-0.215559, 0.114249) (axis cs:-0.205579, 0.112467) (axis cs:-0.195623, 0.110598) (axis cs:-0.185688, 0.108646) (axis cs:-0.175775, 0.106617) (axis cs:-0.165880, 0.104516) (axis cs:-0.156004, 0.102346) (axis cs:-0.146146, 0.100111) (axis cs:-0.136304, 0.097816) (axis cs:-0.126477, 0.095462) (axis cs:-0.116665, 0.093053) (axis cs:-0.106868, 0.090591) (axis cs:-0.097083, 0.088079) (axis cs:-0.087312, 0.085518) (axis cs:-0.077553, 0.082912) (axis cs:-0.067806, 0.080262) (axis cs:-0.058070, 0.077570) (axis cs:-0.048345, 0.074838) (axis cs:-0.038630, 0.072066) (axis cs:-0.028925, 0.069258) (axis cs:-0.019230, 0.066413) (axis cs:-0.009544, 0.063535) (axis cs:0.000133, 0.060623) (axis cs:0.009802, 0.057679) (axis cs:0.019462, 0.054704) (axis cs:0.029114, 0.051700) (axis cs:0.038759, 0.048667) (axis cs:0.048396, 0.045606) (axis cs:0.058026, 0.042519) (axis cs:0.067649, 0.039405) (axis cs:0.077266, 0.036267) (axis cs:0.086876, 0.033105) (axis cs:0.096480, 0.029920) (axis cs:0.106077, 0.026711) (axis cs:0.115669, 0.023481) (axis cs:0.125255, 0.020230) (axis cs:0.134836, 0.016958) (axis cs:0.144411, 0.013666) (axis cs:0.153981, 0.010355) (axis cs:0.163546, 0.007025) (axis cs:0.173106, 0.003676) (axis cs:0.182662, 0.000310) (axis cs:0.192212, -0.003074) (axis cs:0.201758, -0.006475) (axis cs:0.211300, -0.009893) (axis cs:0.220837, -0.013327) (axis cs:0.230370, -0.016776) (axis cs:0.239899, -0.020241) (axis cs:0.249424, -0.023722) (axis cs:0.258945, -0.027217) (axis cs:0.268462, -0.030727) (axis cs:0.277975, -0.034251) (axis cs:0.287484, -0.037789) (axis cs:0.296989, -0.041342) (axis cs:0.306491, -0.044908) (axis cs:0.315989, -0.048488) (axis cs:0.325484, -0.052082) (axis cs:0.334974, -0.055688) (axis cs:0.344462, -0.059309) (axis cs:0.353946, -0.062942) (axis cs:0.363426, -0.066589) (axis cs:0.372903, -0.070248) (axis cs:0.382376, -0.073921) (axis cs:0.391845, -0.077607) (axis cs:0.401312, -0.081306) (axis cs:0.410774, -0.085018) (axis cs:0.420234, -0.088743) (axis cs:0.429689, -0.092482) (axis cs:0.439141, -0.096234) (axis cs:0.448590, -0.099999) (axis cs:0.458034, -0.103778) (axis cs:0.467475, -0.107571) (axis cs:0.476913, -0.111377) (axis cs:0.486346, -0.115198) (axis cs:0.495776, -0.119033) (axis cs:0.505202, -0.122882) (axis cs:0.514624, -0.126747) (axis cs:0.524042, -0.130626) (axis cs:0.533456, -0.134520) (axis cs:0.542865, -0.138430) (axis cs:0.552271, -0.142355) (axis cs:0.561672, -0.146297) (axis cs:0.571069, -0.150255) (axis cs:0.579809, -0.156664) (axis cs:0.569688, -0.155410) (axis cs:0.559571, -0.154139) (axis cs:0.549458, -0.152852) (axis cs:0.539350, -0.151549) (axis cs:0.529246, -0.150230) (axis cs:0.519146, -0.148896) (axis cs:0.509051, -0.147546) (axis cs:0.498959, -0.146182) (axis cs:0.488871, -0.144802) (axis cs:0.478787, -0.143409) (axis cs:0.468707, -0.142001) (axis cs:0.458631, -0.140579) (axis cs:0.448558, -0.139143) (axis cs:0.438489, -0.137693) (axis cs:0.428424, -0.136230) (axis cs:0.418363, -0.134753) (axis cs:0.408304, -0.133263) (axis cs:0.398250, -0.131760) (axis cs:0.388199, -0.130243) (axis cs:0.378152, -0.128713) (axis cs:0.368108, -0.127170) (axis cs:0.358067, -0.125615) (axis cs:0.348030, -0.124046) (axis cs:0.337997, -0.122463) (axis cs:0.327967, -0.120868) (axis cs:0.317941, -0.119260) (axis cs:0.307918, -0.117638) (axis cs:0.297899, -0.116003) (axis cs:0.287883, -0.114354) (axis cs:0.277871, -0.112691) (axis cs:0.267863, -0.111015) (axis cs:0.257859, -0.109325) (axis cs:0.247858, -0.107621) (axis cs:0.237861, -0.105902) (axis cs:0.227868, -0.104169) (axis cs:0.217879, -0.102420) (axis cs:0.207895, -0.100657) (axis cs:0.197914, -0.098878) (axis cs:0.187938, -0.097083) (axis cs:0.177966, -0.095271) (axis cs:0.167998, -0.093444) (axis cs:0.158035, -0.091599) (axis cs:0.148077, -0.089737) (axis cs:0.138123, -0.087857) (axis cs:0.128174, -0.085958) (axis cs:0.118231, -0.084041) (axis cs:0.108292, -0.082104) (axis cs:0.098359, -0.080148) (axis cs:0.088432, -0.078170) (axis cs:0.078510, -0.076172) (axis cs:0.068594, -0.074151) (axis cs:0.058684, -0.072108) (axis cs:0.048781, -0.070041) (axis cs:0.038883, -0.067951) (axis cs:0.028993, -0.065835) (axis cs:0.019109, -0.063694) (axis cs:0.009233, -0.061526) (axis cs:-0.000636, -0.059330) (axis cs:-0.010497, -0.057106) (axis cs:-0.020351, -0.054852) (axis cs:-0.030196, -0.052567) (axis cs:-0.040033, -0.050251) (axis cs:-0.049860, -0.047901) (axis cs:-0.059679, -0.045516) (axis cs:-0.069488, -0.043096) (axis cs:-0.079286, -0.040639) (axis cs:-0.089075, -0.038143) (axis cs:-0.098853, -0.035606) (axis cs:-0.108619, -0.033028) (axis cs:-0.118374, -0.030405) (axis cs:-0.128116, -0.027737) (axis cs:-0.137845, -0.025020) (axis cs:-0.147561, -0.022253) (axis cs:-0.157263, -0.019434) (axis cs:-0.166950, -0.016559) (axis cs:-0.176622, -0.013626) (axis cs:-0.186277, -0.010632) (axis cs:-0.195914, -0.007573) (axis cs:-0.205534, -0.004446) (axis cs:-0.215134, -0.001246) (axis cs:-0.224713, 0.002031) (axis cs:-0.234270, 0.005390) (axis cs:-0.243804, 0.008837) (axis cs:-0.253312, 0.012379) (axis cs:-0.262793, 0.016024) (axis cs:-0.272244, 0.019780) (axis cs:-0.281662, 0.023658) (axis cs:-0.291044, 0.027670) (axis cs:-0.300386, 0.031833) (axis cs:-0.309683, 0.036165) (axis cs:-0.318928, 0.040691) (axis cs:-0.328112, 0.045442) (axis cs:-0.337224, 0.050464) (axis cs:-0.346246, 0.055818) (axis cs:-0.355154, 0.061600) (axis cs:-0.363904, 0.067971) (axis cs:-0.372413, 0.075244) (axis cs:-0.380465, 0.084223) (axis cs:-0.385791, 0.103372)
};

\draw [smooth, ultra thin] plot coordinates {(axis cs:-0.050000, -0.325000) (axis cs:0.050000, -0.325000)
};

\draw [smooth, ultra thin] plot coordinates {(axis cs:0.000000, -0.580000) (axis cs:0.000000, -0.325000)
};

\draw [smooth, ultra thin] plot coordinates {(axis cs:-0.060000, -0.250000) (axis cs:-0.060000, -0.400000)
};

\draw [smooth, ultra thin] plot coordinates {(axis cs:-0.060000, -0.250000) (axis cs:0.060000, -0.250000)
};

\draw [smooth, ultra thin] plot coordinates {(axis cs:0.060000, -0.250000) (axis cs:0.060000, -0.400000)
};

\draw [smooth, ultra thin] plot coordinates {(axis cs:0.000000, -0.250000) (axis cs:0.000000, 0.000000)
};

\node[]    at    (axis cs:0.1, -0.025) {$m_s$};
\node[]    at    (axis cs:0.135, -0.325) {$c_s$};
\draw [smooth, solid, gray, ultra thin] plot coordinates {(axis cs:0.000600, 0.000000) (axis cs:-0.799400, 0.000000)
};

\draw [smooth, solid, gray, ultra thin] plot coordinates {(axis cs:0.000580, -0.000155) (axis cs:-0.772161, 0.206900)
};

\node[]    at    (axis cs:-0.6, 0.075) {$\theta(t)$};
\node[]    at    (axis cs:-0.58, -0.325) {$u_s$};
\draw [smooth, solid, gray, <->, ultra thin] plot coordinates {(axis cs:-0.650000, -0.580000) (axis cs:-0.650000, 0.000000)
};

\end{axis}
\end{tikzpicture}
 \caption{Foil-damper system} \label{FIG: FOIL DAMPER}
\end{figure}
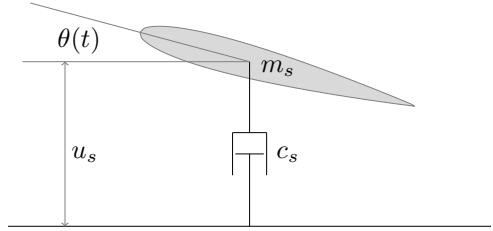
}
The governing equations for the fluid are the isentropic Navier-Stokes
equations:
\begin{align}
\frac{\partial \rho}{\partial t}  + \frac{\partial}{\partial x_i}
(\rho u_i) &= 0, \label{EQ: ns1}\\
\frac{\partial}{\partial t} (\rho u_i) +
\frac{\partial}{\partial x_i} (\rho u_i u_j+ p)  &=
+\frac{\partial \tau_{ij}}{\partial x_j}
\quad\text{for }i=1,2,3, \label{EQ: ns2} \\
\frac{\partial}{\partial t} (\rho E) +
\frac{\partial}{\partial x_i} \left(u_j(\rho E+p)\right) &=
-\frac{\partial q_j}{\partial x_j}
+\frac{\partial}{\partial x_j}(u_j\tau_{ij}), \label{EQ: ns3}
\end{align}
in $\Omega(t)$ where $\rho$ is the fluid density, $u_1,u_2,u_3$ are the
velocity components, and $E$ is the total energy. The viscous stress tensor and
heat flux are given by
\begin{align*}
\tau_{ij} = \mu
\left( \frac{\partial u_i}{\partial x_j} +
\frac{\partial u_j}{\partial x_i} -\frac23
\frac{\partial u_k}{\partial x_k} \delta_{ij} \right)
\qquad \text{ and } \qquad
q_j = -\frac{\mu}{\mathrm{Pr}} \frac{\partial}{\partial x_j}
\left( E+\frac{p}{\rho} -\frac12 u_k u_k \right).
\end{align*}
Here, $\mu$ is the viscosity coefficient and $\mathrm{Pr = 0.72}$ is
the Prandtl number which we assume to be constant. For an ideal gas,
the pressure $p$ has the form
\begin{align}
p=(\gamma-1)\rho \left( E - \frac12 u_k u_k\right), 
\label{EQ: ns5}
\end{align}
where $\gamma$ is the adiabatic gas constant. The isentropic assumption states
the entropy of the system is assumed constant, which is tantamount to the flow
being adiabatic and reversible. For a perfect gas, the entropy is defined as
\begin{equation}\label{eqn:entropy}
  s = p/\rho^\gamma.
\end{equation}
The conservation law defined in (\ref{EQ: ns1})-(\ref{EQ: ns3}) is
reformulated in an ALE framework, i.e., transformed to a reference domain
$\Omega_0$, as described in Section~\ref{SUBSUBSEC: ALE Fluid}. The transformed
conservation law is discretized with a standard high-order discontinuous
Galerkin method using Roe's flux \cite{roe1981approximate} for the inviscid
numerical flux and the Compact DG flux \cite{peraire2008compact} for the
viscous numerical flux. The DG discretization uses a mesh consisting of $3912$
cubic simplex elements. The second-order ODE in \eqnref{EQ: SIMP STRUCT} is
the governing equation for the mass-damper system with mass $m_s$, damping
constant $c_s = 1$, stiffness $k_s = 0$, and external force given from the
fluid as described in Section~\ref{SEC: APP FSI STRUCT}. The mesh motion is
determined from the position and velocity of the structure using the
blending maps introduced in \cite{persson2009discontinuous} and identical to
that used in Section 5.1 of \cite{zahr2016adjoint}. Snapshots of the vorticity
field and motion of the airfoil are shown in Figure~\ref{FIG: FOIL DAMPER VORT}
for a single configuration of the fluid-structure system.
\ifbool{fastcompile}{}{
\begin{figure}
 \includegraphics[width=0.3\textwidth]{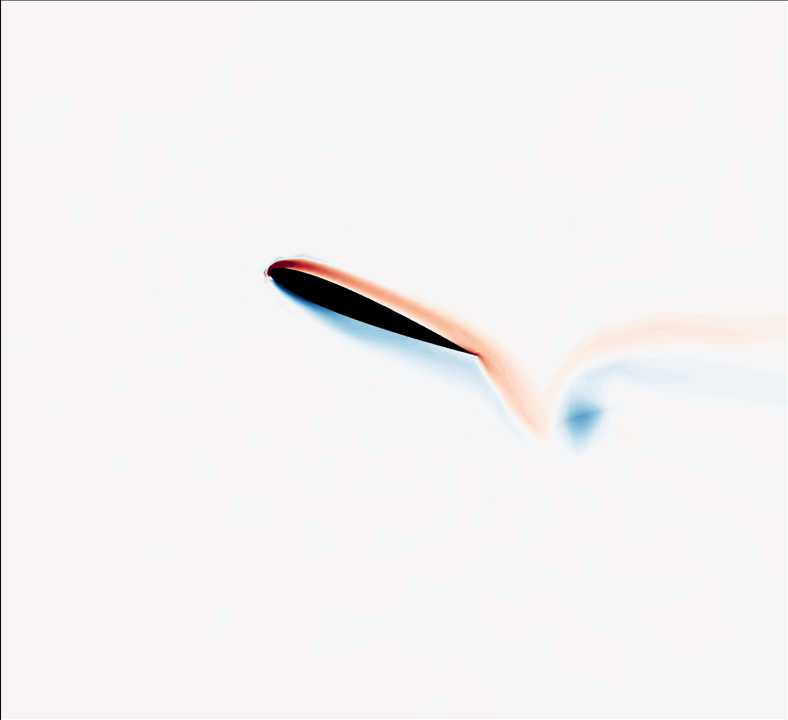} \quad
 \includegraphics[width=0.3\textwidth]{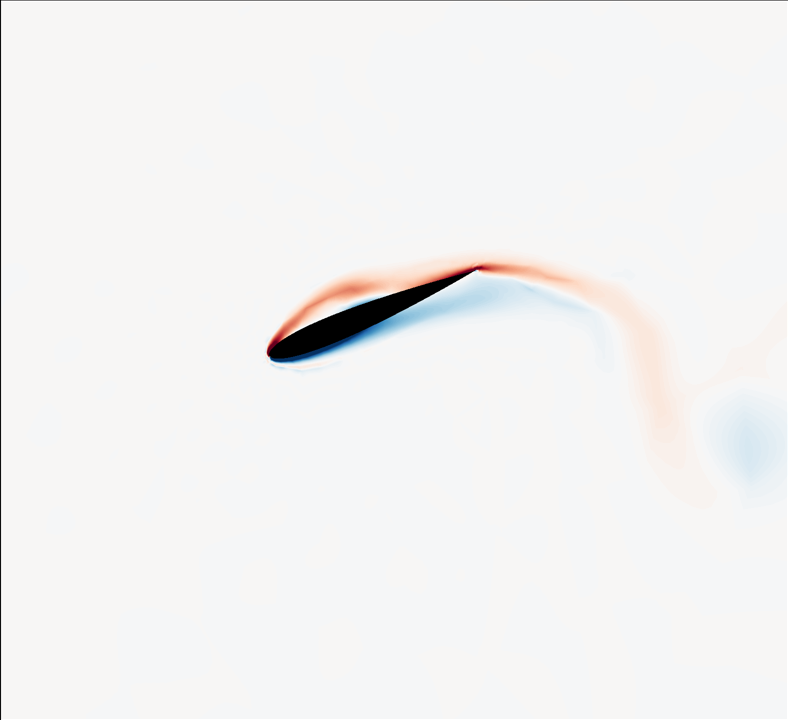} \quad
 \includegraphics[width=0.3\textwidth]{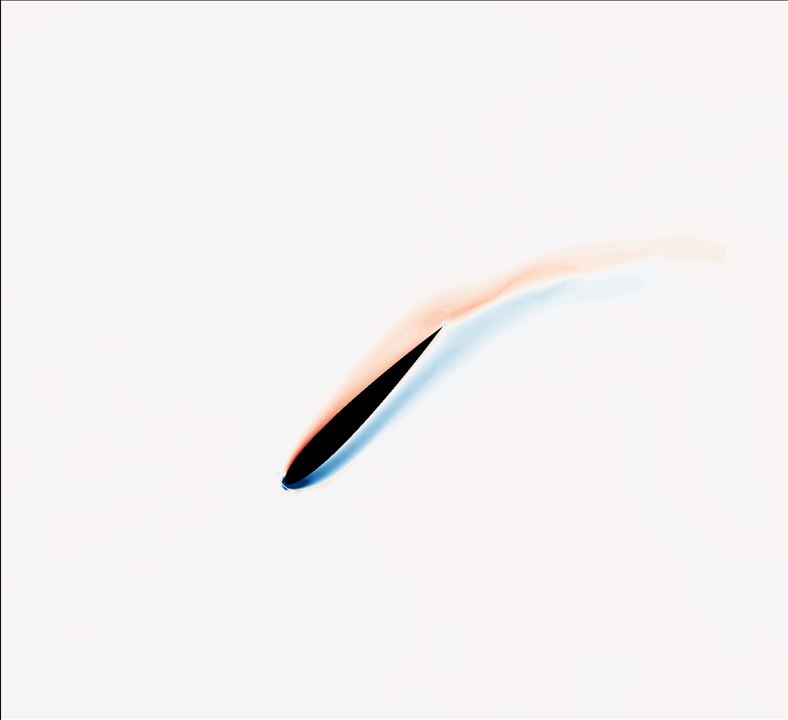} \\
 \includegraphics[width=0.3\textwidth]{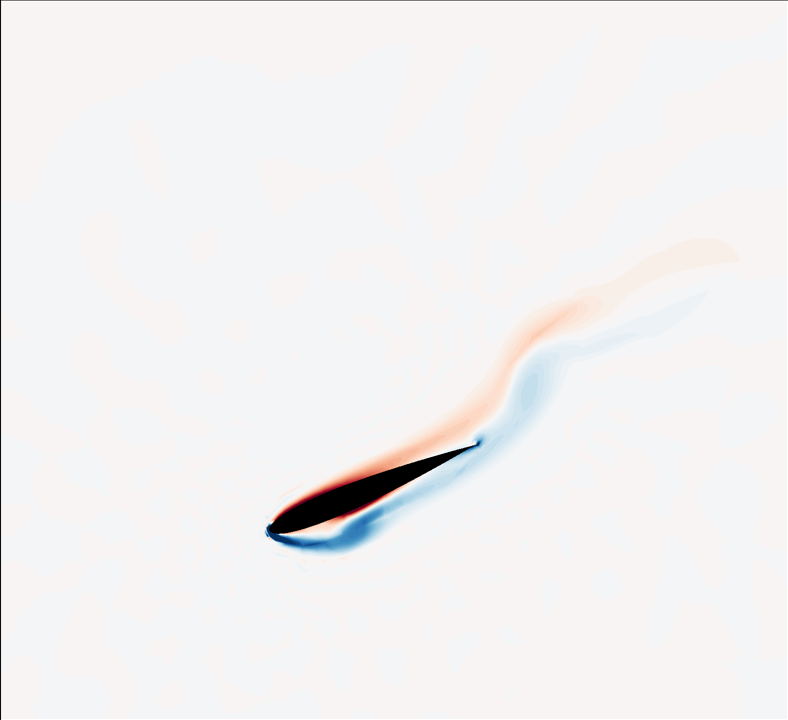} \quad
 \includegraphics[width=0.3\textwidth]{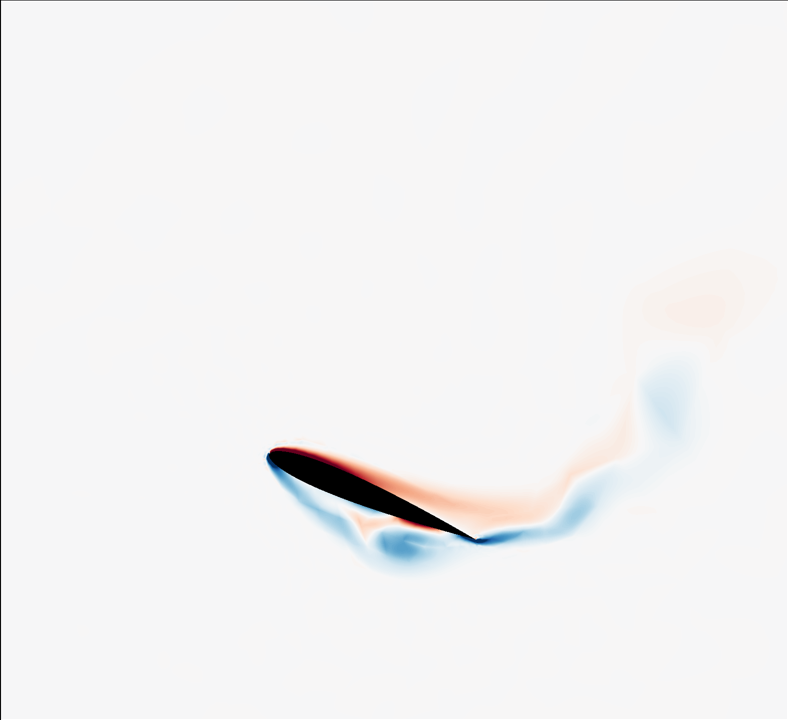} \quad
 \includegraphics[width=0.3\textwidth]{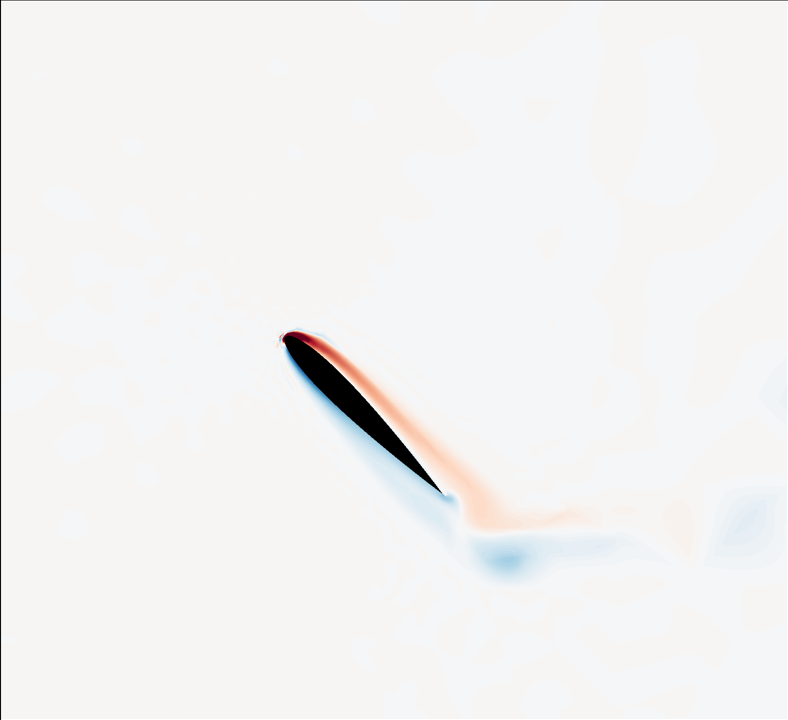}
 \caption{Airfoil motion and flow vorticity corresponding to foil-damper
          system under prescribed rotational motion
          $\theta(t) = \frac{\pi}{4}\cos(2\pi f t)$ with frequency $f = 0.2$
          at various snapshots in time:
          $t = 0.83,\,1.67,\,2.5,\,3.33,\,4.17,\,5.0$
          (\emph{left}-to-\emph{right}, \emph{top}-to-\emph{bottom}).}
 \label{FIG: FOIL DAMPER VORT}
\end{figure}
}

Our first numerical experiment studies the stability of the four proposed
predictors as a function of the mass ratio between the structure and fluid,
an important parameter that can impact the stability of partitioned solvers
as identified in \cite{causin2005added, van2007higher}, and time step size for
IMEX schemes up to fourth order. The mass ratio, $\bar{m}$, considered is the
ratio of the mass of the structure, $m$, to the mass of fluid displaced by the
structure, $\rho A$, where $\rho$ is the density of the fluid and
$A = 0.08221$ is the area of the airfoil. Since the isentropic Navier-Stokes
equations can be seen as an artificial compressibility formulation for
the incompressible Navier-Stokes equations
\cite{lin1995incompressible, desjardins1999incompressible}, we consider the
density to be constant and equal to the freestream $\rho = 1$. Variations in
the mass ratio are achieved by varying the mass of the structure with all
other parameters fixed. The stability results are summarized in
Figure~\ref{FIG: FOIL DAMPER STAB} where ~\ref{line:foil_damper_stable}~
indicates a $(\Delta t,\,\bar{m})$-pair that leads to a stable simulation and
~\ref{line:foil_damper_unstable}~ leads to an unstable one.
\ifbool{fastcompile}{}{
\begin{figure}
 \centering
 \begin{tikzpicture}
\begin{groupplot}[
    group style={
        group name=conv plots,
        group size=4 by 4,
        horizontal sep=1.0cm
    },
    xmode=log, ymode=log,
    xmin=0.003, xmax=0.1,
    ymin=1, ymax=1000.0/0.08221,
    width=4.5cm,
    height=4.5cm,
]

\nextgroupplot[ylabel={Mass ratio $(\bar{m})$}, title={weak Jacobi}]
 \addplot [blue, solid, thick, mark=triangle*, mark size=1.5, only marks, mark options={solid}]  table[x index=1, y expr=\thisrowno{0}/0.08221] {data/foil_damper/nacamsh1ref0p3_Re1000_Cpltyp0_Rk1.stable.dat};
 \addplot [red, solid, thick, mark=square*, mark size=1.5, only marks, mark options={solid}]  table[x index=1, y expr=\thisrowno{0}/0.08221] {data/foil_damper/nacamsh1ref0p3_Re1000_Cpltyp0_Rk1.unstable.dat};

\nextgroupplot[ytick=\empty, title={strong Jacobi}]
 \addplot [blue, solid, thick, mark=triangle*, mark size=1.5, only marks, mark options={solid}]  table[x index=1, y expr=\thisrowno{0}/0.08221] {data/foil_damper/nacamsh1ref0p3_Re1000_Cpltyp1_Rk1.stable.dat};
 \addplot [red, solid, thick, mark=square*, mark size=1.5, only marks, mark options={solid}]  table[x index=1, y expr=\thisrowno{0}/0.08221] {data/foil_damper/nacamsh1ref0p3_Re1000_Cpltyp1_Rk1.unstable.dat};

\nextgroupplot[ytick=\empty, title={weak Gauss-Seidel}]
 \addplot [blue, solid, thick, mark=triangle*, mark size=1.5, only marks, mark options={solid}]  table[x index=1, y expr=\thisrowno{0}/0.08221] {data/foil_damper/nacamsh1ref0p3_Re1000_Cpltyp2_Rk1.stable.dat};
 \addplot [red, solid, thick, mark=square*, mark size=1.5, only marks, mark options={solid}]  table[x index=1, y expr=\thisrowno{0}/0.08221] {data/foil_damper/nacamsh1ref0p3_Re1000_Cpltyp2_Rk1.unstable.dat};

\nextgroupplot[ytick=\empty, title={strong Gauss-Seidel}]
 \addplot [blue, solid, thick, mark=triangle*, mark size=1.5, only marks, mark options={solid}]  table[x index=1, y expr=\thisrowno{0}/0.08221] {data/foil_damper/nacamsh1ref0p3_Re1000_Cpltyp3_Rk1.stable.dat};
 \addplot [red, solid, thick, mark=square*, mark size=1.5, only marks, mark options={solid}]  table[x index=1, y expr=\thisrowno{0}/0.08221] {data/foil_damper/nacamsh1ref0p3_Re1000_Cpltyp3_Rk1.unstable.dat};

\nextgroupplot[ylabel={Mass ratio $(\bar{m})$}]
 \addplot [blue, solid, thick, mark=triangle*, mark size=1.5, only marks, mark options={solid}]  table[x index=1, y expr=\thisrowno{0}/0.08221] {data/foil_damper/nacamsh1ref0p3_Re1000_Cpltyp0_Rk2.stable.dat};
\label{line:foil_damper_stable}
 \addplot [red, solid, thick, mark=square*, mark size=1.5, only marks, mark options={solid}]  table[x index=1, y expr=\thisrowno{0}/0.08221] {data/foil_damper/nacamsh1ref0p3_Re1000_Cpltyp0_Rk2.unstable.dat};
\label{line:foil_damper_unstable}

\nextgroupplot[ytick=\empty]
 \addplot [blue, solid, thick, mark=triangle*, mark size=1.5, only marks, mark options={solid}]  table[x index=1, y expr=\thisrowno{0}/0.08221] {data/foil_damper/nacamsh1ref0p3_Re1000_Cpltyp1_Rk2.stable.dat};
 \addplot [red, solid, thick, mark=square*, mark size=1.5, only marks, mark options={solid}]  table[x index=1, y expr=\thisrowno{0}/0.08221] {data/foil_damper/nacamsh1ref0p3_Re1000_Cpltyp1_Rk2.unstable.dat};

\nextgroupplot[ytick=\empty]
 \addplot [blue, solid, thick, mark=triangle*, mark size=1.5, only marks, mark options={solid}]  table[x index=1, y expr=\thisrowno{0}/0.08221] {data/foil_damper/nacamsh1ref0p3_Re1000_Cpltyp2_Rk2.stable.dat};
 \addplot [red, solid, thick, mark=square*, mark size=1.5, only marks, mark options={solid}]  table[x index=1, y expr=\thisrowno{0}/0.08221] {data/foil_damper/nacamsh1ref0p3_Re1000_Cpltyp2_Rk2.unstable.dat};

\nextgroupplot[ytick=\empty]
 \addplot [blue, solid, thick, mark=triangle*, mark size=1.5, only marks, mark options={solid}]  table[x index=1, y expr=\thisrowno{0}/0.08221] {data/foil_damper/nacamsh1ref0p3_Re1000_Cpltyp3_Rk2.stable.dat};
 \addplot [red, solid, thick, mark=square*, mark size=1.5, only marks, mark options={solid}]  table[x index=1, y expr=\thisrowno{0}/0.08221] {data/foil_damper/nacamsh1ref0p3_Re1000_Cpltyp3_Rk2.unstable.dat};

\nextgroupplot[ylabel={Mass ratio $(\bar{m})$}]
 \addplot [blue, solid, thick, mark=triangle*, mark size=1.5, only marks, mark options={solid}]  table[x index=1, y expr=\thisrowno{0}/0.08221] {data/foil_damper/nacamsh1ref0p3_Re1000_Cpltyp0_Rk3.stable.dat}; 
 \addplot [red, solid, thick, mark=square*, mark size=1.5, only marks, mark options={solid}]  table[x index=1, y expr=\thisrowno{0}/0.08221] {data/foil_damper/nacamsh1ref0p3_Re1000_Cpltyp0_Rk3.unstable.dat}; 

\nextgroupplot[ytick=\empty]
 \addplot [blue, solid, thick, mark=triangle*, mark size=1.5, only marks, mark options={solid}]  table[x index=1, y expr=\thisrowno{0}/0.08221] {data/foil_damper/nacamsh1ref0p3_Re1000_Cpltyp1_Rk3.stable.dat};
 \addplot [red, solid, thick, mark=square*, mark size=1.5, only marks, mark options={solid}]  table[x index=1, y expr=\thisrowno{0}/0.08221] {data/foil_damper/nacamsh1ref0p3_Re1000_Cpltyp1_Rk3.unstable.dat};

\nextgroupplot[ytick=\empty]
 \addplot [blue, solid, thick, mark=triangle*, mark size=1.5, only marks, mark options={solid}]  table[x index=1, y expr=\thisrowno{0}/0.08221] {data/foil_damper/nacamsh1ref0p3_Re1000_Cpltyp2_Rk3.stable.dat};
 \addplot [red, solid, thick, mark=square*, mark size=1.5, only marks, mark options={solid}]  table[x index=1, y expr=\thisrowno{0}/0.08221] {data/foil_damper/nacamsh1ref0p3_Re1000_Cpltyp2_Rk3.unstable.dat};

\nextgroupplot[ytick=\empty]
 \addplot [blue, solid, thick, mark=triangle*, mark size=1.5, only marks, mark options={solid}]  table[x index=1, y expr=\thisrowno{0}/0.08221] {data/foil_damper/nacamsh1ref0p3_Re1000_Cpltyp3_Rk3.stable.dat};
 \addplot [red, solid, thick, mark=square*, mark size=1.5, only marks, mark options={solid}]  table[x index=1, y expr=\thisrowno{0}/0.08221] {data/foil_damper/nacamsh1ref0p3_Re1000_Cpltyp3_Rk3.unstable.dat};

\nextgroupplot[ylabel={Mass ratio $(\bar{m})$}, xlabel={Time step ($\Delta t$)}]
 \addplot [blue, solid, thick, mark=triangle*, mark size=1.5, only marks, mark options={solid}]  table[x index=1, y expr=\thisrowno{0}/0.08221] {data/foil_damper/nacamsh1ref0p3_Re1000_Cpltyp0_Rk4.stable.dat}; 
 \addplot [red, solid, thick, mark=square*, mark size=1.5, only marks, mark options={solid}]  table[x index=1, y expr=\thisrowno{0}/0.08221] {data/foil_damper/nacamsh1ref0p3_Re1000_Cpltyp0_Rk4.unstable.dat}; 

\nextgroupplot[xlabel={Time step ($\Delta t$)}, ytick=\empty]
 \addplot [blue, solid, thick, mark=triangle*, mark size=1.5, only marks, mark options={solid}]  table[x index=1, y expr=\thisrowno{0}/0.08221] {data/foil_damper/nacamsh1ref0p3_Re1000_Cpltyp1_Rk4.stable.dat};
 \addplot [red, solid, thick, mark=square*, mark size=1.5, only marks, mark options={solid}]  table[x index=1, y expr=\thisrowno{0}/0.08221] {data/foil_damper/nacamsh1ref0p3_Re1000_Cpltyp1_Rk4.unstable.dat};

\nextgroupplot[xlabel={Time step ($\Delta t$)}, ytick=\empty]
 \addplot [blue, solid, thick, mark=triangle*, mark size=1.5, only marks, mark options={solid}]  table[x index=1, y expr=\thisrowno{0}/0.08221] {data/foil_damper/nacamsh1ref0p3_Re1000_Cpltyp2_Rk4.stable.dat};
 \addplot [red, solid, thick, mark=square*, mark size=1.5, only marks, mark options={solid}]  table[x index=1, y expr=\thisrowno{0}/0.08221] {data/foil_damper/nacamsh1ref0p3_Re1000_Cpltyp2_Rk4.unstable.dat};

\nextgroupplot[xlabel={Time step ($\Delta t$)}, ytick=\empty]
 \addplot [blue, solid, thick, mark=triangle*, mark size=1.5, only marks, mark options={solid}]  table[x index=1, y expr=\thisrowno{0}/0.08221] {data/foil_damper/nacamsh1ref0p3_Re1000_Cpltyp3_Rk4.stable.dat};
 \addplot [red, solid, thick, mark=square*, mark size=1.5, only marks, mark options={solid}]  table[x index=1, y expr=\thisrowno{0}/0.08221] {data/foil_damper/nacamsh1ref0p3_Re1000_Cpltyp3_Rk4.unstable.dat};

\end{groupplot}
\end{tikzpicture}
 \caption{Behavior of the predictor-based partitioned schemes for a range
          of mass ratios and time steps for IMEX1-IMEX4 (\emph{top} to
          \emph{bottom}) schemes with the weak Jacobi
          predictor (\emph{left}), strong Jacobi predictor
          (\emph{center left}), weak Gauss-Seidel predictor
          (\emph{center right}), and strong Gauss-Seidel predictor
          (\emph{right}). \emph{Legend}: \ref{line:foil_damper_stable}
          indicates a stable simulation and \ref{line:foil_damper_unstable}
          indicates an unstable simulation.}
 \label{FIG: FOIL DAMPER STAB}
\end{figure}
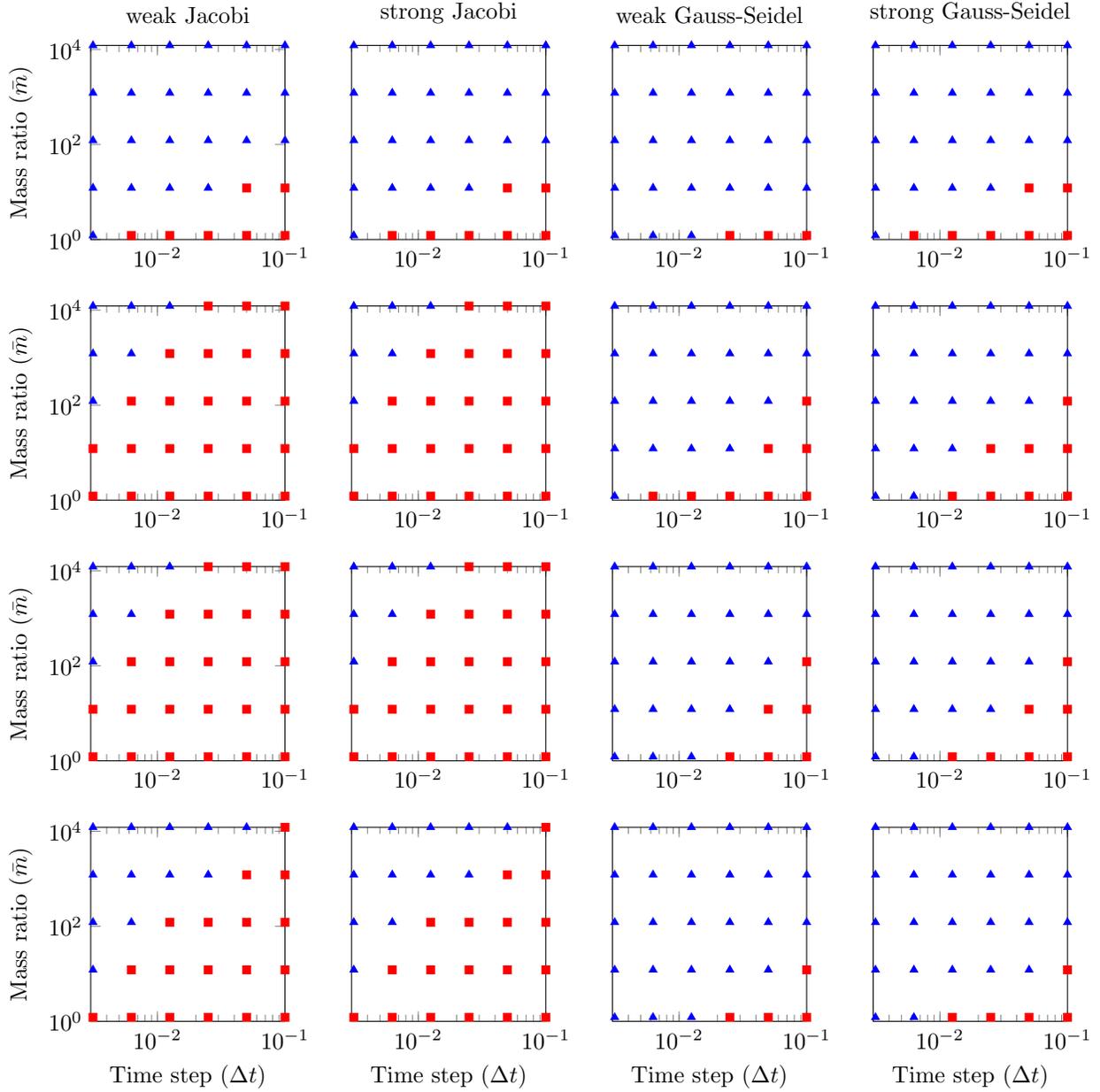
}
This figure shows the weak and strong Jacobi predictors lack robustness beyond
the first-order scheme since they only lead to stable simulations for
small step sizes or large mass ratios, while the Gauss-Seidel predictors are
stable across a larger set of $(\dt{},\,\bar{m})$ pairs. This does not
contradict the stability theory in Section~\ref{SEC: ACCURATE-STABILITY} since
the robustness issues manifest as a nonlinear instabilities that come from
lagging the mesh motion to the previous time step during the fluid solve.
Figure~\ref{FIG: FOIL DAMPER STAB} also shows that all schemes are stable once
the time step is sufficiently small, at least for this range of mass ratios
considered. The first-order IMEX scheme is the most robust, which is expected
given the large amount of numerical dissipation associated with first-order
solvers. This figure also highlights the robustness of the
proposed solver, particularly with the Gauss-Seidel predictor, since the
maximum stable time step is three orders of magnitude larger than the maximum
stable time step of a fluid-only simulation with RK4, indicating the scheme
benefits from treating both subsystems implicitly and the coupling correction
explicitly.

With the stability of the predictors established, we confirm the order of
accuracy for the Gauss-Seidel predictors in \figref{FIG: FOIL_DAMP_ERR}
up to fourth order. The error metric used is the error in the time-integrated
vertical force the fluid exerts on the structure, i.e., the integral of the
fluid stress tensor over the airfoil over time, denoted $e_\text{FSI2}$.
The temporal integral is computed to exactly the same order as the
semi-discrete system by recasting the time integral to an ODE and applying
the same IMEX scheme, i.e., solver-consistent integration of quantities of
interest \cite{zahr2016adjoint}. A reference solution is computed using the
IMEX5 scheme with $\Delta t = 3.125 \times 10^{-3}$ and strong Gauss-Seidel
predictor. 
\ifbool{fastcompile}{}{
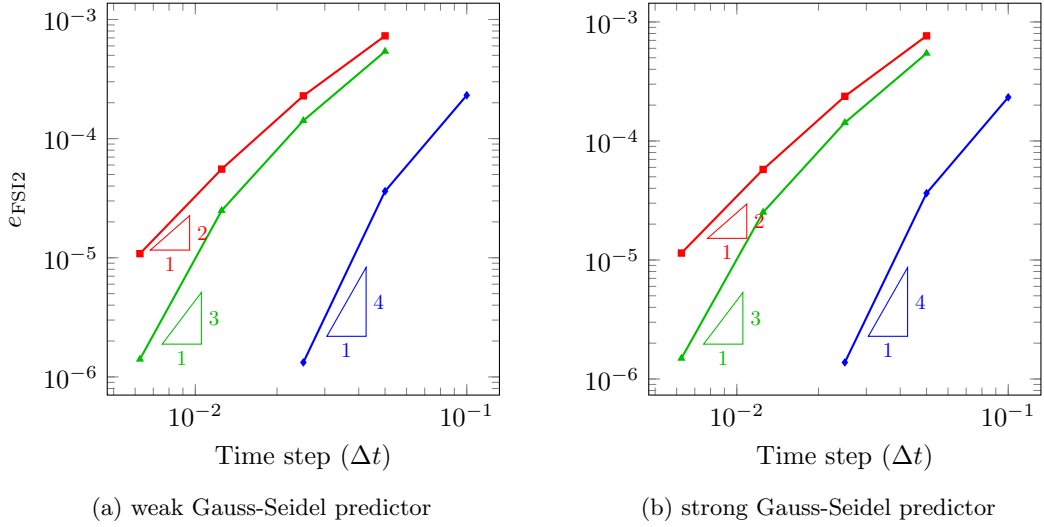
\begin{figure}[ht]
 \centering
 \begin{subfigure}[b]{0.42\textwidth}
     \begin{tikzpicture}
\begin{loglogaxis}[
    width=0.98\textwidth,
    height=0.98\textwidth,
    xlabel={Time step ($\Delta t$)},
    ylabel={$e_\text{FSI2}$}]

\addplot [red, solid, thick, mark=square*, mark size=1, mark options={solid}]  table[x index=0, y index=7] {data/foil_damper/nacamsh1ref0p3_Re1000_Mridx1_Cpltyp2_Rk2.summ.dat}; 
\addplot [green!75!black, solid, thick, mark=triangle*, mark size=1, mark options={solid}]  table[x index=0, y index=7] {data/foil_damper/nacamsh1ref0p3_Re1000_Mridx1_Cpltyp2_Rk3.summ.dat}; 
\addplot [blue, solid, thick, mark=diamond*, mark size=1, mark options={solid}]  table[x index=0, y index=7, select coords between index={0}{2}] {data/foil_damper/nacamsh1ref0p3_Re1000_Mridx1_Cpltyp2_Rk4.summ.dat}; 

\logLogSlopeTriangle{0.21}{0.1}{0.37}{2}{red};
\logLogSlopeTriangle{0.24}{0.1}{0.13}{3}{green!75!black};
\logLogSlopeTriangle{0.66}{0.1}{0.15}{4}{blue};

\end{loglogaxis}

\end{tikzpicture}
     \caption{weak Gauss-Seidel predictor}
 \end{subfigure} \quad
 \begin{subfigure}[b]{0.42\textwidth}
     \begin{tikzpicture}
\begin{loglogaxis}[
    width=0.98\textwidth,
    height=0.98\textwidth,
    xlabel={Time step ($\Delta t$)},
    ylabel={~~}]

\addplot [red, solid, thick, mark=square*, mark size=1, mark options={solid}]  table[x index=0, y index=7] {data/foil_damper/nacamsh1ref0p3_Re1000_Mridx1_Cpltyp3_Rk2.summ.dat}; \label{line:foil_damper_strong_GS:rk2}
\addplot [green!75!black, solid, thick, mark=triangle*, mark size=1, mark options={solid}]  table[x index=0, y index=7] {data/foil_damper/nacamsh1ref0p3_Re1000_Mridx1_Cpltyp3_Rk3.summ.dat}; \label{line:foil_damper_strong_GS:rk3}
\addplot [blue, solid, thick, mark=diamond*, mark size=1, mark options={solid}]  table[x index=0, y index=7, select coords between index={0}{2}] {data/foil_damper/nacamsh1ref0p3_Re1000_Mridx1_Cpltyp3_Rk4.summ.dat}; \label{line:foil_damper_strong_GS:rk4}

\logLogSlopeTriangle{0.25}{0.1}{0.4}{2}{red};
\logLogSlopeTriangle{0.24}{0.1}{0.13}{3}{green!75!black};
\logLogSlopeTriangle{0.66}{0.1}{0.15}{4}{blue};

\end{loglogaxis}

\end{tikzpicture}
     \caption{strong Gauss-Seidel predictor}
 \end{subfigure}
 \caption{Convergence of the IMEX2 (\ref{line:foil_damper_strong_GS:rk2}), IMEX3 (\ref{line:foil_damper_strong_GS:rk3}), and IMEX4 (\ref{line:foil_damper_strong_GS:rk4}) with Gauss-Seidel type predictors as applied to the foil-damper system. Both predictors achieve the design orders and give very similar levels of accuracy.}
 \label{FIG: FOIL_DAMP_ERR}
 \end{figure}
}

\section{Application to particle-laden flows}\label{SEC: APP PARTICLE FLOWS}
Our final application is a two-phase particle-laden flow that is common in
biological flows \cite{kleinstreuer2003laminar}, plasma problems \cite{jacobs2009implicit}, and environmental
flows \cite{ferrante2003physical}, among others. In these flows, momentum
and energy are exchanged between the carrier flow and small, immiscible
particles. This interaction plays an important role in both phases of the
flow and results in complex behavior.

The governing equations for the carrier flow are the unsteady compressible
Navier-Stokes equations (\ref{EQ: FSI GOVERN}) with a source term that accounts for
the momentum and energy the particles contribute to the flow
\begin{equation} \label{EQ: CARRIER GOVERN}
 \begin{aligned}
  \frac{\partial \rho}{\partial t}  + \frac{\partial}{\partial x_i}
  (\rho u_i) &= 0, \\
  \frac{\partial}{\partial t} (\rho u_i) +
  \frac{\partial}{\partial x_i} (\rho u_i u_j+ p)  &=
  f_i + \frac{\partial \tau_{ij}}{\partial x_j}
  \quad\text{for }i=1,2,\dots,d, \\
  \frac{\partial}{\partial t} (\rho E) +
  \frac{\partial}{\partial x_i} \left(u_j(\rho E+p)\right) &=
  f_j u_j
  -\frac{\partial q_j}{\partial x_j}
  +\frac{\partial}{\partial x_j}(u_j\tau_{ij}), \\
 \end{aligned}
\end{equation}
in the spatial-temporal domain $(x,\,t) \in \Omega \times (0,\,T]$, where
$f = \begin{bmatrix} f_1 & \cdots & f_d \end{bmatrix}^T$ is the force the
particles exert on the flow and all quantities are defined in
(\ref{EQ: ns1})-(\ref{EQ: ns5}). The force a system of $M$ particles 
at positions $x_1,\,\dots,\,x_M$ with velocities $v_1,\,\dots,\,v_M$ exert
on the fluid at a position $x$ is approximated as
\begin{equation} \label{EQ: PL_Fluid_Force}
 \begin{aligned}
  f(x) = -\sum_{k=1}^M m_p \frac{u - v_k}{\tau_p}
                           D(\norm{x - x_k}_2),
 \end{aligned}
\end{equation}
where $m_p$ and $\tau_p$ are the mass and response time of the particle,
and $D(r)$ is an approximated Dirac delta function
\begin{equation}
D(r) = \frac{1}{(2\pi \sigma^2)^{d/2} }\exp\left(-\frac{r^2}{2\sigma^2}\right).
\end{equation}
For particles of diameter $d_p$ and density $\rho_p$, Stokes' drag law
gives the following relation for the particle response time
\begin{equation}
 \tau_p = \frac{d_p^2\rho_p}{18\mu},
\end{equation}
where $\mu$ is the dynamic viscosity of the fluid. Finally, the equations of
motion for the system of $M$ particles are derived from Newton's second law
as the following system $2M$ of ODEs
\begin{equation} \label{EQ: PL_ODE}
 \begin{aligned}
  \oder{x_k}{t} &= v_k \\
  m_p \oder{v_k}{t} &=
  m_p \frac{u(x_k,\,t) - v_k}{\tau_p}
 \end{aligned}
\end{equation}
for $k = 1,\,2,\,\dots,\,M$, where $u(x_k,\,t)$ is the flow velocity
at position $x_k$ and time $t$. The system of ODEs is expressed
compactly as
\begin{equation}
 \mass[q]\dot\qbm = \res[q](\qbm,\,\cpl[q]),
\end{equation}
where the mass matrix, $\mass[q]$, is
\begin{equation}
 \mass[q] = \begin{bmatrix} 1 &        &   &     &        & \\
                              & \ddots &   &     &        & \\
                              &        & 1 &     &        & \\
                              &        &   & m_p &        & \\
                              &        &   &     & \ddots & \\
                              &        &   &     &        & m_p \\
            \end{bmatrix}
\end{equation}
and the generalized coordinates, $\qbm$, coupling term, $\cpl[q]$,
and velocity term, $\res[q](\qbm,\,\cpl[q])$ are
\begin{equation}
 \qbm = \begin{bmatrix}
         x_1 \\ \vdots \\ x_M \\ v_1 \\ \vdots \\ v_M
        \end{bmatrix}, \quad
 \cpl[q] = \begin{bmatrix}
            u(x_1,\,t) \\ \vdots \\ u(x_M,\,t)
           \end{bmatrix}, \quad
 \res[q](\qbm,\,\cpl[q]) =
        \begin{bmatrix}
         v_1 \\ \vdots \\ v_M \\
         m_p \frac{u(x_1,\,t)-v_1}{\tau_p} \\ \vdots \\
         m_p \frac{u(x_M,\,t)-v_M}{\tau_p}
        \end{bmatrix}.
\end{equation}

We consider the model problem of $100$ particles in an ideal gas flow in a
rectangular domain $[0,\,20]\times[0,\,15]$. The initial fluid state
(see \figref{FIG: PL_SETUP}) is a superposition of a uniform flow of velocity
$u_{\infty} = 1.0$ and angle $\theta = \arctan(1/2)$ and a vortex centered at
$(x_0,\,y_0) = (5.0,\,5.0)$ with characteristic radius $r_c=1.5$ and strength
parameter $\epsilon = 15$ \cite{persson2009discontinuous}
\begin{equation}
\begin{aligned}
u_1 &= u_{\infty}\left( \cos{\theta} - \frac{\epsilon(y-y_0)}{2\pi r_c}\exp\left(\frac{\varphi(x,\,y)}{2}\right)\right)\\
u_2 &= u_{\infty}\left( \sin{\theta} - \frac{\epsilon(x-x_0)}{2\pi r_c}\exp\left(\frac{\varphi(x,\,y)}{2}\right)\right)\\
\rho &= \rho_{\infty}\left( 1 - \frac{\epsilon^2(\gamma - 1)M_{\infty}^2}{8\pi^2}\exp\left(\frac{\varphi(x,\,y)}{2}\right)\right)^{\frac{1}{\gamma - 1}}\\
p &= p_{\infty}\left( 1 - \frac{\epsilon^2(\gamma - 1)M_{\infty}^2}{8\pi^2}\exp\left(\frac{\varphi(x,\,y)}{2}\right)\right)^{\frac{1}{\gamma - 1}}\\
\end{aligned}
\end{equation}
where $\varphi(x, y) = (1 - (x  - x_0)^2  -   (y - y_0)^2)/r_c^2$,
$M_{\infty} = 0.5$ is the Mach number, $\rho_{\infty} = 1.0$ is the density,
and $p_{\infty} = 1/\gamma M_{\infty}^2$ is the pressure. The specific heat
ratio and Reynolds number are $\gamma = 1.4$ and $Re = 200$, respectively.
The particles are initially at rest and randomly positioned near the
vortex center, i.e., the positions are drawn from the uniform distribution
over the interval $[x_0-3.0,\,x_0+3.0]\times[y_0-3.0,\,y_0+3.0]$.
\ifbool{fastcompile}{}{
\begin{figure}
 \centering
 \includegraphics[width=0.5\textwidth]{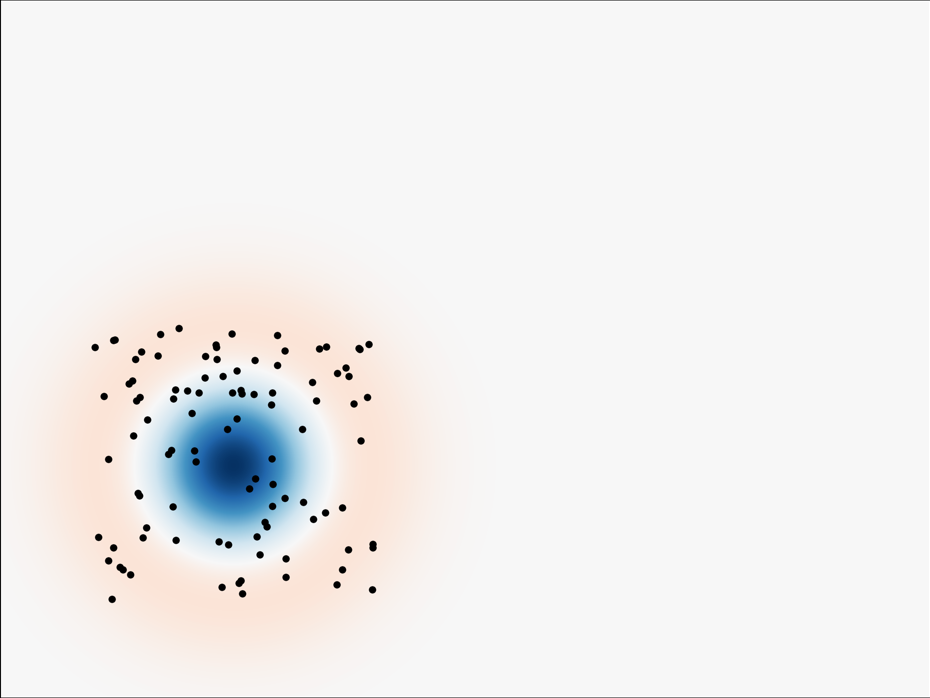}
 \caption{Vorticity profile and particles at $t = 0$}
 \label{FIG: PL_SETUP}
\end{figure}
}

The governing fluid equations (\ref{EQ: CARRIER GOVERN}) are discretized with
a standard high-order discontinuous Galerkin method using Roe's flux
\cite{roe1981approximate} for the inviscid numerical flux and the Compact DG
flux \cite{peraire2008compact} for the viscous numerical flux on a
structured mesh of $2400$ simplex elements. After the DG
spatial discretization is applied, the governing equations reduce to the
following system of ODEs
\begin{equation}
\label{EQ: CARRIER0}
 \mass[f]\stvcdot[f] = \res[f](\stvc[f],\,\cpl[f](\qbm))
\end{equation}
where $\mass[f]$ is the fixed mass matrix, $\stvc[f](t)$ is the semi-discrete
state vector, i.e., the discretization of $U$ on $\Omega$,
$\res[f](\stvc[f])$ is the spatial discretization of the inviscid and viscous
flux terms on $\Omega$, and $\cpl[f] = \qbm$ is the coupling term.

Given the extensive studies of the order of accuracy of the four predictors
in the previous section, we focus this section on stability of the weak
and strong Gauss-Seidel predictors. In the following numerical experiments,
the solution is integrated to time $t = 10.0$ with time step $\dt{} = 0.1$,
regardless of the IMEX-RK scheme used. Since the particle-laden flow is a
two-system multiphysics problem, the coupling structure conforms to the
format in \eqnref{EQ: SPE STRU} and only one coupling predictor is required:
$\cplprd[q]$. The weak Gauss-Seidel predictor lags the fluid velocity and
particle state to the previous time step
\begin{equation}
 \cplprd[q](\qbm,\,\stvc[f],\,\bar\qbm,\,\stvcbar[f]) =
 \begin{bmatrix}
  \bar{u}(\bar{x}_1)^T & 
  \cdots &
  \bar{u}(\bar{x}_M)^T
 \end{bmatrix}^T,
\end{equation}
while the strong Gauss-Seidel lags the fluid velocity to the current time
step, but uses the current particle state
\begin{equation}
 \cplprd[q](\qbm,\,\stvc[f],\,\bar\qbm,\,\stvcbar[f]) =
 \begin{bmatrix}
  \bar{u}(x_1)^T &
  \cdots &
  \bar{u}(x_M)^T
 \end{bmatrix}^T.
\end{equation}

We consider two scenarios:
\begin{inparaenum}[(1)]
\item light particles: $\rho_p = 0.1$, $d_p = 0.01$ and
\item heavy particles: $\rho_p = 1000.0$, $d_p = 0.01$.
\end{inparaenum}
In the first case, the particle mass is about $5.2\times 10^{-8}$ and particle
response time is about $1.1\times10^{-4}$ so the coupled system is stiff
considering the large coefficients in \eqnref{EQ: PL_ODE}. 
 Several simulations with a third-order DG discretization
  (quadratic $p = 2$ elements) with different time steps are
  performed to demonstrate the stability of the proposed high-order
  partitioned solver and predictors. \figref{FIG: PL_light_particles}a
  shows the particle trajectories for both weakly and strongly coupled
  Gauss-Seidel predictors using the second-order temporal discretization
  (IMEX2) with $\dt{} = 0.1$. In this extreme case, the weak
  Gauss-Seidel predictor exhibits a form of instability, which can be
  seen from the oscillations that appear in some particle trajectories
  (Figure~\ref{FIG: PL_light_particles}a); however, the strong Gauss-Seidel
  predictor gives smooth results. For smaller time steps, i.e., $\dt{} = 0.05$,
  the IMEX2 scheme with the weakly coupled Gauss-Seidel predictor leads to
  stable results. Interestingly, the IMEX3 and IMEX4 schemes do not exhibit
  the aforementioned instabilities as both the weakly and strongly coupled
  Gauss-Seidel predictors are stable even with larger time steps. This case
  demonstrates that although the weakly coupled Gauss-Seidel predictor is
  inferior to strongly coupled Gauss-Seidel predictor in some cases, both
  possess good stability properties considering the time step is three orders
  of magnitude larger than the particle response time. The particle
  trajectories for all these cases are provided in
  \figref{FIG: PL_light_particles}b.

  To close this section, a formal convergence study is conducted. The accuracy
  is quantified via the $L_{\infty}$-norm of the error in the flow solution
  at time $t = 10.0$
\begin{equation}
 \begin{aligned}
  e_\text{PL} &= \norm{\pstp[f]{N} - \stvc[f](10.0)}_\infty,
 \end{aligned}
\end{equation}
  where $\stvc[f](10.0)$ is the reference solution computed by using the IMEX4
  scheme with $\Delta t = 1.25\times10^{-2}$ and strong Gauss-Seidel predictor,
  and $\pstp[f]{N}$ is the flow state from the numerical solution at the final
  time step. The convergence result is presented in \figref{FIG: PL_ERR}, which
  illustrates the partitioned solver with both Gauss-Seidel predictors achieves
  the design order of the underlying IMEX scheme.

\ifbool{fastcompile}{}{
\begin{figure}[h]
\centering
\begin{subfigure}[t]{0.49\textwidth}
\centering
    \includegraphics[height=0.8\textwidth]{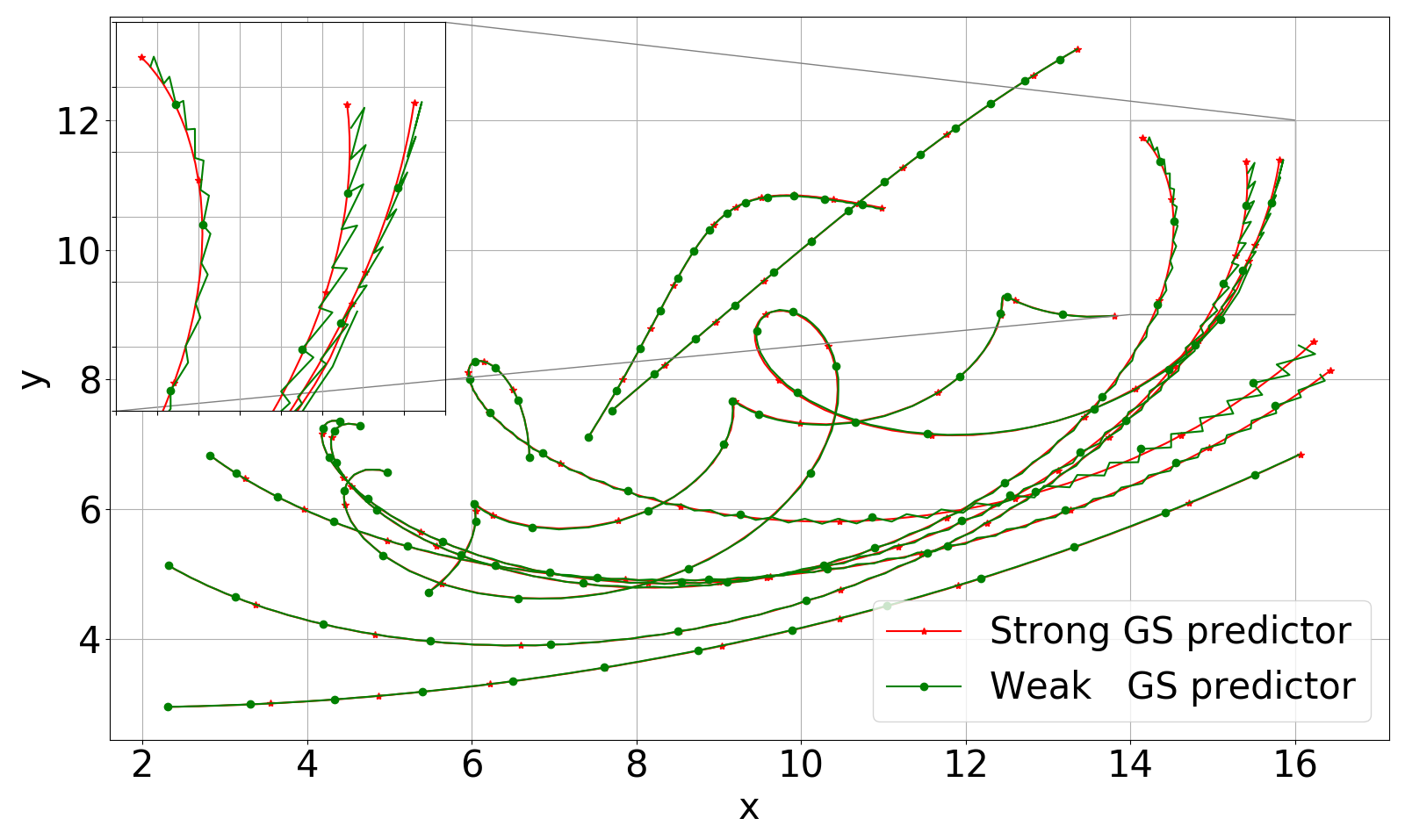}
    \caption{}
    \end{subfigure}
\begin{subfigure}[t]{0.49\textwidth}
\centering
    \includegraphics[height=0.8\textwidth]{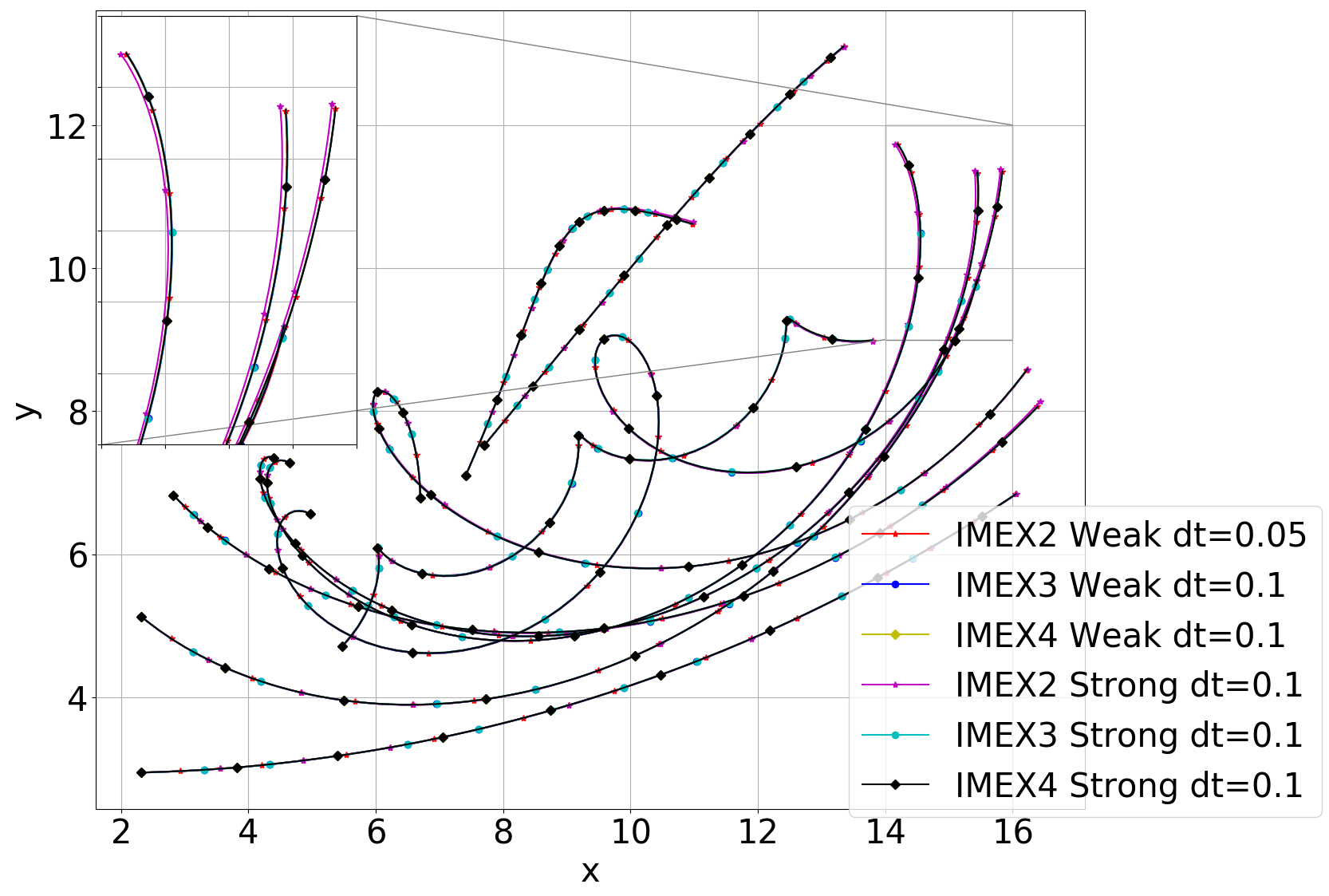}
    \caption{}
   \end{subfigure}
\caption{Sampling light particle trajectories in the unsteady compressible vortex: (a) Comparison of strong/weak GS coupling predictors for IMEX2 of $\Delta t = 0.1$ (b) Comparison of strong/weak GS coupling predictors for different order schemes.}
 \label{FIG: PL_light_particles}
\end{figure}
}

\ifbool{fastcompile}{}{
\begin{figure}[ht]
 \centering
 \begin{subfigure}[b]{0.42\textwidth}
     \begin{tikzpicture}
\begin{loglogaxis}[
    width=0.98\textwidth,
    height=0.98\textwidth,
    xlabel={Time step ($\Delta t$)},
    ylabel={$e_\text{PL}$}]

\addplot [red, solid, thick, mark=square*, mark size=1, mark options={solid}]  table[x index=0, y index=1] {data/particle_laden/Particle_laden_flow_fluid_RK2Weak_GS.dat}; 
\addplot [green!75!black, solid, thick, mark=triangle*, mark size=1, mark options={solid}]  table[x index=0, y index=1] {data/particle_laden/Particle_laden_flow_fluid_RK3Weak_GS.dat}; 
\addplot [blue, solid, thick, mark=diamond*, mark size=1, mark options={solid}]  table[x index=0, y index=1, select coords between index={0}{2}] {data/particle_laden/Particle_laden_flow_fluid_RK4Weak_GS.dat}; 

\logLogSlopeTriangle{0.28}{0.1}{0.37}{2}{red};
\logLogSlopeTriangle{0.38}{0.1}{0.13}{3}{green!75!black};
\logLogSlopeTriangle{0.82}{0.1}{0.25}{4}{blue};

\end{loglogaxis}

\end{tikzpicture}
     \caption{weak Gauss-Seidel predictor}
 \end{subfigure} \quad
 \begin{subfigure}[b]{0.42\textwidth}
     \begin{tikzpicture}
\begin{loglogaxis}[
    width=0.98\textwidth,
    height=0.98\textwidth,
    xlabel={Time step ($\Delta t$)},
    ylabel={~~}]

\addplot [red, solid, thick, mark=square*, mark size=1, mark options={solid}]  table[x index=0, y index=1] {data/particle_laden/Particle_laden_flow_fluid_RK2Strong_GS.dat}; \label{line:particle_laden_strong_GS:rk2}
\addplot [green!75!black, solid, thick, mark=triangle*, mark size=1, mark options={solid}]  table[x index=0, y index=1] {data/particle_laden/Particle_laden_flow_fluid_RK3Strong_GS.dat}; \label{line:particle_laden_strong_GS:rk3}
\addplot [blue, solid, thick, mark=diamond*, mark size=1, mark options={solid}]  table[x index=0, y index=1, select coords between index={0}{2}] {data/particle_laden/Particle_laden_flow_fluid_RK4Strong_GS.dat}; \label{line:particle_laden_strong_GS:rk4}

\logLogSlopeTriangle{0.28}{0.1}{0.48}{2}{red};
\logLogSlopeTriangle{0.52}{0.1}{0.33}{3}{green!75!black};
\logLogSlopeTriangle{0.76}{0.1}{0.15}{4}{blue};

\end{loglogaxis}

\end{tikzpicture}
     \caption{strong Gauss-Seidel predictor}
 \end{subfigure}
 \caption{Convergence of the IMEX2 (\ref{line:particle_laden_strong_GS:rk2}), IMEX3 (\ref{line:particle_laden_strong_GS:rk3}), and IMEX4 (\ref{line:particle_laden_strong_GS:rk4}) with Gauss-Seidel type predictors as applied to the particle-laden flow problem. Both predictors achieve the design orders and give very similar levels of accuracy.}
 \label{FIG: PL_ERR}
 \end{figure}
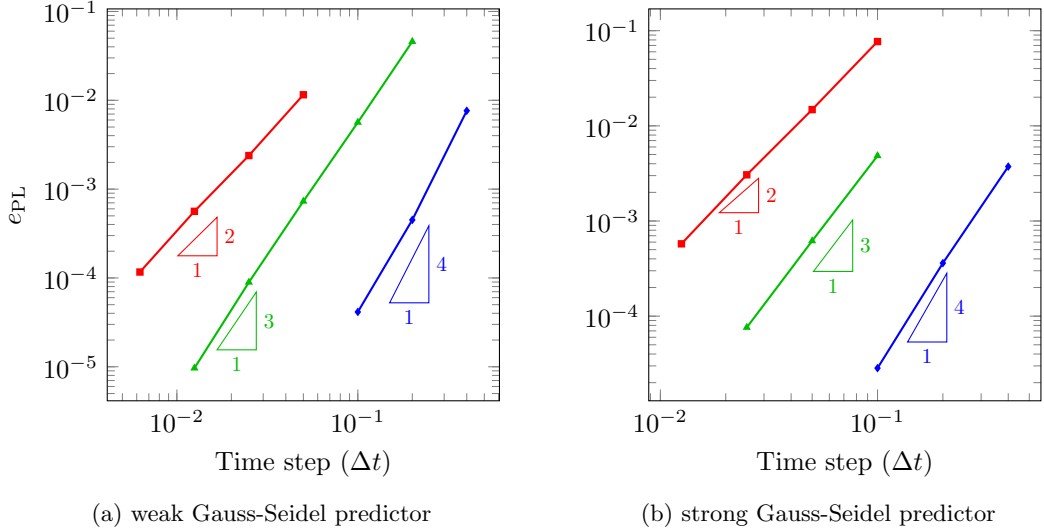
}

For the case with heavy particles, the particle mass is about
$5.2\times 10^{-4}$ and the particle response time is about $1.11$. In this
case, the coupling effect is stronger than for the light particles and
both predictors are stable for all discretization orders considered.
\figref{FIG: PL_particles} shows the vorticity profiles and
particle positions at several time instances for the simulations with
light and heavy particles. Light particles drift with the vortex, while
heavy particle advect with the flow since they are more affected by
inertial forces.
\ifbool{fastcompile}{}{
\begin{figure}
\centering
\begin{subfigure}[t]{0.32\textwidth}
    \includegraphics[width=1.0\textwidth]{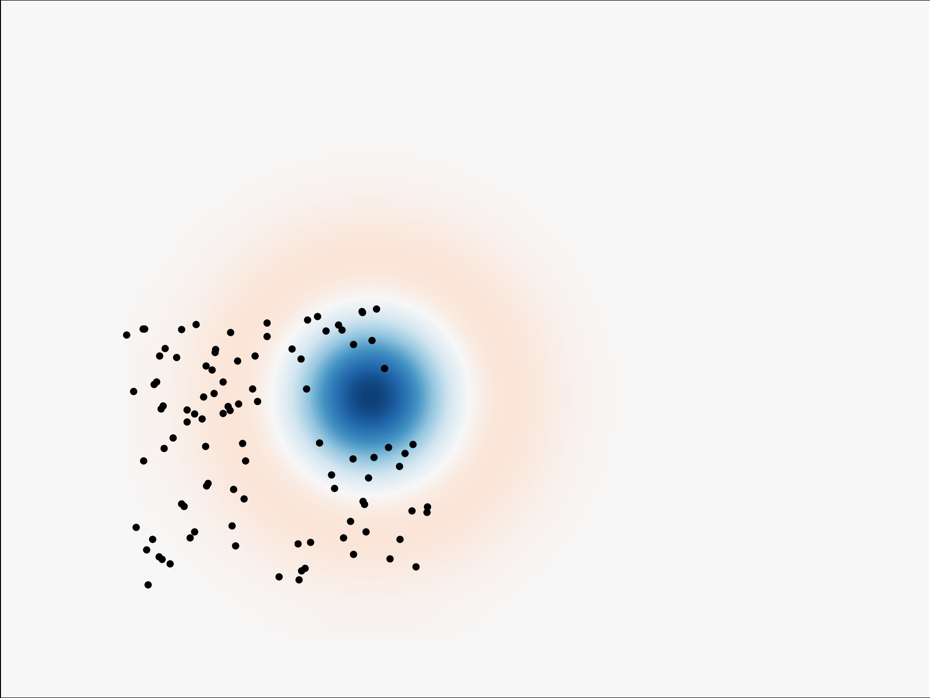}
    \caption{Heavy particle case}
\end{subfigure}
\begin{subfigure}[t]{0.32\textwidth}
    \includegraphics[width=1.0\textwidth]{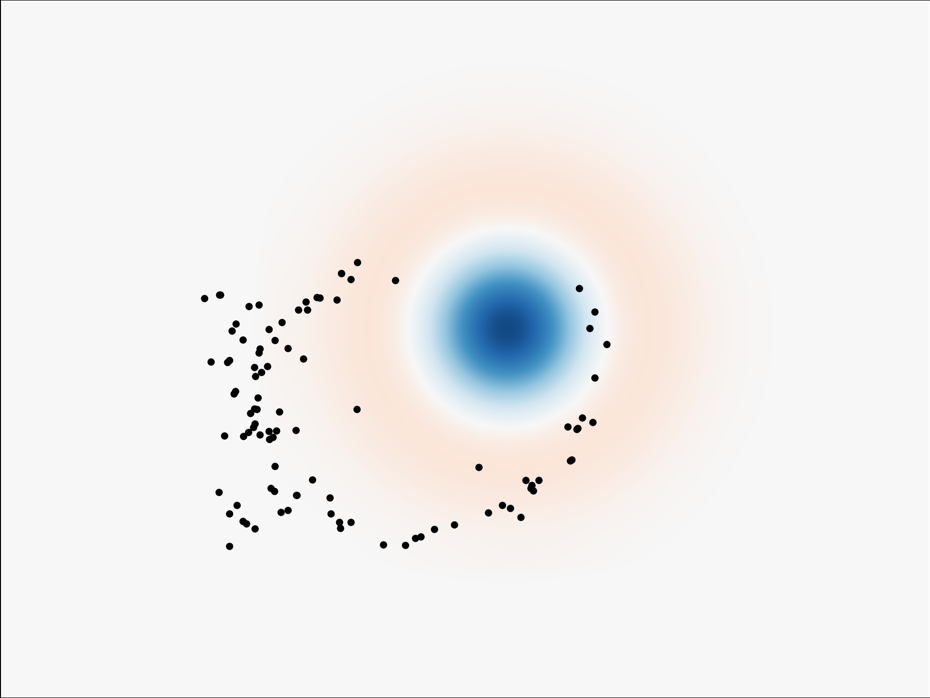}
    \caption{Heavy particle case}
\end{subfigure}
\begin{subfigure}[t]{0.32\textwidth}
    \includegraphics[width=1.0\textwidth]{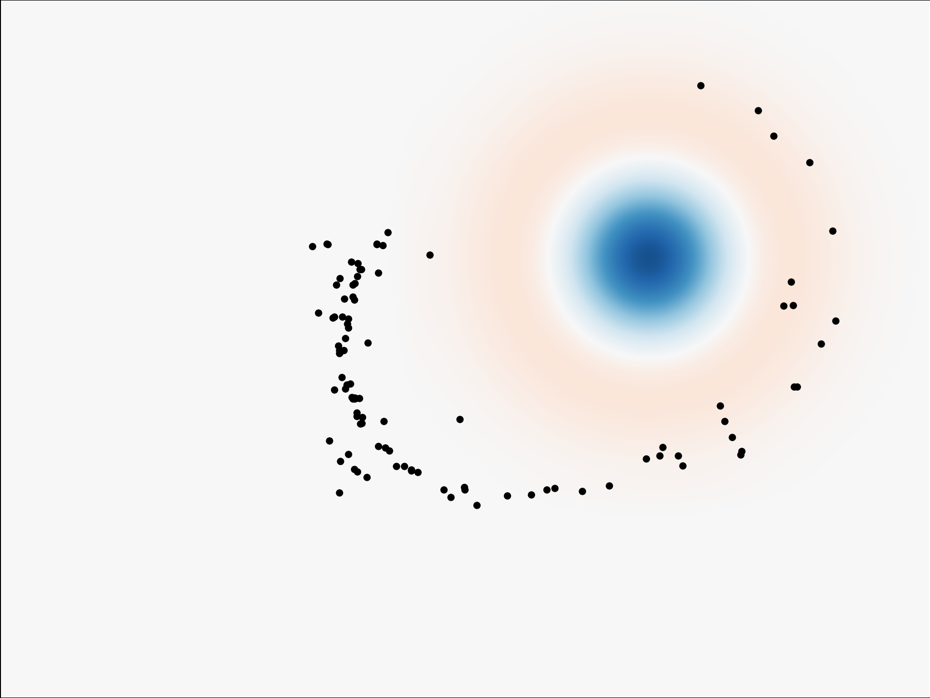}
    \caption{Heavy particle case}
\end{subfigure}

\begin{subfigure}[t]{0.32\textwidth}
    \includegraphics[width=1.0\textwidth]{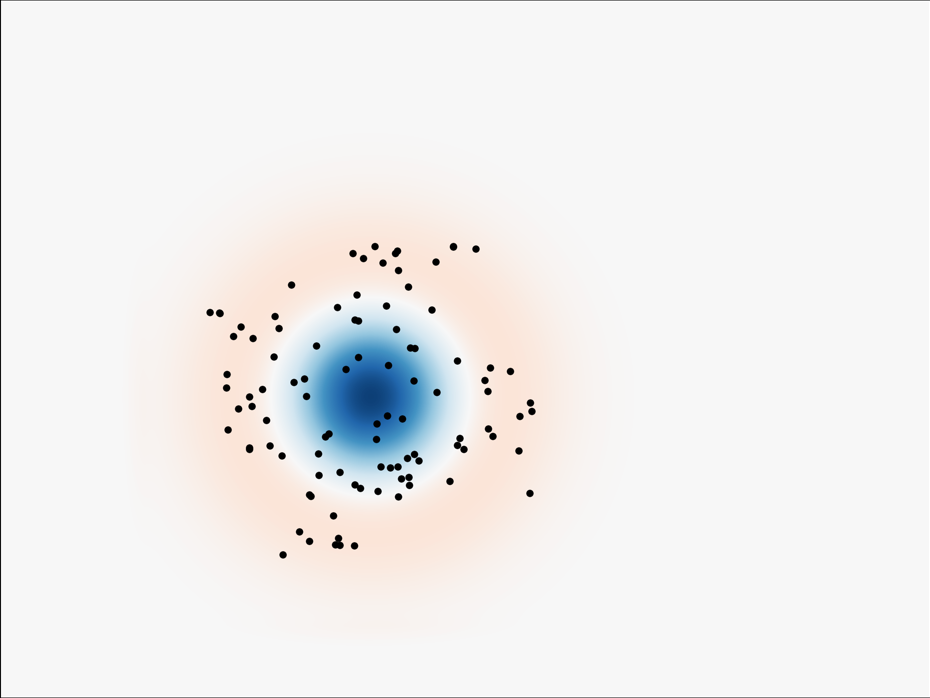}
    \caption{Light particle case}
\end{subfigure}
\begin{subfigure}[t]{0.32\textwidth}
    \includegraphics[width=1.0\textwidth]{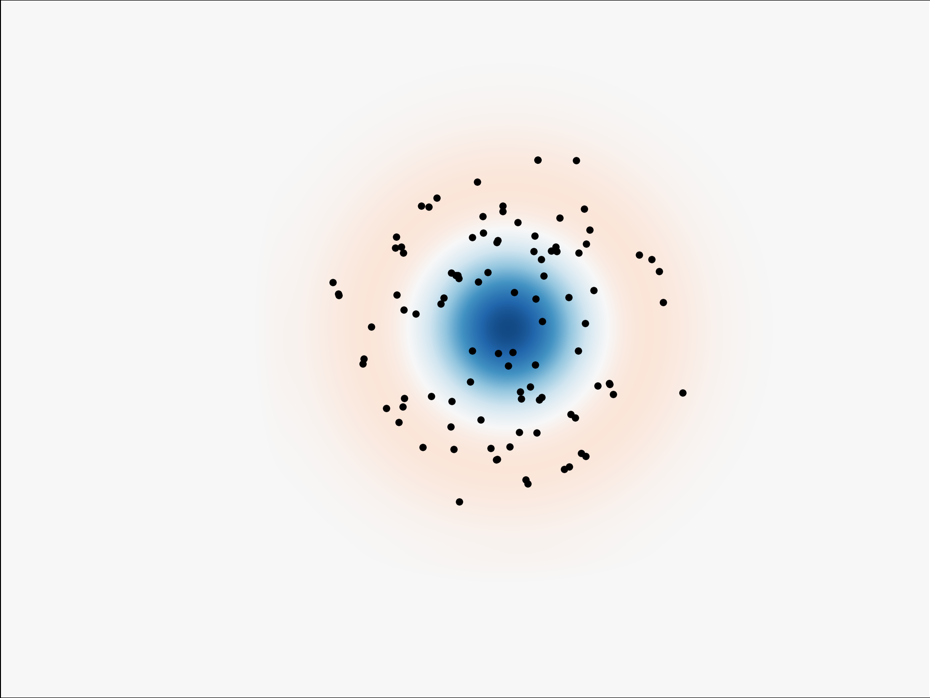}
    \caption{Light particle case}
\end{subfigure}
\begin{subfigure}[t]{0.32\textwidth}
    \includegraphics[width=1.0\textwidth]{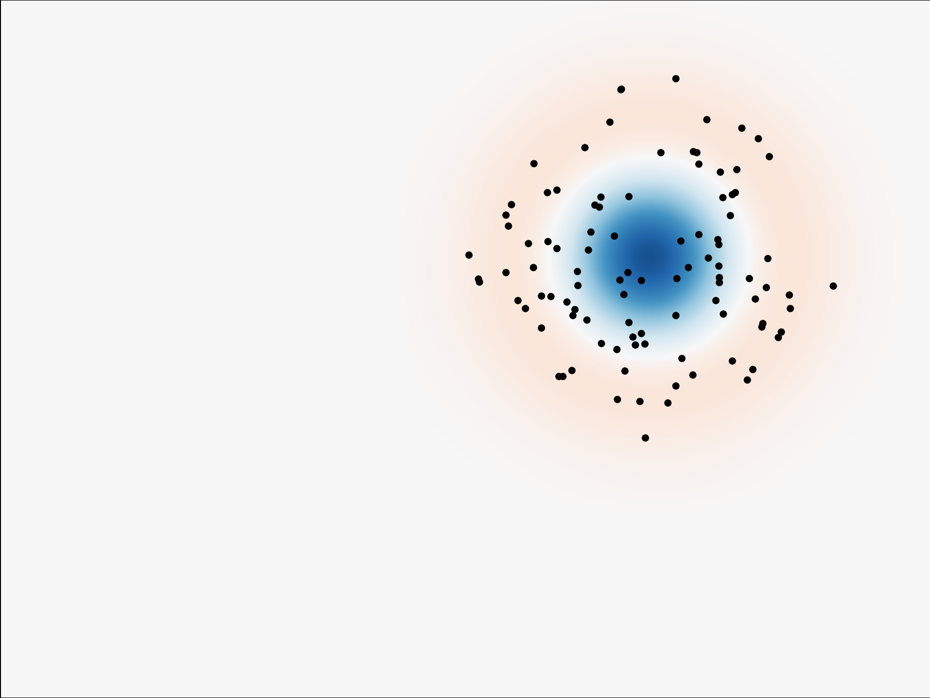}
    \caption{Light particle case}
\end{subfigure}
\caption{Particle laden flow --- vorticity profiles and particle positions at $t= 3.33,\,6.67,\,10$.}
\label{FIG: PL_particles}
\end{figure}
}


\section{Conclusions}
This paper introduces a framework for constructing high-order, linearly
stable, partitioned solvers for general multiphysics problems. The solvers
are constructed from an IMEX-RK discretization applied to the monolithic
system of $n$ systems of ODEs. A specific implicit-explicit decomposition
that introduces the concept of a predictor allows the monolithic systems
to be solved in a partitioned manner if the predictor meets basic
requirements. The four predictors, i.e., weak and strong
Jacobi and Gauss-Seidel predictors, introduced lead to different
IMEX-RK-based partitioned solvers, each with their own advantages and
disadvantages. The weak predictors require the least implementation
effort since they do not require any terms from the Jacobian of the
coupling term and therefore allows for maximal re-use of existing
software, while the strong predictors require the diagonal entries
from the Jacobian of the coupling term, which is unlikely to be available
in existing software. The Jacobi predictors allow for all subsystems to
be solved in parallel at a given stage of a given time step, while the
Gauss-Seidel predictors require the subsystems be solved sequentially.
Despite the simplicity and efficiency of the weak and Jacobi predictors
over the strong and Gauss-Seidel predictors, they have weaker linear
stability properties, which is shown theoretically and experimentally.
It is interesting to note that our linear stability analysis suggests
the strong Gauss-Seidel predictor is unconditionally linearly stable
in the context of the chosen model problem,
despite being a partitioned scheme.
The performance of the four partitioned solvers was demonstrated on a
slew of multiphysics problems, including an advection-diffusion-reaction
system, fluid-structure interaction problems, and particle-laden flow, where we
verified the design order of the IMEX scheme and studied various stability
properties. Future work will consider analysis of the nonlinear stability
of these schemes and derivation of the fully discrete sensitivity and
adjoint equations corresponding to these four solvers so they can be used
for gradient-based optimization of multiphysics systems.


\section*{Acknowledgments}
This work was supported in part by the Luis W. Alvarez Postdoctoral Fellowship
(MZ), by the Director, Office of Science, Office of Advanced Scientific
Computing Research, U.S. Department of Energy under Contract
No. DE-AC02-05CH11231 (MZ, PP), and by the NASA National Aeronautics and Space
Administration under grant number NNX16AP15A (MZ, PP). The content of this
publication does not necessarily reflect the position or policy of any of
these supporters, and no official endorsement should be inferred.


\appendix
\section{Stability analysis of high-order, partitioned IMEX-RK solvers}
\label{SEC: APPENDIX-A}
In this section, we analyze the linear stability of the high-order IMEX-RK
schemes: 2nd-order 2-stage trapezoidal rule,
3rd-order 4-stage ARK3(2)4L[2]SA, and 4th-order 6 stage ARK4(3)6L[2]SA
\cite{christopher2001additive} based on the model problem in
(\ref{EQ: MODEL-2WAY-SIMPLE}) and the predictor-based implicit-explicit
partitions in \tabref{TAB: MODEL-2WAY-SIMPLE PARTITION}.   

The linear stability analysis of the high-order IMEX-RK schemes proceeds
according to the procedure outlined in Section~\ref{SEC: ACCURATE-STABILITY},
namely, the IMEX-RK scheme is written as a one-step update scheme
(\ref{EQ: ITERATIVE IMEX}) and region where the spectral radius of the update
matrix, $\rho(\Cboldcal)$, is less than unity is identified. For brevity, we
directly write the spectral radius and subsequently identify stable regions.

It can be shown that the spectral radius of the one-step IMEX-RK update matrix
corresponding to the 2nd-order 2-stage trapezoidal rule and strong
Gauss-Seidel predictor is
\begin{equation}
 \rho(\Cboldcal) = \max\left\{1,\, \left| \frac{1 +\frac{\dt{} \lambda_1}{2}}{1 - \frac{ \dt{} \lambda_1}{2}}  \frac{1 +  \frac{\dt{} \lambda_2}{2}}{1 -  \frac{\dt{} \lambda_2}{2}} \right| \right\},
\end{equation}
which is independent of $\alpha$, less than unity for all $\dt{}$, and
therefore the scheme is unconditionally stable for all $\alpha$. The
spectral radius corresponding to the weak Gauss-Seidel scheme is
\begin{equation}
 \rho(\Cboldcal) = \max\left\{1,\,\left|\frac{1 +\frac{\dt{} (\lambda_1 + \lambda_2)}{2}(1 + \alpha)  +  \frac{\dt{}^2 \lambda_1\lambda_2}{4} (1 + \alpha)^2 + \frac{\dt{}^2 (\lambda_1^2 + \lambda_2^2)}{2}\alpha +  \frac{\dt{}^3 \lambda_1\lambda_2(\lambda_1 + \lambda_2)}{4}\alpha^2 }{(1 - (1 - \alpha)\frac{ \dt{} \lambda_1}{2})  (1 -  (1 - \alpha)\frac{\dt{} \lambda_2}{2})}  \right| \right\},
\end{equation}
which is unconditionally stable if and only if $\alpha = 0$. Finally, the
spectral radius corresponding to the strong Jacobi scheme is
\begin{equation}
 \rho(\Cboldcal) = \max\left\{1,\,\left| \frac{(1 +\frac{\dt{} \lambda_1}{2})(1 + \frac{ \dt{} \lambda_1}{2}) - \frac{\dt{}^3}{4}(\lambda_1^2\lambda_2 + \lambda_2^2\lambda_1)}{(1 -  \frac{\dt{} \lambda_2}{2})(1 -  \frac{\dt{} \lambda_2}{2})}  \right| \right\},
\end{equation}
which is not unconditionally stable.

For 3rd order 4-stage ARK3(2)4L[2]SA and 4th order 6-stage ARK4(3)6L[2]SA in
\cite{christopher2001additive}, we consider only the strong Gauss-Seidel predictor.
The spectral radius of the update matrices for the third and fourth order schemes
are
\begin{equation}
 \rho(\Cboldcal_3) = \max\left\{ 1,\, \left| \frac{p(\lambda_1\dt{},\,\lambda_2\dt{})}{q(\lambda_1\dt{},\,\lambda_2\dt{})}  \right| \right\}, \qquad
 \rho(\Cboldcal_4) = \max\left\{ 1,\, \left| \frac{\bar{p}(\lambda_1\dt{},\,\lambda_2\dt{})}{\bar{q}(\lambda_1\dt{},\,\lambda_2\dt{})}  \right| \right\},
\end{equation}
respectively, where $p$ and $q$ are 6th order polynomials and $\bar{p}$ and
$\bar{q}$ are 10th order polynomials
\begin{equation} 
 \begin{aligned}
  p(x_1, x_2) &= \sum_{i,j=0}^3 p_{ij}x_1^ix_2^j \\
  q(x_1, x_2) &= (1 - a x_1)^3(1 - a x_2)^3 \\
  \bar{p}(x_1, x_2) &= \sum_{i,j=0}^5 \bar{p}_{ij}x_1^ix_2^j \\
  \bar{q}(x_1, x_2) &= (1 - \bar{a} x_1)^5(1 - \bar{a} x_2)^5
 \end{aligned}
\end{equation}
and $a = 0.4358665216$ and $\bar{a} = 0.25$ are the coefficients of the second
entry on the diagonal of the implicit Runge-Kutta Butcher tableau for the
ARK3(2)4L[2]SA and ARK4(3)6L[2]SA schemes, respectively. From the coefficients of
$p$ and $\bar{p}$ in \tabref{TAB: f3} and \tabref{TAB: f4}, we observe that
\begin{equation}
|p_{ij}| \leq (-1)^{i+j} q_{ij},
\end{equation}
which implies the Gauss-Seidel predictors lead to unconditionally stable schemes
when $\lambda_1 \leq 0$ and $\lambda_2 \leq 0$. 
\begin{table}[!htb]
\begin{tabular}{c|cccc}
&   $0$  &  $1$ &   $2$  &  $3$\\
\hline
$0$&  $1$                   &  $- 0.307599564300000$      & $- 0.237660691030414$  & $0$\\
$1$&  $- 0.307599564300000$ &  $ 0.0946174918786356$       & $0.0731043252393467$   & $0$    \\
$2$&  $-0.237660691030414$  &  $ 0.0731043252393467$      & $0$                    & $ -0.0138993203184737$\\
$3$&  $0$                   &  $0$                        & $-0.0138993202982233$  & $ 0.00685679356380471$    \\
\end{tabular}
\caption{Coefficients of $p_{ij}$}
\label{TAB: f3}
\end{table}
\begin{table}[!htb]
\scalebox{0.65}{
\begin{tabular}{c|cccccc}
&   $0$          &  $1$ &   $2$  &  $3$ &  $4$   & $5$\\
\hline
$0$&  $1.0$                   &  $-0.25$                  & $-0.125$                & $  0.0104166666865151$    & $ 0.00911458332517619$  & $ 0.0 $  \\
$1$&  $-0.25$               &  $0.06245$                & $0.03125$               & $ -0.00260416668514596$     & $ -0.00407734171291718 $  & $ 0.0 $  \\
$2$&  $-0.125$              &  $0.03125$                & $ 0.015625$             & $  0.00606937261393480 $    & $ -0.00389797283406001$  & $1.71399137262393\textrm{e-}4  $  \\
$3$&  $0.0104166666865151$  &  $- 0.00260416668792851$  & $0.00606937262572686$   & $ -0.00535453941337523$     & $ 0.00164424787309041$  & $-8.30991742855789\textrm{e-}5  $  \\
$4$&  $0.00911458332517619$ &  $- 0.00407734171177262$  & $- 0.00389797283635686$ & $  0.00164424787343589$     & $ -8.99044034172063\textrm{e-}5$  & $ 7.55399650866135\textrm{e-}6 $  \\
$5$&  $0.0$                   &  $0.0$                      & $1.71399137131092\textrm{e-}4$  & $ -8.30991742950535\textrm{e-}5$  & $ 7.55399650329534\textrm{e-}6$  & $ 9.53674314457072\textrm{e-}7$ \\
\end{tabular}
}
\caption{Coefficients of $\bar{p}_{ij}$}
\label{TAB: f4}
\end{table}

Finally, we consider a more general linear system of ODEs
\begin{equation}
 \oder{\ubm}{t} = \Acal\ubm,
\end{equation}
where $\Acal = \Lcal + \Dcal + \Ucal$ is an $n \times n$ matrix, $\Lcal$ is
the lower triangular part of $\Acal$, $\Ucal$ is the upper triangular part
of $\Acal$, and $\Dcal$ is the diagonal of $\Acal$. In the remainder of this
section, we show that if $\Acal$ is diagonally dominant with negative diagonal
entries and the coupling term is taken as $\cbm(\ubm) = (\Lcal+\Ucal)\ubm$,
both the Jacobi and Gauss-Seidel predictors are unconditionally stable
for the forward-backward Euler IMEX scheme (\eqnref{EQ: IMEX-O1}). The
update matrix for the weak/strong Jacobi predictor takes the form
\begin{equation}
 \Cboldcal^J = (\Ical - \dt{} \Dcal)^{-1}(\Ical + \dt{} \Ucal + \dt{} \Lcal)
\end{equation}
and its spectral radius is
\begin{equation}
 \rho(\Cboldcal^J) \leq
 ||\Cboldcal^J||_\infty =
 \max_i \frac{\sum_{j\neq i} |\dt{} a_{i,j}|+1}{1-\dt{} a_{i,i}} \leq 1,
\end{equation}
where the first inequality follows from the Gershgorin circle theorem
and the last uses the assumption of diagonal dominance and negative
diagonal entries. This confirms that, under the stated assumptions,
the weak/strong Jacobi predictor with the IMEX1 scheme is unconditionally
stable. The update matrix for the weak/strong Gauss-Seidel predictor takes
the form
\begin{equation}
 \Cboldcal^{GS} = (\Ical - \dt{}\Lcal - \dt{}\Dcal)^{-1}
                  (\Ical + \dt{}\Ucal).
\end{equation}
Any of its eigenpairs $(\lambda, \xbm)$ satisfy the relation
\begin{equation}
 (\Ical - \dt{}\Lcal - \dt{}\Dcal)^{-1}
                   (\Ical + \dt{}\Ucal) \xbm =
 \lambda \xbm,
\end{equation}
which can be re-arranged as
\begin{equation}
 (\Ical + \dt{}\Ucal) \xbm = (\Ical - \dt{}\Lcal - \dt{}\Dcal) \lambda \xbm
\end{equation}
or written as components as
\begin{equation}
 x_i + \dt{}\sum_{j>i} a_{ij}x_j + \dt{}\lambda \sum_{j<i} a_{ij} x_j =
 \lambda x_i - \dt{}\lambda a_{ii}x_i
\end{equation}
for $i = 1,\dots,N$. Application of the triangular inequality and division by
$|x_i|$ leads to the relation
\begin{equation} \label{eqn:app:eq0}
 1 + \dt{}\sum_{j>i} |a_{ij}| \frac{|x_j|}{|x_i|} +
 \dt{}|\lambda| \sum_{j<i} |a_{ij}| \frac{|x_j|}{|x_i|} =
 |\lambda| |1 - \dt{} a_{ii}|.
\end{equation}
The assumption of diagonal dominance and negative diagonal entries leads to
the following bound
\begin{equation} \label{eqn:app:ineq1}
 |\lambda| |1 - \dt{} a_{ii}| =
 |\lambda| (1 + \dt{} |a_{ii}|) \geq
 |\lambda| (1 + \dt{} \sum_{j\neq i} |a_{ij}|).
\end{equation}
On the other hand, if $i = \arg\max_{1\leq j\leq n} |x_j|$, (\ref{eqn:app:eq0})
leads to
\begin{equation} \label{eqn:app:ineq2}
 |\lambda| |1 - \dt{} a_{ii}| \leq
 1 + \dt{}\sum_{j>i} |a_{ij}| + \dt{}|\lambda| \sum_{j<i} |a_{ij}|.
\end{equation}
Combining (\ref{eqn:app:ineq1}) and (\ref{eqn:app:ineq2}), we arrive at
\begin{equation}
 |\lambda|(1 + \dt{} \sum_{j>i} |a_{ij}|) \leq 1 + \dt{} \sum_{j>i}|a_{ij}|
\end{equation}
for $i = \arg\max_{1\leq j\leq n} |x_j|$, which leads to the desired
result
\begin{equation}
 \rho(\Cboldcal^{GS}) \leq 1
\end{equation}
and confirms that, under the stated assumptions, the weak/strong Gauss-Seidel
predictor with the IMEX1 scheme is unconditionally stable.

\bibliographystyle{unsrt}
\bibliography{myref}

\end{document}